\documentclass[11pt,a4paper]{article}
\usepackage{geometry}

\UseRawInputEncoding
\usepackage{bbold}
\usepackage[utf8]{inputenc}
\usepackage[T1]{fontenc}
\usepackage{amsmath}
\usepackage{amsthm}
\usepackage{amssymb}
\usepackage{graphicx}
\usepackage{caption}
\usepackage{multirow}
\usepackage{mathrsfs}
\usepackage{tabularx}
\usepackage{dsfont}
\usepackage{hyperref}
\usepackage[dvipsnames]{xcolor}
\usepackage{verbatim}
\usepackage[shortlabels]{enumitem}

\definecolor{linkcol}{rgb}{0,0,0.4} 
\definecolor{citecol}{rgb}{0.5,0,0} 

\hypersetup{colorlinks=true, linkcolor=linkcol, urlcolor=linkcol, citecolor=citecol}

\usepackage[ruled,vlined]{algorithm2e}
\SetKwProg{Init}{Initialisation}{}{}
\usepackage{float}
 \usepackage{ulem}
\usepackage{array,multirow,makecell}
\setcellgapes{1pt}
\makegapedcells
\newcolumntype{R}[1]{>{\raggedleft\arraybackslash }b{#1}}
\newcolumntype{L}[1]{>{\raggedright\arraybackslash }b{#1}}
\newcolumntype{C}[1]{>{\centering\arraybackslash }b{#1}}
\makeatletter
\newcommand\justify{%
  \let\\\@centercr
  \rightskip\z@skip
  \leftskip\z@skip}
\makeatother

\newtheorem{thm}{Theorem}[section]

\newtheorem{prop}[thm]{Proposition}

\newtheorem{ass}{Assumption}[section]

\theoremstyle{definition}
\newtheorem{defi}{Definition}[section]

\newtheorem{rmq}{Remark}[section]

\newenvironment{preuve}{\begin{proof} \rm}{\end{proof}}

\newtheorem*{prop*}{Proposition}
\newtheorem*{rmq*}{Remark}

\numberwithin{equation}{section}

\makeatletter
\newcommand{\pushright}[1]{\ifmeasuring@#1\else\omit\hfill$\displaystyle#1$\fi\ignorespaces}
\newcommand{\pushleft}[1]{\ifmeasuring@#1\else\omit$\displaystyle#1$\hfill\fi\ignorespaces}
\makeatother
\newcommand{\p}{\mathbb{P}}
\newcommand{\R}{\mathbb{R}}

\newcommand{\e}{\mathbb{E}}

\newcommand{\n}{\mathbb{N}}

\newcommand{\f}{\mathbb{F}}

\newcommand{\U}{\mathcal{U}}
\newcommand{\Hr}{\mathcal{H}}

\newcommand{\Or}{\mathcal{O}}

\newcommand{\K}{\mathcal{K}}

\newcommand{\Z}{\mathcal{Z}}

\newcommand{\E}{\mathcal{E}}

\newcommand{\F}{\mathcal{F}}

\newcommand{\A}{\mathcal{A}}
\newcommand{\Lr}{\mathcal{L}}

\newcommand{\ABS}[1]{{\left| #1 \right|}} % |1|
\newcommand{\CRO}[1]{{\left< #1 \right>}} % <1>
\newcommand{\PAR}[1]{{\left(#1\right)}} % (1)
\newcommand{\SBRA}[1]{{\left[#1\right]}} % [1]
\newcommand{\BRA}[1]{{\left\{#1\right\}}} % {1}
\newcommand{\NRM}[1]{{\left\Vert #1\right\Vert}} % ||1||

\newcommand{\ACC}[1]{{\left\{
    \begin{array}{ll}
    #1
    \end{array}
\right.}}

\DeclareMathOperator{\Tr}{Tr}

\DeclareMathOperator{\drift}{drift}

\DeclareMathOperator{\loc}{loc}

\title{Optimal investment and consumption under forward utilities with relative performance concerns\thanks{Acknowledgements:  The first and the second authors research is part of the ANR project DREAMeS (ANR-21-CE46-0002) and benefited from the support of respectively the "Chair Risques Emergents en Assurance (RE2A)"  and  "Chair Impact de la Transition Climatique en Assurance (ITCA)"  under the aegis of Fondation du Risque, a joint initiative by  Risk and Insurance Institute of  Le Mans,  and MMA-Cov\'ea and Groupama respectively. }}  
\author{Guillaume Broux-Quemerais\textsuperscript{1} \and Anis Matoussi\textsuperscript{1} \and Chao Zhou\textsuperscript{2}}

\begin{document}

\maketitle
\rightskip=2cm
\leftskip=2cm

% Affiliations
\footnotetext[1]{\footnotesize Laboratoire Manceau de Math\'ematiques, membre de la Fédération de Recherche CNRS Henri Lebesgue Bretagne - Pays de la Loire FR 2962 du CNRS et partenaire du LabEx Centre Henri Lebesgue. Institut du Risque et de l'Assurance, Le Mans University}
\footnotetext[2]{\footnotesize Department of Mathematics and Risk Management Institute, National University of Singapore.}

\begin{small}
\begin{center}
\section*{Abstract}
\end{center}
We study an $n$-player and mean-field portfolio optimization problem under relative performance concerns with non-zero volatility for both wealth and consumption. Time consistency of forward relative performance criteria is ensured through martingale optimality along a relative performance metric. This yields a stochastic partial differential equation (SPDE) of Hamilton-Jacobi-Bellman (HJB) type and highlights a structural link between wealth and consumption utilities. In the constant relative risk aversion (CRRA) case, under a suitable condition on the drift of the consumption utility, we derive the best-response policy and prove its admissibility. Characterizing $n$-player Nash equilibria is more delicate due to the strong coupling across agents induced by the consistency HJB-SPDE. This motivates the analysis of the mean-field regime, where interactions arise through the conditional law of the state given the common noise. In this setting, we show that the equilibrium consumption process solves an McKean-Vlasov stochastic differential equation (SDE). Under CRRA specifications, we provide sufficient conditions ensuring admissibility of the resulting mean-field Nash equilibrium when the consumption dynamics are independent of the conditional law. Finally, we present numerical experiments for the equilibria and argue that non-zero volatility acts as an additional preference parameter that is central in inverse preference problems.
\end{small}

\rightskip=0cm
\leftskip=0cm

\vspace{1cm}
%\tableofcontents

\newpage
\section*{Introduction}

We study $n$-player and mean-field investment-consumption optimization problems with relative performance concern in an asset specialization setting with common noise, within the forward preference framework. In this setting, each agent chooses her portfolio-consumption strategy while having concerns about the average wealth and consumption of other investors. Asset specialization refers to a situation in which investment managers focus on specific assets, due, for example, to sector familiarity, trading friction, liquidity costs (\cite{wagner2011systemic}) and ambiguity aversion (\cite{mukerji2001ambiguity}). Competitive investment-consumption problems have been analyzed under expected utility maximization in, among others, \cite{lacker2019mean}, \cite{lacker2020many}, \cite{hu2022n}, \cite{fu2023mean}, \cite{becherer2024common}. Forward relative performance criteria have been introduced and studied in \cite{anthropelos2022competition}, \cite{dos2021forward} and \cite{dos2022forward}, but, to the best of our knowledge, only in the locally riskless (zero utility volatility) class.

Forward utilities have been introduced and developed in \cite{musiela2006investments} to address limitations of the classical expected utility, which is specified through a fixed horizon and yields strategies that depend strongly on pre-specified terminal preferences. Forward criteria enable the dynamic adjustment of decision rules, starting from preferences specified at an initial time and remaining time-consistent thereafter, rather than imposing a potentially distant and arbitrary time horizon. Preferences are described by a utility random field ($U(t,\cdot)$), and time consistency within the given decision-making context is expressed through the martingale optimality principle; if $X_t^{\pi, c}$ is the observable process (typically the wealth) resulting from the admissible  decision/strategy $(\pi, c)$, then the preference process $U(t, X_t^{\pi, c})$ is a supermartingale, and there exists an optimal strategy such that the preference process is a martingale. In the It\^o setting, \cite{musiela2010stochastic}, \cite{nicole2013exact} established that time consistency is enforced through a stochastic HJB equation: the martingale optimality property translates into an SPDE constraint on the utility random field $U$. A key feature of the forward approach is that the utility volatility appears as an agent-specific input, while the drift is constrained by the HJB-SPDE. This flexibility implies an inherent non-uniqueness when only the initial specification is given, and characterizing volatility inputs that yield well-behaved solutions remains delicate. In other words, utility volatility acts as an additional preference parameter, beyond the usual curvature (or risk aversion for CRRA utilities), controlling how sensitively the agent's criterion reacts to the realized environment. In particular, restricting to the zero-volatility class as in \cite{dos2022forward}, \cite{zariphop24PUQR}, produces a tractable family in which randomness enters only through the market input. Allowing for a non-vanishing utility volatility yields genuinely stochastic preference criteria, with richer pathwise responsiveness to the model's randomness, see \cite{ng2024optimal}. Forward criteria for investment and consumption have been developed in \cite{el2018consistent} (see also \cite{kallblad2020black}), where dynamics of wealth and consumption utilities are linked by the consistency HJB-SPDE. Forward utilities have found diverse applications over recent years, including but not limited to option valuation, insurance, mean-field games (\cite{dos2021forward}, \cite{zariphop24PUQR}), long term interest rate modeling (\cite{el2022ramsey}), risk measures (\cite{chong2019pricing}), or more recently pension design (\cite{hillairet2022time},\cite{ng2024optimal}).

In competitive portfolio problems, utility functions include a relative performance metric, to evaluate the impact of a strategy relative to others. Early contributions under classical expected utility include \cite{espinosa2015optimal}, which studies Nash equilibrium in a Black-Scholes model for exponential preferences. Asset specialization models under CARA and CRRA preferences have been developed in \cite{lacker2019mean}, and extended to include consumption optimization in \cite{lacker2020many}, see also \cite{hu2022n} for more general It\^o-diffusion environment. In the mean-field specialization problem with random coefficients, \cite{fu2023mean} characterizes equilibria via a  quadratic growth FBSDE, and \cite{becherer2024common} studies competitive investment with common noise allowing both Brownian and jump sources. Optimal investment consumption problems with relative concerns have also been studied in \cite{dianetti2024optimal}, under the framework of recursive utilities of type Epstein-Zin. Within the forward preference framework, incorporating a relative performance metric leads to forward relative performance criteria, initiated in \cite{anthropelos2022competition} for two investors with homothetic, locally riskless forward utilities in general It\^o-diffusion markets. This line of work has been extended to $n$-player and mean-field formulations in \cite{zariphop24PUQR}, and related investment-consumption problems under forward relative criteria are studied in \cite{dos2021forward} and \cite{dos2022forward}. The discrete-time, predictable version of forward preferences introduced in \cite{angoshtari2020predictable} has also been extended to relative performance evaluation in a binomial setting in \cite{liang2023predictable}. To the best of our knowledge, all existing analysis of forward relative performance criteria are restricted to the locally riskless (null utility volatility) class.

Our objective is to analyze competitive investment-consumption under forward relative performance criteria with non-zero volatility. We consider $n$ agents investing in a financial market, whose randomness arises from $n+1$ independent Brownian motions $(B_{t})_{t \geq 0}$ and $\PAR{(W_{t}^{i})_{t \geq 0}}_{1 \leq i \leq n}$. Each agent trades exclusively in a private asset, whose price is driven by her idiosyncratic noise $W^{i}$ and the noise common to the entire population $B$. We introduce relative wealth and consumption rate metrics $(\widehat{x}^{i}, \, \widehat{c}^{i})$ and consider, for agent $i$, preference criteria of the form
\begin{eqnarray*}
    Q^{i}(t, \widehat{x}^{i} ; (\pi_{s}^{-i}, c_{s}^{-i})_{s \leq t}) = U^{i}(t, \widehat{x}^{i}) + \int_{0}^{t} V^{i}(s, \widehat{c_{s}}^{i} \widehat{x}^{i}) \, ds,
\end{eqnarray*}
where $U^{i}$ and $V^{i}$ are utility random fields from wealth and consumption, and $(\pi_{s}^{-i}, c_{s}^{-i})_{s \leq t}$ denotes the strategies of all agents except agent \( i \) up to time $t$. Assuming the wealth utility to be an It\^o random field, we characterize time consistency of the forward performance criterion $Q^{i}$ through an HJB-SPDE on $(U^{i}, V^{i})$, and derive the associated best-response policy, that is assuming competitors strategies admissible and given. Then, we focus on separable utilities of CRRA type, of the form
\begin{eqnarray*}
U^{(\alpha_{i})}(t,x) = Z_{t}^{i}u^{(\alpha_{i})}(x), \quad V^{(\alpha_{i})}(t,c) = \phi_{t}^{i} u^{(\alpha_{i})}(c), \quad i=1, ..., n,
\end{eqnarray*}
where $u^{(\alpha_{i})}$ denotes the standard power utility function with risk aversion $\alpha_{i} \in (0, 1) \cup (1, +\infty)$, and $(Z_{t}^{i})_{t \geq 0}$ and $(\phi_{t}^{i})_{t \geq 0}$ are continuous stochastic processes with log-normal dynamics. In this case, the consistency SPDE reduces to a drift restriction for $Z^{i}$. Under additional structural assumptions on $\phi^{i}$, we derive an explicit best-response strategy. The optimal portfolio depends on the volatility characteristics of the wealth utility factor $Z^{i}$, while the optimal consumption rate follows a logistic like dynamics for which we derive an explicit representation. Deriving an admissible Nash equilibrium for the $n$-player optimization problem is much more challenging since the competitors' strategies are no longer fixed, and the equilibrium is characterized by a system of coupled consumption processes. We identify the candidate Nash equilibrium associated with CRRA specifications, and establish its admissibility only in a restricted setting.

This limitation naturally motivates the mean-field formulation, in which the interaction among competitors arises though the conditional law of the state given the common noise. Indeed, the study of optimization problems with a large number of interacting agents naturally leads to the situation of an infinite number of players, formalized through the mean-field game theory introduced in \cite{huang2006large}, \cite{lasry2007mean} and generalized to the common noise setting in \cite{carmona2016}, \cite{carmona2018probabilistic}. The mean-field game formulation relies on a random type vector representation introduced in \cite{lacker2019mean}, representing the distribution of preferences and model parameters of a continuum of agents. We investigate the example of CRRA separable utilities, by considering two utility random fields from wealth and consumption $U$ and $V$, respectively driven by stochastic processes $Z$ and $\phi$. In this setting, the random type vector includes the stock price parameters, as well as the risk aversion and competition parameters. In the forward framework, an additional dynamic preference factor captures the evolution, across the population, of the distribution of time-dependent random utility parameters, namely the volatility local characteristics of $Z$ and $\phi$, while the drift of $Z$ is constrained by a mean-field HJB-SPDE that ensures time consistency. The mean-field forward performance criterion $Q$ generated by $U$ and $V$, is evaluated along relative wealth and consumption processes $\widehat{X_{t}} = \frac{X_{t}}{\PAR{\overline{X_{t}}}^{\theta}}, \quad \widehat{c_{t}} = \frac{c_{t}}{\PAR{\overline{C_{t}}}^{\theta}}$
where $\overline{X_{t}}$ and $\overline{C_{t}}$ are geometric average wealth and consumption process of the continuum of agents. The strategy for solving the mean-field forward optimization problem is the following. First fix stochastic processes $(\overline{X_{t}})_{t \geq 0}$ and $(\overline{C_{t}})_{t \geq 0}$, and solve the martingale optimality problem on $Q(t, \widehat{X_{t}})$ using the metric relative to these fixed average processes. Then given this optimal control, determine the law of $\overline{X_{t}^{*}}$ conditionally on the common noise $B$. Since $\overline{X_{t}}$ and $\overline{C_{t}}$ model the geometric average wealth  and consumption of the continuum of agents, determining the equilibrium amounts to find a fixed point $\overline{X_{t}}, \, \overline{C_{t}}$ satisfying 
\begin{eqnarray*}
    \overline{X_{t}} = \exp \PAR{\e \SBRA{\log X_{t}^{*} | \F_{t}^{B}}}, \quad  \overline{C_{t}} = \exp \PAR{\e \SBRA{\log c_{t}^{*} | \F_{t}^{B}}}.
\end{eqnarray*}
This compatibility condition reflects that the optimal strategy conditionally on the common noise must be typical of the population. The case of separable utilities of CRRA type is tractable; we show that under structural assumptions on processes $Z$ and $\phi$ ensuring consistency of the relative performance criterion, the MF Nash equilibrium consumption is characterized by a McKean-Vlasov SDE with common noise. Typical results of existence and uniqueness of a solution to such equations, see \cite{kumar2022well}, ensure that the resulting optimal consumption is well-defined and admissible, provided that $\phi$ exhibits a mean reversion term with respect to optimal consumption as does $Z$, with one weight common across the population. This assumption reduces the characterization of the equilibrium consumption to a standard SDE, for which we derive the explicit solution. 

The main novelty of our approach is twofold. First, the optimal portfolio strategy explicitly depends on wealth utility volatility. Second, the optimal consumption rate process is stochastic, characterized by the relative importance the agents assign to wealth utility compared to consumption utility. This contrasts with deterministic equilibria obtained in the locally riskless setting in \cite{dos2022forward}, \cite{zariphop24PUQR}, and shows that optimal strategies associated with forward utilities with non zero volatility include preference randomness. We also highlight the role of non-zero utility volatility in the inverse preference problem: given a target investment portfolio, the wealth utility volatility can be used as a calibration parameter to match the desired behavior. Our results then provide the corresponding equilibrium consumption and ensure that the resulting preference criterion remains dynamically consistent.

The paper is organized as follows. In Section \ref{sect:FRPP}, we introduce best-response forward relative performance criteria and establish the HJB-SPDE ensuring time consistency. We derive an explicit equilibrium for separable utilities of CRRA type. In Section \ref{sect:nagent}, we investigate the existence of a Nash equilibrium in the $n$-player situation. The mean-field optimization problem against a continuum of agent is studied in Section \ref{sect:MF}. Finally, we illustrate our results with numerical examples in Section \ref{sect:discMF}. 

\newpage
\noindent\textbf{Notations:}\\
All stochastic processes in the sequel are defined on a standard probability space $\PAR{\Omega, \F, \f, \p}$ supporting $n+1$ independent Brownian motions $B$, $W^{1}, ..., W^{n}$, where $(\F_{t})_{t \geq 0}$ is the natural augmented filtration they generate.  For $x, y \in \R^{d}$, we denote by $x^{\top}$ the transpose of vector $x$, $x.y = x y^{\top}$ the scalar product between $x$ and $y$, and $\NRM{.}$ the usual norm $\NRM{x} = \Tr(x x^{\top})^{\frac{1}{2}}$. For $p \in \n^{*}$ and $\Or \subset \R^{d}$, we define the following spaces
\begin{eqnarray*}
\Hr^{p}_{\loc}(\Or) &=& \BRA{\varphi \, -\text{prog. measurable process valued in} \, \Or \, \text{ s.t. for all } T > 0, \, \e \SBRA{\int_{0}^{T} \NRM{\varphi_{s}}^{p} ds} < \infty}. \\
L^{\infty}_{\loc}(\Or) &=& \BRA{\varphi \, -\text{prog. measurable process valued in} \, \Or \, \text{ s.t. for all } T > 0, \, \underset{0 \leq s \leq T}{\sup} \ABS{\varphi_{s}} < \infty, \ a.s}
\end{eqnarray*}
Let $\U_{std}$ denote the set of standard deterministic utility functions. We next introduce spaces for studying the regularity of semimartingales in terms of their local characteristics, following \cite{kunita1997stochastic}, \cite{ikeda2014stochastic}.
\paragraph{Definition of seminorms -} Let $\beta$ be an $\R^{k}$-valued random field of class \hfill \linebreak
$C^{m, \delta}\PAR{]0, +\infty[}$, with $m$ a nonnegative integer and $\delta$ a number in $(0, 1]$, i.e. $\beta$ is $m$ times differentiable in $x$ and its $m^{\text{th}}$ derivative is $\delta$-Hölder, for any $t$, almost surely. We introduce the following family of random Hölder $K$-seminorms aiming to control asymptotic behavior of $\beta$ and the regularity of its Hölder derivatives, for any $K \subset ]0, +\infty[$
\begin{eqnarray*}
\NRM{\beta}_{m, K}(t, \omega) &=& \underset{x \in K}{\sup} \frac{\NRM{\beta(t, x, \omega)}}{x} + \sum_{1 \leq j \leq m} \underset{x \in K}{\sup} \NRM{\partial_{x}^{j} \beta(t, x, \omega)} \\
\NRM{\beta}_{m, \delta, K}(t, \omega) &=& \NRM{\beta}_{m, K}(t, \omega) + \underset{x, y \in K}{\sup} \frac{\NRM{\partial_{x}^{m} \beta(t, x, \omega) - \partial_{x}^{m} \beta(t, y, \omega)}}{\ABS{x-y}^{\delta}}.
\end{eqnarray*}
As mentioned in \cite{nicole2013exact}, the first term of these random semi-norms is divided by $x$, in order to control behavior in the neighborhood of $x=0$, by imposing at most linear vanishing at the boundary. This normalization also preserves traditional asymptotic results from \cite{kunita1997stochastic}.

\paragraph{Associated function spaces -} The previous norms are related to the space parameter. We add the temporal dimension by requiring these seminorms (or their square) to be integrable in time with respect to Lebesgue measure on $[0, T]$. We then define the following sets:
\begin{enumerate}
    \item $\K^{m}_{\loc}$ (resp. $\overline{\K}^{m}_{\loc}$) denotes the set of $C^{m}$-random fields $\beta$ such that $\frac{\beta}{x}$ and $\partial_{x}^{k} \beta$ for $k \leq m$ are $\Lr^{1}$ (resp. $\Lr^{2}$)-locally bounded, that is for any compact $K \subset ]0, +\infty[$ and any $T$, 
    \begin{eqnarray*}
    \int_{0}^{T} \NRM{\beta}_{m, K}(t, \omega)dt < \infty, \quad \PAR{\text{resp. } \int_{0}^{T} \NRM{\beta}_{m, K}^{2}(t, \omega)dt < \infty}.
    \end{eqnarray*}
    \item $\K^{m, \delta}_{\loc}$ (resp. $\overline{\K}^{m, \delta}_{\loc}$) denotes the set of $C^{m, \delta}$-random fields such that for any compact $K \subset ]0, +\infty[$ and any $T$, 
    \begin{eqnarray*}
    \int_{0}^{T} \NRM{\beta}_{m, \delta, K}(t, \omega)dt < \infty, \quad \PAR{\text{resp. } \int_{0}^{T} \NRM{\beta}_{m, \delta, K}^{2}(t, \omega)dt < \infty}.
    \end{eqnarray*}
    \item When these norms are defined on the whole space $]0, +\infty[$, the derivatives up to a certain order are bounded in the spatial parameter, with an integrable (resp. square integrable) random bound, so that we use the notation $\K_{b}^{m}, \, \overline{\K}_{b}^{m}$ or $\K_{b}^{\delta, m}, \, \overline{\K}_{b}^{m, \delta}$.
\end{enumerate}

\newpage

\section{Forward relative performance criteria in the $n$-player setting} \label{sect:FRPP}

We investigate investment-consumption optimization in an asset specialization setting with a finite population of agents. Each manager trades exclusively in a private asset, and stock prices are driven by an idiosyncratic noise and a noise common to the entire population. Preferences are modeled with utility random fields evaluated along relative wealth and consumption metrics, to reflect their will to compete against each other. We first introduce the relative performance metrics and the resulting dynamics of relative wealth. We then derive the consistency SPDE and the associated best response strategies in a general It\^o setting. Finally, we specialize to CRRA-type forward relative performance criteria and discuss the admissibility of optimal consumption.

\subsection{Asset specialization and relative performance}

Consider a market consisting in one riskless bond and $n$ risky assets. In the asset specialization model, stock $i$ is traded exclusively by agent $i$. Assuming the bond to be the numeraire, the discounted price $(S_{t}^{i})$ of stock $i$ solves
\begin{eqnarray}
\frac{dS_{t}^{i}}{S_{t}^{i}} = \mu_{i}dt + \nu_{i}dW_{t}^{i} + \sigma_{i}dB_{t}, \quad S_{0}^{i} = s_{0}^{i} > 0,
\end{eqnarray}
where $\mu_{i} \in \R$, and $\sigma_{i}, \, \nu_{i} \geq 0$ satisfy $\sigma_{i}^{2} + \nu_{i}^{2} > 0$. The Brownian motion $B$ is called the common noise since it induces a correlation across stocks, while $W_{t}^{i}$ are idiosyncratic noises, specific to each agent $i$. We assume $B, \, W^{1}, \dots, W^{n}$ are independent.

\paragraph{Agent's wealth -} For $i=1,...,n$, agent $i$ uses a self-financing strategy $(\pi_{t}^{i})_{t \geq 0}$ representing the fraction of wealth invested in stock $i$. Let $(c_{t}^{i})_{t \geq 0}$ be the consumption rate per unit of wealth and let the initial wealth of agent $i$ satisfy $X_{0}^{i} = x^{i} > 0$. The $i^{th}$ agent's wealth dynamics writes as
\begin{eqnarray} \label{eq:wealth}
dX_{t}^{i} = \pi_{t}^{i}X_{t}^{i}\PAR{\mu_{i}dt + \nu_{i}dW_{t}^{i} + \sigma_{i}dB_{t}} - c_{t}^{i}X_{t}^{i}dt, \quad X_{0}^{i} = x^{i}.
\end{eqnarray}
An investment-consumption strategy $(\pi^{i}, c^{i})$ is admissible if it is progressively measurable, valued in $\R \times (0, \infty)$, such that $\e \SBRA{\int_{0}^{t} (\ABS{\pi_{s}^{i}}^{2} + \ABS{c_{s}^{i}})ds} < \infty$ for all $t > 0$. We denote $\A$ the set of admissible strategies.

\begin{rmq}
The specialization model above coincides with the standard lognormal framework used in competitive portfolio problems (see \cite{lacker2019mean}, \cite{lacker2020many}, \cite{dos2022forward}). It is, in general, incomplete for each agent, since stock $i$ is driven by the two-dimensional noise $(W^{i}, B)$, while only one risky asset is traded. For later use, set $\Sigma_{i} = (\nu_{i}, \, \sigma_{i})$ and $\overline{W_{t}^{i}}= (W_{t}^{i} \, , B_{t})$. Then
\begin{eqnarray*}
    \frac{d S_{t}^{i}}{S_{t}^{i}} = \Sigma_{i}. \PAR{d \overline{W_{t}^{i}} + \eta_{i} dt}, \quad \eta_{i} = \frac{\mu_{i}}{\sigma_{i}^{2} + \nu_{i}^{2}} \Sigma_{i},
\end{eqnarray*}
and the wealth dynamics become $dX_{t}^{i} = X_{t}^{i} \SBRA{\pi_{t}^{i} \Sigma_{i}.(d\overline{W_{t}^{i}} + \eta_{i}dt) - c_{t}^{i}dt}$. Portfolio constraints can be incorporated by restricting the portfolio rescaled by the volatility $\pi_{t}^{i} \Sigma_{i}$ to a progressive constraint set, in which case, the risk premium $\eta_{i}$ enters only through its projection onto the constraint set, see \cite{el2018consistent}.
\end{rmq}

\paragraph{Agent's interaction and relative performance -} Each agent measures her performance in comparison to her competitors. We introduce the relative wealth and consumption rate metrics as
\begin{eqnarray}
\widehat{X}^{i} &=& \frac{X^{i}}{\PAR{\overset{\sim}{X}^{(-i)}}^{\theta_{i}}}, \quad \text{where} \quad \overset{\sim}{X}^{(-i)} = \PAR{\prod_{k \neq i}^{n} X^{k}}^{\frac{1}{n-1}}, \label{relativewealth} \\
\widehat{c}^{i} &=& \frac{c^{i}}{\PAR{\overset{\sim}{c}^{(-i)}}^{\theta_{i}}}, \quad \text{where} \quad \overset{\sim}{c}^{(-i)} = \PAR{\prod_{k \neq i}^{n} c^{k}}^{\frac{1}{n-1}}, \label{relativecons}
\end{eqnarray}
where $\theta_{i} \in \, \SBRA{0, 1}$ is the relative concern parameter. The closer $\theta_{i}$ is to one, the more agent $i$ cares about the geometric average wealth and consumption of her competitors. Applying It\^o's formula yields the dynamics of the geometric average wealth of other agents
\begin{multline}
\frac{d \overset{\sim}{X_{t}}^{(-i)}}{\overset{\sim}{X_{t}}^{(-i)}} = \PAR{\overline{\mu \pi_{t}}^{(-i)} - \frac{1}{2} \PAR{\overline{\Sigma \pi_{t}^{2}}^{(-i)} - \PAR{\overline{\sigma \pi_{t}}^{(-i)}}^{2} - \frac{1}{n-1} \overline{\PAR{\nu \pi_{t}}^{2}}^{(-i)}} - \overline{c_{t}}^{(-i)}}dt \\
+ \frac{1}{n-1} \sum_{k \neq i}^{n} \nu_{k} \pi_{t}^{k} dW_{t}^{k} + \overline{\sigma \pi_{t}}^{(-i)}dB_{t}, 
\end{multline}
where the bar notation denotes empirical averages over $k \neq i$, namely
\begin{eqnarray*}
\overline{\mu \pi_{t}}^{(-i)} &=& \frac{1}{n-1} \sum_{k \neq i}^{n} \mu_{k} \pi_{t}^{k}, \quad \overline{(\nu \pi_{t})^{2}}^{(-i)} = \frac{1}{n-1} \sum_{k \neq i}^{n} (\nu_{k} \pi_{t}^{k})^{2}, \quad \overline{\sigma \pi_{t}}^{(-i)} = \frac{1}{n-1} \sum_{k \neq i}^{n} \sigma_{k} \pi_{t}^{k}, \\
\overline{\Sigma \pi_{t}^{2}}^{(-i)} &=& \frac{1}{n-1} \sum_{k \neq i}^{n} \Sigma_{k}(\pi_{t}^{k})^{2}, \quad \overline{c_{t}}^{(-i)} = \frac{1}{n-1} \sum_{k \neq i}^{n} c_{t}^{k}, \quad \Sigma_{k} = \sigma_{k}^{2} + \nu_{k}^{2}.
\end{eqnarray*}

\noindent
Consequently, the relative wealth process satisfies
\begin{eqnarray} \label{xhatdyn}
\hspace{-0.3cm} \frac{d \widehat{X_{t}}^{i}}{\widehat{X_{t}}^{i}} = \xi_{i}dt - \PAR{c_{t}^{i} - \theta_{i} \overline{c_{t}}^{(-i)}}dt + \PAR{\nu_{i} \pi_{t}^{i} dW_{t}^{i} - \theta_{i} \frac{1}{n-1} \sum_{k \neq i}^{n} \nu_{k} \pi_{t}^{k} dW_{t}^{k}} + \PAR{\sigma_{i} \pi_{t}^{i} - \theta_{i} \overline{\sigma \pi_{t}}^{(-i)}}dB_{t},
\end{eqnarray}
where
\begin{eqnarray}
    \xi_{i} = \mu_{i} \pi_{t}^{i} - \theta_{i} \overline{\mu \pi_{t}}^{(-i)} + \frac{\theta_{i}}{2} \overline{\Sigma \pi_{t}^{2}}^{(-i)}  - \frac{\theta_{i}^{2}}{2} \PAR{(\overline{\sigma \pi_{t}}^{(-i)})^{2} + \frac{1}{n-1} \overline{(\nu \pi_{t})^{2}}^{(-i)}} - \theta_{i} \sigma_{i} \pi_{t}^{i} \overline{\sigma \pi_{t}}^{(-i)}.
\end{eqnarray}

\subsection{Forward relative performance criteria in the It\^o framework} \label{sect:bestresponsegen}

Each manager $i = 1,..., n$ evaluates her performance relative to her peers using a relative performance criteria, modeled as an $\F_{t}$-progressively measurable random field $Q^{i}: \Omega \times (0, \infty) \times [0, \infty) \to \R$. This preference captures respectively the utility from wealth and consumption through utility random fields $U(t, x)$ and $V(t, c)$ on $(\R^{+})^{2}$. Formally, a utility random field $U = (t,x,\omega) \in \R^+ \times \R^{+} \times \Omega \to \R$, is a collection of random utility functions such that:
\begin{itemize}
\item for all $t \geq 0$, and $x \in \R^{+}$, $U (t, x)$ is $\mathcal F_t$-measurable.
\item the functions $x \mapsto U(t,x,\omega)$ are non-negative, strictly concave increasing functions of class $\mathcal C^{2}$ on $]0, \infty[$, $(\omega,t)$ a.s.
\item $u_0:= U(0, \cdot)$ is a standard (deterministic) utility function.
\end{itemize}

\begin{defi}\label{FRPP}
Consider an agent $i \in \BRA{1, ..., n}$ and assume that each manager $j \neq i$ follows an admissible strategy $(\pi^{j}, c^{j})$. For $(\pi^{i}, c^{i}) \in \A$, define
\begin{eqnarray} \label{exprQ}
    Q^{i}(t, \widehat{x}^{i} ; (\pi_{s}^{-i}, c_{s}^{-i})_{s \leq t}) = U^{i}(t, \widehat{x}^{i}) + \int_{0}^{t} V^{i}(s, \widehat{c_{s}}^{i} \widehat{x}^{i}) \, ds,
\end{eqnarray}
where $U^{i}$ and $V^{i}$ are utility random fields from wealth and consumption respectively, and \( \pi_{t}^{-i} = (\pi_{t}^{j})_{j \neq i} \) and $c_{t}^{-i} = (c_{t}^{j})_{j \neq i}$ denote the strategies of all agents except agent \( i \). \\
The random field $Q^{i}$ is a best-response forward relative performance criterion for manager $i$ if, for all $t \geq 0$:
\begin{itemize}
\item \textit{Martingale optimality -} For any admissible $(\pi^{i}, c^{i})$, $Q^{i}(t, \widehat{X_{t}}^{i})$ is a local supermartingale, and there exists an admissible strategy $(\pi^{i, *}, c^{i, *})$ such that $Q^{i}(t, \widehat{X_{t}}^{i, *})$ is a local martingale.
\end{itemize}
\end{defi}
In the literature, such preference criteria are often referred to as best-response criteria (see, for instance, \cite{anthropelos2022competition}, \cite{zariphop24PUQR}), since they characterize manager \( i \)'s optimal response given the others' strategies. In what follows, we omit the dependence on the trajectory of the other agents' strategies $(\pi_{s}^{-i}, c_{s}^{-i})_{s \leq t})$ when no confusion can arise. In addition, the martingale optimality condition is the standard dynamic programming characterization of market and time consistency. Hence, corresponding preferences are commonly referred to as consistent, see for instance \cite{el2008backward}, \cite{el2018consistent}.

As extensively discussed in \cite{anthropelos2022competition}, \cite{zariphop24PUQR}, two features distinguish the forward-preference framework from classical expected-utility theory. First, at time $t$, an agent’s best-response problem only requires knowledge of her competitors’ policies up to time $t$; future competitor strategies need not be fixed in advance. In contrast, solving the standard utility-maximization problem on a fixed horizon $\SBRA{0,T}$ typically requires specifying competitors’ strategies over the whole interval. Second, the forward best-response criterion can be formulated without committing ex ante to the full future market dynamics. By comparison, the expected-utility problem on $\SBRA{0,T}$ requires such pre-specification in order to evaluate terminal wealth. In our setting, the forward criterion $Q^{i}$ preserves time consistency with the market environment via the martingale optimality principle. Consequently, the resulting forward best-response policies are adaptive and not tied to arbitrary prior choices of horizon, model, or competitors’ future strategies.

Assuming It\^o dynamics for the wealth utility random field, the martingale optimality induces a drift restriction on $Q^{i}(t, \widehat{X_{t}}^{i})$. Under sufficient regularity conditions, this leads a stochastic HJB-type SPDE characterizing time consistency of relative performance criteria. To analyze the criterion's drift variations, we recall the Fenchel-Legendre transform, which is well defined under the regularity assumptions on $V^{i}$.

\begin{defi}
Let $V: \Omega \times (0, \infty) \times [0, \infty) \to \R$ be a random field, strictly concave in $x$. The Fenchel-Legendre transform of $V$, denoted $\overset{\sim}{V}$ is defined by
\begin{eqnarray}
\overset{\sim}{V}(t, x') = \underset{x > 0}{\sup} \BRA{V(t, x) - x'x}, \quad x' > 0, \, t \geq 0.
\end{eqnarray}
\end{defi}

\noindent
We work within the regular random field framework used to study differentiability properties of Itô random fields, see also \cite{nicole2013exact}. Regularity of a utility random field is strongly tied to the regularity of its local characteristics. Specifically, for $i = 1, \dots, n$, we assume that $U^i$ is an It\^o random field with local characteristics $\PAR{\beta^{i}, (\gamma^{i, W}, \gamma^{i, B})} \in \K^{2, \epsilon}_{\loc} \times \overline{\K}^{2, \epsilon}_{\loc}$ for some $\epsilon > 0$, that is,
\begin{eqnarray} \label{Itodynut}
    dU^{i}(t, x) = \beta^{i}(t, x)dt + \gamma^{i, W}(t, x)dW_{t}^{i} + \gamma^{i, B}(t, x)dB_{t}, \quad \text{a.s}.
\end{eqnarray}
By Theorem 2.2 in \cite{nicole2013exact}, $U^{i}$ is then an $\K^{2, \epsilon'}_{\loc}$ semimartingale, for any $\epsilon' < \epsilon$. We recall assumptions on the local characteristics ensuring $U^{i}(t, .)$ remains increasing and concave.
\begin{ass} \label{ass:regconc}
The random fields $(\beta^{i}, \gamma^{i})$ are $C^{2}$ in $x$. Moreover, there exists random bounds $D_{t}^{i} \in \Hr_{\loc}^{1}$ and $G_{t}^{i} \in \Hr_{\loc}^{2}$ such that a.s. for all $t \geq 0$, $x > 0$:
\begin{alignat}{2}
\ABS{\beta^{i}_{x}(t, x)} &\leq D_{t}^{i} \ABS{U^{i}_{x}(t, x)}, &\quad \NRM{\gamma^{i}_{x}(t, x)} &\leq G_{t}^{i} \ABS{U^{i}_{x}(t, x)}, \\
\ABS{\beta^{i}_{xx}(t, x)} &\leq D_{t}^{i} \ABS{U^{i}_{xx}(t, x)}, &\quad \NRM{\gamma^{i}_{xx}(t, x)} &\leq G_{t}^{i} \ABS{U^{i}_{xx}(t, x)}.
\end{alignat}
\end{ass}
\noindent
Under Assumption \ref{ass:regconc}, $U^{i}$ is a well-defined utility random fields, according to Corollary 1.3 in \cite{nicole2013exact}. Some regularity conditions are required in our verification Theorem \ref{SPDEprop}, to ensure that the optimal relative wealth process $\widehat{X_{t}}^{i, *}$ admits a unique strong global solution. The first is a growth condition on the risk tolerance ratio, which appears in the optimal portfolio policy.
\begin{ass} \label{ass:admport}
There exists an adapted random bound $R_{t}^{i} \in \Hr_{\loc}^{2}$ such that a.s, for all $t \geq 0$, $x > 0$:
\begin{eqnarray*}
    \ABS{\frac{U_{x}^{i}(t, x)}{U_{xx}^{i}(t, x)}} \leq R_{t}^{i} \ABS{x}.
\end{eqnarray*}
\end{ass}
\noindent
Second, we introduce a domination condition relating the marginal utility from consumption to the marginal utility from wealth. It yields a linear growth bound for the best-response consumption.
\begin{ass} \label{ass:admcons}
There exists an adapted random process $K_{t}^{i} \in \Hr_{\loc}^{1}$ such that, a.s, for all $t \geq 0$, $x > 0$:
\begin{eqnarray}
    V_{x}^{i}(t, K_{t}^{i} x (\overset{\sim}{c_{t}}^{(-i)})^{-\theta_{i}}) \geq U_{x}^{i}(t, x).
\end{eqnarray}
\end{ass}
\noindent
Note that it mirrors the assumption used in Theorem 4.8 in \cite{el2018consistent} for the single agent investment-consumption problem, with an additional scaling factor $(\overset{\sim}{c_{t}}^{(-i)})^{-\theta_{i}}$ due to competition. For ease of notation, we introduce the random field $\varphi$ defined as
\begin{small}
\begin{multline} \label{phidef}
\varphi(t, x) =  U_{x}^{i}(t, x)x \left( - \theta_{i} \overline{\mu \pi_{t}}^{(-i)}  + \frac{\theta_{i}}{2} \overline{\Sigma \pi_{t}^{2}}^{(-i)} + \frac{\theta_{i}^{2}}{2} \PAR{(\overline{\sigma \pi_{t}}^{(-i)})^{2} + \frac{1}{n-1} \overline{(\nu \pi_{t})^{2}}^{(-i)}} \right. \\
\left. + \theta_{i} \overline{c_{t}}^{(-i)} \right) + \frac{1}{2} U_{xx}^{i}(t, x) (x)^{2} \PAR{\frac{\theta_{i}^{2}}{n-1} \overline{(\nu \pi_{t})^{2}}^{(-i)} + (\theta_{i} \overline{\sigma \pi_{t}}^{(-i)})^{2}} - \gamma_{x}^{i, B}(t, x) \theta_{i} \overline{\sigma \pi_{t}}^{(-i)} x.
\end{multline}
\end{small}
\noindent
When convenient, we denote the optimal strategies by $\pi_{t}^{i, *}(x)$ and $c_{t}^{i, *}(x)$, to emphasize their dependence on the state variable $x$. Applying the It\^o-Ventzel's formula to the relative performance criterion $Q^{i}$ along the relative wealth process $\widehat{X}_{t}^{i}$ allows one to analyze the martingale optimality from Definition \ref{FRPP}. This yields a sufficient condition on the drift of the wealth utility $U^{i}$ ensuring time consistency of the criterion, as well as explicit best-response controls.

\begin{thm} \label{SPDEprop}
Consider an agent $i \in \BRA{1, ..., n}$ and assume that each manager $j \neq i$ follows an admissible strategy $(\pi^{j}, c^{j}) \in \A$. Let $(U^{i},V^{i})$ be a candidate random utility system with initial specification $(u_{0}^{i}, v_{0}^{i})$ for agent $i$, satisfying Assumptions \ref{ass:regconc}, \ref{ass:admport} and \ref{ass:admcons}. Suppose that the couple $(U^{i}, V^{i})$ solves the SPDE
\begin{multline} \label{utspde_first}
dU^{i}(t, x) = \PAR{-  \varphi(t, x) + \frac{1}{2} U_{xx}^{i}(t, x) x^{2} (\nu_{i}^{2} + \sigma_{i}^{2}) \PAR{\pi_{t}^{i, *}(x)}^{2} - \widetilde{V}^{i}(t, (\overset{\sim}{c_{t}}^{(-i)})^{\theta_{i}} U_{x}^{i}(t, x))}dt \\
+ \gamma^{i, W}(t, x)dW_{t}^{i} + \gamma^{i, B}(t, x)dB_{t},
\end{multline}
where 
\begin{multline} \label{opt_port}
    \pi_{t}^{i, *}(x) = \frac{1}{\nu_{i}^{2} + \sigma_{i}^{2}} \left( \sigma_{i} \theta_{i} \overline{\sigma \pi_{t}}^{(-i)} - \frac{1}{U_{xx}^{i}(t, x) x} \left(\gamma_{x}^{i, W}(t, x) \nu_{i} + \gamma_{x}^{i, B}(t, x) \sigma_{i} \right. \right. \\
    \left. \left. + (\mu_{i} - \theta_{i} \sigma_{i} \overline{\sigma \pi_{t}}^{(-i)})U_{x}^{i}(t, x) \right) \right)
\end{multline}
\begin{flalign}\label{optcons}
    c_{t}^{i, *}(x) = \frac{(V_{x}^{i})^{-1} \PAR{t, U_{x}^{i}(t, x)(\overset{\sim}{c_{t}}^{(-i)})^{\theta_{i}}} (\overset{\sim}{c_{t}}^{(-i)})^{\theta_{i}}}{x}.
\end{flalign}
Then $Q^{i}$ defined by \eqref{exprQ} is a best-response forward relative performance criterion for manager $i$, and the policy $(\pi^{i, *}, c^{i, *}) \in \A$ is optimal in the sense of Definition \ref{FRPP}.
\end{thm}

By Definition \ref{FRPP}, forward preferences are endogenous, self-generating. Theorem \ref{SPDEprop} is therefore understood as a verification result: starting from a candidate system $(U^{i}, V^{i})$ that satisfies the stochastic HJB equation \eqref{utspde_first} together with the feedback optimal controls \eqref{opt_port}--\eqref{optcons}, we obtain that the induced criterion is time consistent and that the associated best-response policy is optimal in the sense of Definition \ref{FRPP}. Assumptions \ref{ass:regconc}, \ref{ass:admport} and \ref{ass:admcons} are only used as sufficient conditions to ensure that the feedback controls are well-defined and admissible. In the next section, we construct explicit separable forward utilities of power type and verify they satisfy the above conditions, which shows that the class characterized by Theorem \ref{SPDEprop} is non-empty.

\begin{rmq}
When $\theta^{i} = 0$, the relative performance component disappears and we recover the standard SPDE characterization of forward investment-consumption criteria, together with the corresponding optimal controls, as in \cite{el2018consistent}. If the wealth utility volatility vanishes, $\gamma^{i}(t, x) = (0, 0)$, we recover the optimal portfolio and the same expression for the optimal consumption \eqref{optcons} as in \cite{dos2022forward}. A key difference is that, in our formulation, the utility from consumption $V^{i}$ is linked to the wealth utility field $U^{i}$ through the SPDE \eqref{utspde_first}, while \cite{dos2022forward} derives the analogous relation from a PDE with random coefficient. This distinction stems from the use of It\^o-Ventzel's formula instead of It\^o's to compute the dynamics of $Q^{i}(t, \widehat{X_{t}}^{i})$: the resulting drift involves the drift characteristics $\beta^{i}$ of the utility random field $U^{i}$ rather than its time derivative.
\end{rmq}

\begin{preuve}
The proof of this result is in three steps. First, we apply It\^o-Ventzel's formula to obtain the dynamics of $Q^{i}(t, \widehat{X_{t}}^{i})$. Then, we derive the optimal strategy processes using the first order condition on the drift of this compound process. Finally, we establish martingale optimality and state the condition on the drift of the wealth utility $\beta$ that ensures time consistency of the resulting preference criterion, which leads to the HJB-SPDE \eqref{utspde_first}.
 
\noindent
\textbf{Preference criterion's dynamics -} Applying It\^o-Ventzel's formula yields
\begin{align}
    dQ^{i}(t, \widehat{X_{t}}^{i}) & = \PAR{\beta^{i}(t, \widehat{X_{t}}^{i}) + V^{i}(t, \widehat{c_{t}}^{i} \widehat{X_{t}}^{i})}dt + \gamma^{i, W}(t, \widehat{X_{t}}^{i})dW_{t}^{i} + \gamma^{i, B}(t, \widehat{X_{t}}^{i})dB_{t} + U_{x}^{i}(t, \widehat{X_{t}}^{i})d\widehat{X_{t}}^{i} \nonumber \\
    &\quad + \frac{1}{2} U_{xx}^{i}(t, \widehat{X_{t}}^{i})d \CRO{\widehat{X_{t}}^{i}, \widehat{X_{t}}^{i}} + \CRO{\gamma_{x}^{i, W}(t, \widehat{X_{t}}^{i}) dW_{t}^{i} + \gamma_{x}^{i, B}(t, \widehat{X_{t}}^{i})dB_{t}, d\widehat{X_{t}}^{i}} \nonumber \\
    & = \PAR{\beta^{i}(t, \widehat{X_{t}}^{i}) + V^{i}(t, \widehat{c_{t}}^{i} \widehat{X_{t}}^{i}) + U_{x}^{i}(t, \widehat{X_{t}}^{i}) \widehat{X_{t}}^{i} \PAR{\xi_{i} - (c_{t}^{i} - \theta_{i} \overline{c_{t}}^{(-i)})}}dt + \gamma^{i, B}(t, \widehat{X_{t}}^{i})dB_{t} \nonumber \\
    &\quad + U_{x}^{i}(t, \widehat{X_{t}}^{i}) \widehat{X_{t}}^{i} \PAR{\PAR{\sigma_{t} \pi_{t}^{i} - \theta_{i} \overline{\sigma \pi_{t}}^{(-i)}} dB_{t} + \PAR{\nu_{i} \pi_{t}^{i} dW_{t}^{i} - \theta_{i} \frac{1}{n-1} \sum_{k \neq i}^{n} \nu_{k} \pi_{t}^{k} dW_{t}^{k}}} \nonumber \\
    &\quad + \frac{1}{2} U_{xx}^{i}(t, \widehat{X_{t}}^{i}) (\widehat{X_{t}}^{i})^{2} \PAR{(\nu_{i} \pi_{t}^{i})^{2} + \frac{\theta_{i}^{2}}{n-1} \overline{(\nu \pi_{t})^{2}}^{(-i)} + (\sigma_{i} \pi_{t}^{i} - \theta_{i} \overline{\sigma \pi_{t}}^{(-i)})^{2}}dt \label{eq:gammQIto} \\
    &\quad + \gamma_{x}^{i, W}(t, \widehat{X_{t}}^{i}) \nu_{i} \pi_{t}^{i} \widehat{X_{t}}^{i} dt + \gamma_{x}^{i, B}(t, \widehat{X_{t}}^{i}) \PAR{\sigma_{i} \pi_{t}^{i} - \theta_{i} \overline{\sigma \pi_{t}}^{(-i)}} \widehat{X_{t}}^{i} dt + \gamma^{i, W}(t, \widehat{X_{t}}^{i})dW_{t}^{i}. \nonumber
\end{align}
In particular, the drift of the above random field takes the form
\begin{align*}
\drift Q^{i}(t, \widehat{X_{t}}^{i}) & = \beta^{i}(t, \widehat{X_{t}}^{i}) + V^{i}(t, \widehat{c_{t}}^{i} \widehat{X_{t}}^{i}) + \gamma_{x}^{i, B}(t, \widehat{X_{t}}^{i}) \PAR{\sigma_{i} \pi_{t}^{i} - \theta_{i} \overline{\sigma \pi_{t}}^{(-i)}} \widehat{X_{t}}^{i}  \\
&\quad + \gamma_{x}^{i, W}(t, \widehat{X_{t}}^{i}) \nu_{i} \pi_{t}^{i} \widehat{X_{t}}^{i} + U_{x}^{i}(t, \widehat{X_{t}}^{i}) \left( \mu_{i} \pi_{t}^{i} - \theta_{i} \overline{\mu \pi_{t}}^{(-i)} + \frac{\theta_{i}}{2} \overline{\Sigma \pi_{t}^{2}}^{(-i)} \right. \\
&\quad \left. + \frac{\theta_{i}^{2}}{2} \PAR{(\overline{\sigma \pi_{t}}^{(-i)})^{2} + \frac{1}{n-1} \overline{(\nu \pi_{t})^{2}}^{(-i)}} - \theta_{i} \sigma_{i} \pi_{t}^{i} \overline{\sigma \pi_{t}}^{(-i)} - (c_{t}^{i} - \theta_{i} \overline{c_{t}}^{(-i)}) \right) \\
&\quad + \frac{1}{2} U_{xx}^{i}(t, \widehat{X_{t}}^{i}) (\widehat{X_{t}}^{i})^{2} \PAR{(\nu_{i} \pi_{t}^{i})^{2} + \frac{\theta_{i}^{2}}{n-1} \overline{(\nu \pi_{t})^{2}}^{(-i)} + (\sigma_{i} \pi_{t}^{i} - \theta_{i} \overline{\sigma \pi_{t}}^{(-i)})^{2}}. 
\end{align*}
\textbf{Optimal strategies -} The martingale optimality principle requires that, for any admissible control, the drift rate of $Q^{i}(t, \widehat{X_{t}}^{i})$ be non positive, and vanish at the optimum. This can be formulated as a first order condition on $\drift Q^{i}(t, \widehat{X_{t}}^{i})$. The resulting optimal portfolio is similar to that in \cite{dos2022forward}, with an additional term involving $\gamma_{x}^{i}$, while the optimal consumption process has the same form, namely
\begin{align*}
\pi_{t}^{i, *}(x) &= \frac{1}{\nu_{i}^{2} + \sigma_{i}^{2}} \left( \sigma_{i} \theta_{i} \overline{\sigma \pi_{t}}^{(-i)} - \frac{1}{U_{xx}^{i}(t, x) x} \left( \gamma_{x}^{i, W}(t, x) \nu_{i} + \gamma_{x}^{i, B}(t, x) \sigma_{i} \right. \right. \\
& \pushright{ \left. \left. +  (\mu_{i} - \theta_{i} \sigma_{i} \overline{\sigma \pi_{t}}^{(-i)})U_{x}^{i}(t, x) \right) \right)} \\
c_{t}^{i, *}(x) &= \frac{(V_{x}^{i})^{-1} \PAR{t, U_{x}^{i}(t, x)(\overset{\sim}{c_{t}}^{(-i)})^{\theta_{i}}} (\overset{\sim}{c_{t}}^{(-i)})^{\theta_{i}}}{x}.
\end{align*}
Assumptions \ref{ass:regconc}, \ref{ass:admport} and \ref{ass:admcons} are sufficient to control the singular behavior of these best-response policy near zero. Since the competitors' strategies $(\pi_{t}^{k}, c_{t}^{k})_{k \neq i}$ are assumed admissible, the integrability of the random bounds appearing in these assumptions implies that the optimal controls are well defined with $\pi_{t}^{i, *}(x) \in \Hr_{\loc}^{2}$ and $c_{t}^{i, *}(x) \in \Hr_{\loc}^{1}$. Hence $(\pi_{t}^{i, *}(x), c_{t}^{i,  *}(x)) \in \A$.

\noindent
\textbf{Martingale optimality -} It remains to verify that the induced criterion $Q^{i}$ satisfies the martingale optimality property of Definition \ref{FRPP}. Using the notation $\varphi$ introduced in \eqref{phidef}, the drift of $Q^{i}$ can be rewritten, for any admissible strategy $\pi_{t}^{i} \in \A$, as
\begin{align*}
\drift Q^{i}(t, \widehat{X_{t}}^{i}) & = \beta^{i}(t, \widehat{X_{t}}^{i}) + \varphi(t, \widehat{X_{t}}^{i}) + \widetilde{V}^{i}(t, (\overset{\sim}{c_{t}}^{(-i)})^{\theta_{i}} U_{x}^{i}(t, \widehat{X_{t}}^{i}))  \\
&\quad + \frac{1}{2} U_{xx}^{i}(t, \widehat{X_{t}}^{i}) (\widehat{X_{t}}^{i})^{2} \left( (\nu_{i}^{2} + \sigma_{i}^{2}) (\pi_{t}^{i})^{2} + 2(\mu_{i} - \theta_{i} \sigma_{i} \overline{\sigma \pi_{t}}^{(-i)}) \pi_{t}^{i} \frac{U_{x}^{i}(t, \widehat{X_{t}}^{i})}{U_{xx}^{i}(t, \widehat{X_{t}}^{i}) \widehat{X_{t}}^{i}} \right.\\ 
&\quad \left. - 2 \theta_{i} \sigma_{i} \overline{\sigma \pi_{t}}^{(-i)}. \pi_{t}^{i} + \frac{1}{U_{xx}^{i}(t, \widehat{X_{t}}^{i}) \widehat{X_{t}}^{i}} \PAR{2 \gamma_{x}^{i, W}(t, \widehat{X_{t}}^{i}) \nu_{i} \pi_{t}^{i} + 2 \gamma_{x}^{i, B}(t, \widehat{X_{t}}^{i}) \sigma_{i} \pi_{t}^{i}} \right) \\
&\quad + V^{i}(t, \widehat{c_{t}}^{i} \widehat{X_{t}}^{i}) - c_{t}^{i} U_{x}(t, \widehat{X_{t}}^{i}) \widehat{X_{t}}^{i} - \widetilde{V}^{i}(t, (\overset{\sim}{c_{t}}^{(-i)})^{\theta_{i}} U_{x}^{i}(t, \widehat{X_{t}}^{i}))  \\
& = \beta^{i}(t, \widehat{X_{t}}^{i}) + \varphi(t, \widehat{X_{t}}^{i}) + \widetilde{V}^{i}(t, (\overset{\sim}{c_{t}}^{(-i)})^{\theta_{i}} U_{x}^{i}(t, \widehat{X_{t}}^{i})) + V^{i}(t, \widehat{c_{t}}^{i} \widehat{X_{t}}^{i}) - c_{t}^{i} U_{x}(t, \widehat{X_{t}}^{i}) \widehat{X_{t}}^{i} \\
&\quad - \widetilde{V}^{i}(t, (\overset{\sim}{c_{t}}^{(-i)})^{\theta_{i}} U_{x}^{i}(t, \widehat{X_{t}}^{i})) + \frac{1}{2} U_{xx}^{i}(t, \widehat{X_{t}}^{i}) (\widehat{X_{t}}^{i})^{2} (\nu_{i}^{2} + \sigma_{i}^{2}) \PAR{  (\pi_{t}^{i})^{2} - 2 \pi_{t}^{i}. \pi_{t}^{i, *}(\widehat{X_{t}}^{i})}
\end{align*}
By the first order condition, the above drift is minimal for $\pi_{t}^{i} = \pi_{t}^{i, *}(\widehat{X_{t}}^{i})$. With a slight abuse of notation, let $\pi_{t}^{i, *} = \pi_{t}^{i, *}(\widehat{X_{t}}^{i, *})$ and $ c_{t}^{i, *} = c_{t}^{i, *}(\widehat{X_{t}}^{i, *})$. The corresponding minimum, evaluated along the optimal relative wealth process, is then equal to
\begin{multline} \label{cond_itocase}
\drift Q^{i}(t, \widehat{X_{t}}^{i, *}) = \beta^{i}(t, \widehat{X_{t}}^{i, *}) + \varphi(t, \widehat{X_{t}}^{i, *}) + \widetilde{V}^{i}(t, (\overset{\sim}{c_{t}}^{(-i)})^{\theta_{i}} U_{x}^{i}(t, \widehat{X_{t}}^{i, *})) + V^{i}(t, \widehat{c_{t}}^{i, *} \widehat{X_{t}}^{i, *}) \\
- c_{t}^{i, *} U_{x}(t, \widehat{X_{t}}^{i, *}) \widehat{X_{t}}^{i, *} - \widetilde{V}^{i}(t, (\overset{\sim}{c_{t}}^{(-i)})^{\theta_{i}} U_{x}^{i}(t, \widehat{X_{t}}^{i, *})) - \frac{1}{2} U_{xx}^{i}(t, \widehat{X_{t}}^{i, *}) (\widehat{X_{t}}^{i, *})^{2} (\nu_{i}^{2} + \sigma_{i}^{2}) \PAR{\pi_{t}^{i, *}}^{2}. 
\end{multline}
Moreover, by definition of the Fenchel-Legendre transform $\widetilde{V}^{i}$ and using \eqref{optcons}, we have for any $x > 0$,
\begin{eqnarray*}
\widetilde{V}^{i}(t, (\overset{\sim}{c_{t}}^{(-i)})^{\theta_{i}} U_{x}^{i}(t, x)) = V^{i}(t, \widehat{c_{t}}^{i, *} x) - c_{t}^{i, *} U_{x}(t, x).
\end{eqnarray*}
Requiring \eqref{cond_itocase} to vanish is sufficient to ensure the martingale optimality of the relative performance criterion $Q^{i}$. By concavity of $U^{i}$, the drift of $Q^{i}(t, \widehat{X_{t}}^{i})$ is then non-positive for any admissible strategy and zero under the optimal control $(\pi_{t}^{i, *}, c_{t}^{i, *})$. A sufficient condition can therefore be expressed on the drift $\beta^{i}$, namely
\begin{eqnarray} \label{eq:driftconstraint}
\beta^{i}(t, x) = -  \varphi(t, x) + \frac{1}{2} U_{xx}^{i}(t, x) x^{2} (\nu_{i}^{2} + \sigma_{i}^{2}) \PAR{\pi_{t}^{i, *}(x)}^{2} - \widetilde{V}^{i}(t, (\overset{\sim}{c_{t}}^{(-i)})^{\theta_{i}} U_{x}^{i}(t, x)).
\end{eqnarray}
Substituting this identity into the It\^o decomposition of $U^{i}$ yields the SPDE \eqref{utspde_first}.

\paragraph{Admissibility -} We now verify that Assumptions \ref{ass:regconc}, \ref{ass:admport} and \ref{ass:admcons} are sufficient to ensure that the optimal relative wealth process $\widehat{X_{t}}^{i, *}$ admits a strong global solution, by an application of Theorem 1.4 from \cite{lan2014new}. From the relative wealth dynamics \eqref{xhatdyn}, the It\^o decomposition of $\widehat{X_{t}}^{i, *}$ has drift and volatility local characteristics
\begin{align}
\beta^X(t,x) &= x\Bigl(\xi_{i} + \theta \overline{c_{t}}^{(-i)} +\pi_t^{i, *}(x)\bigl(\mu_i-\theta_i\,\sigma_i\,\overline{\sigma\pi_t}^{\,(-i)}\bigr)- c_t^{i, *}(x)\Bigr) \\
\gamma^{X}(t, x) &= x \PAR{\nu_{i} \pi_{t}^{i, *}(x), \, \sigma_{i} \pi_{t}^{i, *}(x) - \theta_{i} \overline{\sigma \pi_{t}}^{(-i)}, \, \PAR{\frac{\theta_{i}}{n-1} \nu_{k} \pi_{t}^{k}}_{k \neq i}}.
\end{align}
From the optimal policy \eqref{opt_port} and \eqref{optcons}, for any $x > 0$,
\begin{eqnarray*}
\ABS{x\pi_{t}^{i, *}(x)} &\leq& \frac{1}{\nu_{i}^{2} + \sigma_{i}^{2}} \sigma_{i} \theta_{i} \overline{\sigma \pi_{t}}^{(-i)} x + \ABS{\Sigma_{i}.\frac{\gamma_{x}^{i}(t, x)}{U_{xx}^{i}(t, x)}}
+ \ABS{\mu_{i} - \theta_{i} \sigma_{i} \overline{\sigma \pi_{t}}^{(-i)}}\ABS{\frac{U_{x}^{i}(t, x)}{U_{xx}^{i}(t, x)}} \\
\ABS{x c_{t}^{i, *}(x)} &=& \ABS{(V_{x}^{i})^{-1} \PAR{t, U_{x}^{i}(t, x)(\overset{\sim}{c_{t}}^{(-i)})^{\theta_{i}}} (\overset{\sim}{c_{t}}^{(-i)})^{\theta_{i}}}.
\end{eqnarray*}
We next obtain a quadratic growth bound for $\NRM{\gamma^{X}(t, x)}^{2} + 2 x \beta^{X}(t, x)$. First
\begin{eqnarray*}
\NRM{\gamma^{X}(t, x)}^{2} &\leq& \nu_{i}^{2}\ABS{x\pi_{t}^{i, *}(x)}^{2} + 2 \sigma_{i}^{2} \ABS{x \pi_{t}^{i, *}(x)}^{2} + 2\ABS{\theta_{i} \overline{\sigma \pi_{t}}^{(-i)} x}^{2} + \sum_{k \neq i} \ABS{x\frac{\theta_{i}}{n-1} \nu_{k} \pi_{t}^{k}}^{2}. \\
x\beta^{X}(t, x) &\leq& \ABS{\xi_{i}} x^{2} + \theta \overline{c_{t}}^{(-i)} x^{2} + \ABS{\mu_i-\theta_i\,\sigma_i\,\overline{\sigma\pi_t}^{\,(-i)}} \ABS{\pi_t^{i, *}(x)} x^{2} + \ABS{c_{t}^{i, *}(x)}x^{2}.
\end{eqnarray*}
Since $(\pi_{t}^{i, *}(x), c_{t}^{i, *}(x)) \in \A$, and $(\pi_{t}^{k}, c_{t}^{k})_{k\neq i}$ are admissible by assumption, there exist adapted processes $f^{\gamma}_{t}, \, f^{\beta}_{t} \in \Hr^{1}_{\loc}$ such that a.s for all $t \geq 0$ and $x > 0$:
\begin{eqnarray*}
    \NRM{\gamma^{X}(t, x)}^{2} + x\beta^{X}(t, x) &\leq& (f_{t}^{\gamma} + f^{\beta}_{t}) (x^{2}+1).
\end{eqnarray*}
Consequently, under Assumptions \ref{ass:regconc}, \ref{ass:admport} and \ref{ass:admcons}, Theorem 1.4 from \cite{lan2014new} applies and ensures global existence and uniqueness of a strong solution $\widehat{X_{t}}^{i, *}$.

\noindent
\textbf{Local martingale optimality -} To complete the martingale optimality argument, we verify that $Q^{i}(t, \widehat{X_{t}}^{i, *})$ is indeed a local martingale. From the It\^o computation \eqref{eq:gammQIto}, the diffusion coefficient of this composite process is
\begin{equation}
    \gamma^{Q^{i}}(t, \widehat{X_{t}}^{i, *}) = \begin{pmatrix}
\widehat{X_{t}}^{i, *} U_{x}(t, \widehat{X_{t}}^{i, *}) \nu_{i} \pi_{t}^{i, *} + \gamma^{i, W}(t, \widehat{X_{t}}^{i}) \\ 
\PAR{- \widehat{X_{t}}^{i, *}U_{x}(t, \widehat{X_{t}}^{i, *}) \frac{\theta_{i} \nu_{k} \pi_{t}^{k}}{n-1} }_{k \neq i} \\ 
\widehat{X_{t}}^{i, *}U_{x}(t, \widehat{X_{t}}^{i, *}) \PAR{\sigma_{i} \pi_{t}^{i, *} - \theta_{i}\overline{\sigma \pi_{t}}^{(-i)}} + \gamma^{i, B}(t, \widehat{X_{t}}^{i, *})
\end{pmatrix}
\end{equation}
By the arguments in the preceding paragraph and the regularity imposed on the utility random field, the process $\widehat{X_{t}}^{i, *} U_{x}(t, \widehat{X_{t}}^{i, *})$ is continuous. Moreover, under Assumption \ref{ass:regconc}, the volatility characteristics $\gamma^{i} = (\gamma^{i, W}, \, \gamma^{i, B})$ are $C^{2}$ in the state variable. In addition, $\pi_{t}^{i, *} \in \A$ by the above and $\pi_{t}^{k} \in \A$ for $k \neq i$ by assumption. Consequently, the integrand $\gamma^{Q^{i}}(t, \widehat{X_{t}}^{i, *})$ is square integrable on every compact time interval, so the stochastic integral term is a well-defined local martingale. Hence, $Q^{i}(t, \widehat{X_{t}}^{i, *})$ is a local martingale.
\end{preuve}

\begin{rmq}
\begin{enumerate}
\item Strengthening the above conclusion to a true martingale property for $Q^{i}(t, \widehat{X_{t}}^{i, *})$ would require additional uniform integrability, namely $\gamma^{i}(t, \widehat{X_{t}}^{i, *}) \in \Hr_{\loc}^{2}$, together with
\begin{eqnarray} \label{eq:martingale_gen_ass}
\e \SBRA{\int_{0}^{t} \PAR{\widehat{X_{t}}^{i, *} U_{x}(t, \widehat{X_{t}}^{i, *}) \pi_{t}^{i, *}}^{2} dt} < \infty.
\end{eqnarray}
One possible route is to impose polynomial growth of $x U_{x}(t, x)$. However, this would also require higher moment bounds for the relative wealth process and stronger admissibility requirements on the optimal allocation strategy, as boundedness for example. In general, establishing true martingale remains delicate. In Section \ref{sect:consutCRRA}, we discuss sufficient conditions ensuring this martingale property for separable CRRA utility random fields.
\item Duality‐based admissibility arguments in, e.g., \cite{nicole2013exact} and \cite{el2018consistent} apply to the single agent setting and reduce admissibility to conditions expressed on the marginal volatility. In the competitive framework, however, both the relative wealth dynamics and the state price density dynamics acquire additional interaction terms, which prevents simplifications and makes the calculations intractable.
\end{enumerate}
\end{rmq}

\begin{rmq} \label{rmq:prefdesign}
Beyond the best-response problem for a fixed forward criterion, one may consider a \textit{preference-design} problem for a single agent. Let $\A_{pref}^{i}$ denote the class of admissible pairs of random utility fields $(U^{i}, V^{i})$, with fixed initial specification $(u_{0}^{i}, v_{0}^{i})$ satisfying the consistency SPDE \eqref{utspde_first}. For each $(U^{i}, V^{i}) \in \A_{pref}^{i}$, the best-response control $(\pi_{t}^{i, *}, c_{t}^{i, *})$ is uniquely determined by Theorem \ref{SPDEprop}, and the associated relative performance criterion denoted $Q^{i, (U^{i}, V^{i})}$ satisfies the martingale optimality property. It is then natural to ask wether one can select $(U^{i}, V^{i}) \in \A_{pref}^{i}$ that maximizes the expected payoff of agent $i$ at all future dates, namely
\begin{eqnarray*}
    \e \SBRA{Q^{i, (U^{i}, V^{i})}(t, \widehat{X_{t}}^{i, *}) | \F_{s}}, \quad s \leq t.
\end{eqnarray*}
However this formulation is typically degenerate: a criterion that dominates all others for all times can exist only when $\A_{pref}^{i}$ is a singleton. A more meaningful formulation of preference design problem is therefore to cast it as an inverse investment (or utility-recovery) problem, see e.g \cite{kallblad2020black}, \cite{elkaroui2024}. Specifically, given an observed (or target) best-response behavior, such as a strategy or a wealth process, one seeks admissible forward criteria $(U^{i}, V^{i}) \in \A^{i}_{pref}$ whose induced optimal policy reproduces that behavior. This perspective also highlights the role of non-zero utility volatility in preference design problems. Indeed, in the best-response portfolio \eqref{opt_port}, $\gamma^{i}$ can serve as a calibration parameter to match a target investment profile. Whenever there exists an admissible random utility pair $(U^{i}, V^{i})$ solving the consistency SPDE, Theorem \ref{SPDEprop} yields the associated best-response consumption and ensures that the resulting preference criterion remains dynamically consistent. In this sense, non-zero volatility is an additional preference parameter, and characterizing the admissible solution set of the consistency SPDE is key to understanding how the forward utility volatility parametrizes both optimal strategies and preferences.
\end{rmq}

\subsection{Forward relative performance criteria of CRRA type} \label{sect:consutCRRA}

In the following, we focus on forward relative performance criteria of CRRA type, with separable time and space dependence. Let $u^{(\alpha)}(x) = \frac{x^{1 - \alpha}}{1 - \alpha}$ denote the standard power utility function with risk aversion parameter $\alpha \in \PAR{0, 1} \cup (1, +\infty)$. When the wealth utility has power-type dependence in the state variable, the drift restriction \eqref{eq:driftconstraint} becomes explicit. In particular, if the wealth utility $U^{\alpha}$ is of separable power type, the couple $(U^{\alpha}, V)$ generates a forward relative performance criterion provided that the consumption utility $V$ is also of separable power type with the same risk aversion parameter, as observed in the single agent setting of \cite{el2018consistent}. In this case, the HJB-SPDE from Theorem \ref{SPDEprop} reduces to a condition linking the time-dependent factors in $U^{\alpha}$ and $V$, which we make explicit below.

Consider an agent $i$ with a pair of utility random fields $(U^{(\alpha_{i})}, V^{i})$. We assume that the wealth utility is of the form $U^{(\alpha_{i})}(t, x) = Z_{t}^{i} u^{(\alpha_{i})}(x)$, where $\alpha_{i} \in \PAR{0, 1}\cup (1, +\infty)$ and $Z_{t}^{i}$ is an It\^o process. Moreover, we assume that $Z_{t}^{i}$ has lognormal dynamics, 
\begin{eqnarray} \label{eq:dynZ}
    dZ_{t}^{i} = Z_{t}^{i} \PAR{b^{Z^{i}}(t)dt + \delta^{Z^{i}, W}(t) dW_{t}^{i} + \delta^{Z^{i}, B}(t) dB_{t}}, \quad \text{with} \quad Z_{0}^{i} = z_{0}^{i} > 0,
\end{eqnarray}
with $b^{Z^{i}} \in \Hr_{\loc}^{1}$ and $\delta^{Z^{i}, W}, \, \delta^{Z^{i}, B} \in L^{\infty}_{\loc}$. The local characteristics of $U^{(\alpha_{i})}$ of type \eqref{Itodynut} are then given by
\begin{eqnarray} \label{local_char}
    \beta^{(\alpha_{i})}(t, x) = b^{Z^{i}}(t) U^{(\alpha_{i})}(t, x) \quad \text{and} \quad \gamma^{(\alpha_{i})}(t, x) = \PAR{\delta^{Z^{i}, W}(t) \,, \, \delta^{Z^{i}, B}(t)} U^{(\alpha_{i})}(t, x).
\end{eqnarray}
In the following, we may denote the two-dimensional volatility component by $\delta^{Z^{i}} = (\delta^{Z^{i}, W}, \, \delta^{Z^{i}, B})$. For ease of notation, we introduce for $i = 1, ..., n$,
\begin{eqnarray}
    f^{i}(t) &=&  - \theta_{i} \overline{\mu \pi_{t}}^{(-i)}  + \frac{\theta_{i}}{2} \overline{\Sigma \pi_{t}^{2}}^{(-i)} + \frac{\theta_{i}^{2}}{2} \PAR{(\overline{\sigma \pi_{t}}^{(-i)})^{2} + \frac{1}{n-1} \overline{(\nu \pi_{t})^{2}}^{(-i)}} \label{feq} \\
    g^{i}(t) &=& \frac{\theta_{i}^{2}}{n-1} \overline{(\nu \pi_{t})^{2}}^{(-i)} + (\theta_{i} \overline{\sigma \pi_{t}}^{(-i)})^{2}. \label{geq} 
\end{eqnarray}
In the CRRA setting, the optimal strategies are independent from the state variable $x$. Therefore, in this section, we drop the dependence in $x$ from the notation. We first make explicit the link between wealth and consumption utilities implied by the consistency SPDE \eqref{utspde_first}.

\begin{prop}
    Let $(U^{(\alpha_{i})}, V^{i})$ be a couple of utility random fields satisfying the SPDE \eqref{utspde_first}, where the wealth utility is assumed to be of CRRA type $U^{(\alpha_{i})}(t, x) = Z_{t}^{i} u^{(\alpha_{i})}(x)$. Then, the consumption utility is necessarily also separable of CRRA type, with the same risk aversion parameter, of the form $V^{(\alpha_{i})}(t, c) = \phi_{t}^{i} u^{(\alpha_{i})}(c)$. 
\end{prop}

\begin{preuve}
Substituting \eqref{local_char} into \eqref{eq:driftconstraint} yields
\begin{eqnarray} \label{lienVU}
    \widetilde{V}^{i}(t, (\overset{\sim}{c_{t}}^{(-i)})^{\theta_{i}} U_{x}^{(\alpha_{i})}(t, x)) = (\overline{v_{t}^{i}} - (1 - \alpha_{i})\theta_{i} \overline{c_{t}}^{(-i)}) U^{(\alpha_{i})}(t, x),
\end{eqnarray}
where $\overline{v_{t}^{i}} = - b^{Z^{i}}(t) + (1 - \alpha_{i}) \PAR{ \delta^{Z^{i}, B}(t) \theta_{i} \overline{\sigma \pi_{t}}^{(-i)} + f^{i}(t) - \frac{1}{2} \alpha_{i}(g^{i}(t) - (\nu_{i}^{2} + \sigma_{i}^{2})\pi_{t}^{i, *}) }$. Differentiating \eqref{lienVU} with respect to the space variable gives
\begin{eqnarray*}
\frac{ d \widetilde{V}^{i}\PAR{t, (\overset{\sim}{c_{t}}^{(-i)})^{\theta_{i}} U_{x}^{(\alpha_{i})}(t, x)} }{dx} = \frac{(1 - \alpha_{i})\theta_{i} \overline{c_{t}}^{(-i)} - \overline{v_{t}^{i}}}{\alpha_{i}( \overset{\sim}{c_{t}}^{(-i)})^{\theta_{i}} } x.
\end{eqnarray*}
Recalling that $V_{x}(t, -\overset{\sim}{V}_{y}(t, y)) = y$, and integrating yields
\begin{eqnarray} \label{V_power}
    V^{i}(t, x) = \PAR{\frac{\overline{v_{t}^{i}} - (1 - \alpha_{i})\theta_{i} \overline{c_{t}}^{(-i)}}{\alpha_{i}}}^{\alpha_{i}} (\overset{\sim}{c_{t}}^{(-i)})^{\theta_{i}(1 - \alpha_{i})} U^{(\alpha_{i})}(t, x).
\end{eqnarray}
Thus, $V^{i}$ is necessarily of power type, proportional to the wealth utility $U^{(\alpha_{i})}$. 
\end{preuve}

\vspace{0.3cm}
\begin{rmq} \label{rmq:linkVU}
The link between $U^{(\alpha_{i})}$ and $V^{i}$ highlighted in \eqref{V_power} holds here only for CRRA utilities, since it relies on the fact that the first and second derivatives of $U^{(\alpha_{i})}$ can be expressed as functions of $U^{(\alpha_{i})}$ itself. In general, define the random field $H$ by
\begin{small}
\begin{eqnarray}
H(t, x) = \frac{-1}{(\overset{\sim}{c_{t}}^{(-i)})^{\theta_{i}} U_{xx}^{i}(t, x)} \PAR{-\beta_{x}^{i}(t, x) - \varphi_{x}(t, x) +\frac{1}{2}(\nu_{i}^{2} + \sigma_{i}^{2}) \PAR{\pi_{t}^{i, *}(x)}^{2}\PAR{ x^{2}U_{xxx}^{i}(t, x) + 2xU_{xx}^{i}(t, x)}}
\end{eqnarray}
\end{small}
Then, provided that the pair $(U^{(\alpha_{i})}, V^{i})$ solves the HJB-SPDE \eqref{utspde_first} and that $H$ admits an inverse with respect to the space variable, denoted $H^{-1}$, the utility from consumption is characterized by
\begin{eqnarray*}
V_{x}^{i}(t, x) = (\overset{\sim}{c_{t}}^{(-i)})^{\theta_{i}} U_{x}^{i}(t, H^{-1}(t, x)).
\end{eqnarray*}
\end{rmq}

From now on, we set $V^{(\alpha_{i})}(t, c) = \phi_{t}^{i} u^{(\alpha_{i})}(c)$, for some It\^o process $\phi_{t}^{i}$. The time consistency constraint \eqref{eq:driftconstraint} reduces to a condition on the drift coefficient $b^{Z^{i}}$, which depends on $Z_{t}^{i}$ and $\phi_{t}^{i}$. We make this dependence explicit below and show that the resulting dynamics are well-posed, with a unique global strong solution.

\begin{prop} \label{prop:VfctUnjoueurs}
Let $(U^{(\alpha_{i})}, V^{(\alpha_{i})})$ be a pair of separable CRRA utility random fields, as defined above. Then $(U^{(\alpha_{i})}, V^{(\alpha_{i})})$ satisfy the HJB-SPDE \eqref{utspde_first} if and only if the drift local characteristic of $Z_{t}^{i}$ takes the form
\begin{eqnarray} \label{driftcond1}
b^{Z^{i}}(t, Z^{i}_{t}, \phi_{t}^{i}) = \overline{b^{Z^{i}}}(t) - \frac{\alpha_{i}}{(\overset{\sim}{c_{t}}^{(-i)})^{\frac{\theta_{i}(1 - \alpha_{i})}{\alpha_{i}}}} \PAR{\frac{\phi_{t}^{i}}{Z_{t}^{i}}}^{\frac{1}{\alpha_{i}}},
\end{eqnarray}
where 
\begin{eqnarray} \label{eq:defA_{i}condb}
\overline{b^{Z^{i}}}(t) = (1 - \alpha_{i}) \PAR{ \delta^{Z^{i}, B}(t) \theta_{i} \overline{\sigma \pi_{t}}^{(-i)} - f^{i}(t) + \frac{1}{2} \alpha_{i} \PAR{g^{i}(t) - (\nu_{i}^{2} + \sigma_{i}^{2}) \pi_{t}^{i, *}} - \theta_{i} \overline{c_{t}}^{(-i)}}.
\end{eqnarray}
Then, the SDE \eqref{eq:dynZ} admits a unique global strong solution $Z_{t}^{i}$.
\end{prop}

\begin{preuve}
Substituting $V^{i}(t, c) = \phi_{t}^{i} u^{\alpha_{i}}(t, c)$ in \eqref{V_power} yields 
\begin{eqnarray*} \label{barvnagent}
\overline{v_{t}^{i}} = (1 - \alpha_{i})\theta_{i} \overline{c_{t}}^{(-i)} + \frac{\alpha_{i}}{(\overset{\sim}{c_{t}}^{(-i)})^{\frac{\theta_{i}(1 - \alpha_{i})}{\alpha_{i}}}} \PAR{\frac{\phi_{t}^{i}}{Z_{t}^{i}}}^{\frac{1}{\alpha_{i}}},
\end{eqnarray*}
which gives the desired characterization of $b^{Z^{i}}$. Plugging \eqref{driftcond1} into \eqref{eq:dynZ}, we obtain
\begin{eqnarray*}
dZ_{t}^{i} = \SBRA{\overline{b^{Z^{i}}}(t) Z_{t}^{i} - \frac{\alpha_{i}}{(\overset{\sim}{c_{t}}^{(-i)})^{\frac{\theta_{i}(1 - \alpha_{i})}{\alpha_{i}}}} (\phi_{t}^{i})^{\frac{1}{\alpha_{i}}} (Z_{t}^{i})^{1 - \frac{1}{\alpha_{i}}}}dt + Z_{t}^{i} \delta^{Z^{i}}(t).d\overline{W_{t}^{i}}.
\end{eqnarray*}
Let $L_{t}^{i} = (Z_{t}^{i})^{\frac{1}{\alpha_{i}}}$. Applying It\^o's formula yields 
\begin{eqnarray}
dL_{t}^{i} = \SBRA{\PAR{\frac{1}{\alpha_{i}} \overline{b^{Z^{i}}}(t) + \frac{1}{2} \frac{1}{\alpha_{i}} (\frac{1}{\alpha_{i}} -1) \NRM{\delta^{Z^{i}}(t)}^{2}} L_{t}^{i} - \frac{1}{(\overset{\sim}{c_{t}}^{(-i)})^{\frac{\theta_{i}(1 - \alpha_{i})}{\alpha_{i}}}} (\phi_{t}^{i})^{\frac{1}{\alpha_{i}}}}dt + L_{t}^{i}\frac{1}{\alpha_{i}} \delta^{Z^{i}}(t).d\overline{W_{t}^{i}}. \label{eq:dynY}
\end{eqnarray}
Define $\displaystyle\E_{t}^{i} = \exp \PAR{\int_{0}^{t} (\frac{1}{\alpha_{i}} \overline{b^{Z^{i}}}(s) - \frac{1}{2} \frac{1}{\alpha_{i}} (\frac{1}{\alpha_{i}} + 1) \NRM{\delta^{Z^{i}}(s)}^{2})ds + \int_{0}^{t} \frac{\delta^{Z^{i}}(s)}{\alpha_{i}}.d\overline{W_{s}^{i}}}$. Since \eqref{eq:dynY} is linear, it admits the explicit solution
\begin{eqnarray*} \label{eq:explinvZ}
    L_{t}^{i} = \E_{t}^{i} \PAR{Y_{0}^{i} - \int_{0}^{t} \frac{1}{(\tilde{c_{s}}^{(-i)})^{\frac{\theta_{i}(1 - \alpha_{i})}{\alpha_{i}}}} (\phi_{s}^{i})^{\frac{1}{\alpha_{i}}} (\E_{s}^{i})^{-1}ds}.
\end{eqnarray*}
Finally $Z_{t}^{i} = (L_{t}^{i})^{\alpha_{i}}$, concluding on global existence of a solution to \eqref{eq:dynZ}.
\end{preuve}

Despite global existence and uniqueness of a strong solution $Z_{t}^{i}$ under condition \eqref{driftcond1}, there is no guarantee that the solution remains strictly positive at all times. A minimal condition for $Z_{t}^{i} > 0$ for all $t \geq 0$ is
\begin{eqnarray*} \label{eq:minadmcondZ}
\int_{0}^{t} \frac{1}{(\tilde{c_{s}}^{(-i)})^{\frac{\theta_{i}(1 - \alpha_{i})}{\alpha_{i}}}} (\phi_{s}^{i})^{\frac{1}{\alpha_{i}}} (\E_{s}^{i})^{-1}ds < L_{0}^{i}, \quad \text{for all } t \geq 0.
\end{eqnarray*}
Since the stochastic exponential $\E_{t}^{i}$ is random and typically unbounded, deriving sufficient conditions for this inequality is difficult. A tractable approach is to impose additional structure on $\phi^{i}$, which allows direct control of the optimal consumption process.

We are now in a position to state the main result of this section on best-response strategies in the separable CRRA framework. Observe that the optimal consumption process depends on the ratio $\PAR{\frac{\phi_{t}^{i}}{Z_{t}^{i}}}^{\frac{1}{\alpha_{i}}}$. Moreover, obtaining that the optimal preference criterion $Q^{i}(t, \widehat{X_{t}}^{i, *})$ is a true martingale requires additional integrability, which we enforce through assumptions on the consumption utility factor $\phi^{i}$. Specifically, we work with a consumption utility random field of the form $V^{(\alpha_{i})}(t, c) = \phi_{t}^{i} u^{(\alpha_{i})}(c)$ where
\begin{eqnarray} \label{eq:dynphi}
    d\phi_{t}^{i} = \phi_{t}^{i} \PAR{b^{\phi^{i}}(t, Z_{t}^{i}, \phi_{t}^{i})dt + \delta^{\phi^{i}, W}(t) dW_{t}^{i} + \delta^{\phi^{i}, B}(t) dB_{t}}, \quad \text{with} \quad \phi_{0}^{i} > 0.
\end{eqnarray}

\begin{thm} \label{thm:optstratpowerN}
Consider an agent $i \in \BRA{1, \dots, n}$, with a pair of separable utility random fields of CRRA type $U^{(\alpha_{i})}(t, x) = Z_{t}^{i} u^{(\alpha_{i})}(x)$ and $V^{(\alpha_{i})}(t, c) = \phi_{t}^{i} u^{(\alpha_{i})}(c)$. Assume that $\PAR{U^{(\alpha_{i})}, V^{(\alpha_{i})}} $ satisfies the HJB-SPDE \eqref{utspde_first}. Then, the best-response strategy $(\pi_{t}^{i, *}, c_{t}^{i, *})$ takes the form
\begin{eqnarray}
 \pi_{t}^{i, *} &=& \frac{1}{\nu_{i}^{2} + \sigma_{i}^{2}} \left( \sigma_{i} \theta_{i} \overline{\sigma \pi_{t}}^{(-i)} + \frac{1}{\alpha_{i}} \left(\delta^{Z^{i}, W}(t) \nu_{i} + \delta^{Z^{i}, B}(t) \sigma_{i} + (\mu_{i} - \theta_{i} \sigma_{i} \overline{\sigma \pi_{t}}^{(-i)}) \right) \right) \\
c_{t}^{i, *} &=& (\overset{\sim}{c_{t}}^{(-i)})^{\theta_{i} (1 - \frac{1}{\alpha_{i}})} \PAR{\frac{\phi_{t}^{i}}{Z_{t}^{i}}}^{\frac{1}{\alpha_{i}}}. \label{eq:optconsCRRAn}
\end{eqnarray}
\textbf{Admissibility -} Assume in addition that the drift of the consumption utility factor $\phi^{i}$ contains a mean-reverting term of the form 
\begin{eqnarray} \label{eq:condphi}
    b^{\phi^{i}}(t, Z_{t}^{i}, \phi_{t}^{i}) = \overline{b^{\phi^{i}}}(t) - \lambda \frac{\alpha_{i}}{(\overset{\sim}{c_{t}}^{(-i)})^{\frac{\theta_{i}(1 - \alpha_{i})}{\alpha_{i}}}} \PAR{\frac{\phi_{t}^{i}}{Z_{t}^{i}}}^{\frac{1}{\alpha_{i}}},
\end{eqnarray}
where $\overline{b^{\phi^{i}}}$ is locally bounded and $\lambda \geq 1$. Assume further that:
\begin{itemize}
    \item if $\alpha_{i} \in (0, 1)$, there exists $\eta_{i} > \frac{\theta_{i}(1 - \alpha_{i})}{\alpha_{i}}$ such that $(\overset{\sim}{c_{t}}^{(-i)})^{- \eta_{i}} \in \Hr_{\loc}^{1}$,
    \item if $\alpha_{i} > 1$, there exists $\eta_{i} > \frac{\theta_{i}(\alpha_{i} - 1)}{\alpha_{i}}$ such that $(\overset{\sim}{c_{t}}^{(-i)})^{ \eta_{i}} \in \Hr_{\loc}^{1}$.
\end{itemize}
Then, the best-response strategy $(\pi_{t}^{i, *}, c_{t}^{i, *})$ is admissible.  

\noindent
\textbf{Martingale optimality -} Assume in addition that:
\begin{itemize}
    \item if $\alpha_{i} \in (0, 1)$, there exists $m_i>2(1-\alpha_i)\theta_i$ such that for all $t \geq 0$
    \[
    \e\!\left[\exp\!\left(\int_0^t m_i\,\overline c_s^{(-i)}\,ds\right)\right]<\infty.
    \]
    \item if $\alpha_i\in(1,\infty)$, and $\lambda > 1$, there exists $m_i > \frac{\theta_{i}(\alpha_{i} - 1)}{\alpha_{i}} \max \PAR{1, \frac{2(\alpha_{i} -1)}{(\lambda - 1)\alpha_{i}}}$ such that $(\tilde{c_{t}}^{(-i)})^{m_i} \in \Hr_{\loc}^{1}$.
\end{itemize} 
Under these conditions, the relative performance criterion at optimum $Q^{i}(t, \widehat{X_{t}}^{i, *})$ is a true martingale.
\end{thm}
Note that the optimal portfolio $\pi_{t}^{i, *}$ associated with separable CRRA forward relative criteria depends only on model coefficients, the diffusion terms of $Z^{i}$ assumed bounded and the given population term $\overline{\sigma \pi_{t}}^{(-i)}$, so that admissibility is straightforward.

\begin{preuve}

\textbf{Well-posedness and admissibility -} First, the optimal portfolio $\pi_{t}^{i, *}$ is admissible provided that volatility coefficients $\delta^{Z^{i}, W}$ and $\delta^{Z^{i}, B}$ are $\Hr_{\loc}^{2}$. Next, we address admissibility of consumption. In the separable CRRA case, Assumption \ref{ass:admcons} is equivalent to
\begin{equation*}
 \e \SBRA{ \int_{0}^{T}\Bigl(\frac{\phi_s^{i}}{Z_s^{i}}\Bigr)^{1/\alpha_{i}} (\overset{\sim}{c_{s}}^{(-i)})^{\theta_{i}}ds}\,<\infty,\qquad\forall T>0.
\end{equation*}
Fix $p > 1$ and let $q = \frac{p}{p-1}$. By Hölder's inequality, 
\begin{eqnarray} \label{eq:holderAss1.3CRRA}
\e \SBRA{ \int_{0}^{T}\Bigl(\frac{\phi_s^{i}}{Z_s^{i}}\Bigr)^{1/\alpha_{i}} (\overset{\sim}{c_{s}}^{(-i)})^{\theta_{i}}}\,ds \leq \PAR{\e \SBRA{ \int_{0}^{T}\Bigl(\frac{\phi_s^{i}}{Z_s^{i}}\Bigr)^{p/\alpha_{i}}ds}}^{\frac{1}{p}} \PAR{\e \SBRA{ \int_{0}^{T} (\overset{\sim}{c_{s}}^{(-i)})^{q \theta_{i}}\,ds}}^{\frac{1}{q}}.
\end{eqnarray}
To exploit the assumed $\Hr_{\loc}^{1}$ integrability of competitors' consumption, choose $p$ large enough so that $q \theta_{i} \leq 1$, i.e $p \geq \frac{1}{1 - \theta_{i}}$. It then suffices to show that $Y_{t}^{i} = \Bigl(\frac{\phi_s^{i}}{Z_s^{i}}\Bigr)^{1/\alpha_{i}} \in \Hr_{\loc}^{p}$ to ensure that \eqref{eq:holderAss1.3CRRA} is finite. Applying It\^o's formula yields
\begin{eqnarray} \label{eq:logistic_R}
d Y_{t}^{i} = Y_{t}^{i} \PAR{ \SBRA{b^{\phi^{i}}(t, Z_{t}^{i}, \phi_{t}^{i}) + \overline{b^{Y^{i}}}(t) + \frac{\alpha_{i}}{(\overset{\sim}{c_{t}}^{(-i)})^{\frac{\theta_{i}(1 - \alpha_{i})}{\alpha_{i}}}} Y_{t}^{i}}dt + \delta^{Y^{i}}(t). d\overline{W_{t}^{i}}}
\end{eqnarray}
where 
\begin{eqnarray*}
    \overline{b^{Y^{i}}}(t) &=&  - \overline{b^{Z^{i}}}(t)  - \NRM{\delta^{Z^{i}}(t)}^{2} - \delta^{\phi^{i}}.\delta^{Z^{i}}(t) + \frac{\alpha_{i}}{2} \PAR{\frac{1}{\alpha_{i}} - 1} \NRM{\delta^{Y^{i}}(t)}^{2} \\
    \delta^{Y^{i}}(t) &=& \frac{1}{\alpha_{i}}(\delta^{\phi^{i}}(t) - \delta^{Z^{i}}(t))
\end{eqnarray*}
Equation \eqref{eq:logistic_R} has a positive quadratic drift, characterizing non-autonomous logistic SDEs. 

\noindent
\textbf{State-feedback condition -} A sufficient condition for global existence and uniqueness of a solution to \eqref{eq:logistic_R} is that the drift of the consumption factor $\phi^{i}$ includes a mean-reverting condition, of type \eqref{eq:condphi}. Under this assumption, equation \eqref{eq:logistic_R} becomes 
\begin{eqnarray} \label{eq:logYfin}
d Y_{t}^{i} = Y_{t}^{i} \PAR{ \SBRA{\overline{b^{Y^{i}}}(t) - (\lambda - 1)\frac{\alpha_{i}}{(\overset{\sim}{c_{t}}^{(-i)})^{\frac{\theta_{i}(1 - \alpha_{i})}{\alpha_{i}}}} Y_{t}^{i}}dt + \delta^{Y^{i}}(t). d\overline{W_{t}^{i}}}.
\end{eqnarray}
For $\lambda \geq 1$, the quadratic term in the drift is non-positive, and \eqref{eq:logYfin} is a logistic SDE. Such equations, originating in population dynamics, have been extensively studied, see, for instance, \cite{jiang2005note} and \cite{giet2015logistic}. Global existence and uniqueness for \eqref{eq:logYfin} follow from an inverse transformation. Applying It\^o's formula yields
\begin{eqnarray} \label{eq:linearizedY}
d \PAR{\frac{1}{Y_{t}^{i}}} = \SBRA{(\lambda - 1) \frac{\alpha_{i}}{(\tilde{c_{t}})^{\frac{\theta_{i}(1 - \alpha_{i})}{\alpha_{i}}}} + \PAR{\frac{1}{Y_{t}^{i}}} \PAR{\NRM{\delta^{Y^{i}}(t)}^{2} - \overline{b^{Y^{i}}}(t)}}dt - \PAR{\frac{1}{Y_{t}^{i}}} \delta^{Y^{i}}(t).d\overline{W_{t}^{i}}.
\end{eqnarray}
Recall that $\delta^{Y^{i}}$ and $\overline{b^{Y^{i}}}$ belong to $L_{\loc}^{\infty}$ by assumption. The above linear SDE admits an explicit strong solution provided that $(\tilde{c}_{t}^{(-i)})^{- \frac{\theta_{i}(1 - \alpha_{i})}{\alpha_{i}}} \in \Hr_{\loc}^{1}$. Let $\E_{t}^{i}$ be defined by 
\begin{eqnarray*}
\E_{t}^{i} = \exp \PAR{\int_{0}^{t} \PAR{ \frac{1}{2} \NRM{\delta^{Y^{i}}(s)}^{2} - \overline{b^{Y^{i}}}(s) }ds - \int_{0}^{t} \delta^{Y^{i}}(s).d\overline{W_{s}^{i}}}.
\end{eqnarray*}
Then \eqref{eq:linearizedY} has the unique global strong solution
\begin{eqnarray} \label{eq:explsolinvY}
\frac{1}{Y_{t}^{i}} = \E_{t}^{i} \PAR{ \frac{1}{Y_{0}^{i}} + \int_{0}^{t} (\E_{s}^{Y^{i}})^{-1} (\lambda - 1) \frac{\alpha_{i}}{(\tilde{c_{s}})^{\frac{\theta_{i}(1 - \alpha_{i})}{\alpha_{i}}}} ds}.
\end{eqnarray}
Consequently, taking the inverse yields
\begin{eqnarray} \label{eq:explY}
Y_{t}^{i} = \frac{(\E_{t}^{i})^{-1}}{\displaystyle \frac{1}{Y_{0}^{i}} + \int_{0}^{t} (\E_{s}^{Y^{i}})^{-1} (\lambda - 1) \frac{\alpha_{i}}{(\tilde{c_{s}})^{\frac{\theta_{i}(1 - \alpha_{i})}{\alpha_{i}}}} ds }.
\end{eqnarray}
Since $\lambda \geq 1$ and $\tilde{c_{s}}^{(-i)} > 0$, for any $r > 0$
\begin{eqnarray} \label{eq:Ysquaredexpl}
(Y_{t}^{i})^{r} \leq (Y_{0}^{i})^{r} (\E_{t}^{i})^{-r} = (Y_{0}^{i})^{r} \exp \PAR{\int_{0}^{t} \PAR{r \overline{b^{Y^{i}}}(s) - \frac{r}{2}\NRM{\delta^{Y^{i}}(s)}^{2} }ds + r \int_{0}^{t} \delta^{Y^{i}}(s).d\overline{W_{s}^{i}}},
\end{eqnarray}
If $\delta^{Z^{i}}$ and $\delta^{\phi^{i}}$ belong to $L_{\loc}^{\infty}$, then $\overline{b^{Y^{i}}}$ and $\delta^{Y^{i}}$ also belong to $L_{\loc}^{\infty}$. In particular, \eqref{eq:Ysquaredexpl} ensures that $Y_{t}^{i} \in \Hr_{\loc}^{r}$. Consequently, the logistic SDE \eqref{eq:logistic_R} admits a unique global positive strong solution, and the optimal consumption strategy \eqref{eq:optconsCRRAn} associated with separable CRRA utilities is admissible. It remains to verify that $c_{t}^{i, *} \in \Hr_{\loc}^{1}$. From \eqref{eq:optconsCRRAn}, Hölder's inequality yields 
\begin{eqnarray*}
\e \SBRA{\int_{0}^{t} c_{s}^{i, *} ds} \leq \PAR{\e \SBRA{\int_{0}^{t} (\tilde{c_{s}}^{(-i)})^{p \theta_{i} (1 - \frac{1}{\alpha_{i}})}ds}}^{\frac{1}{p}} \PAR{\e \SBRA{\int_{0}^{t} \PAR{\frac{\phi_{s}^{i}}{Z_{s}^{i}}}^{\frac{q}{\alpha_{i}}}ds }}^{\frac{1}{q}},
\end{eqnarray*}
for conjugate exponents $p, \, q > 1$. From the integrability assumptions on $(\overset{\sim}{c_{t}}^{(-i)})$, specified according to whether $\alpha_{i} < 1$ or $\alpha_{i} > 1$, there exists $p$ sufficiently close to $1$, such that $(\overset{\sim}{c_{t}}^{(-i)})^{p \theta_{i} \PAR{1 - \frac{1}{\alpha_{i}}}} \in \Hr_{\loc}^{1}$. Moreover, since $Y_{t}^{i} \in \Hr_{\loc}^{r}$ for all $r > 0$, the second factor is finite, and therefore $c_{t}^{i, *} \in \Hr_{\loc}^{1}.$

\noindent
\textbf{Martingale optimality -} We conclude by proving that the optimal relative performance criterion $Q^{i}(t, \widehat{X_{t}^{i, *}})$ is a true martingale. In the CRRA setting, the sufficient martingale condition \eqref{eq:martingale_gen_ass} reduces to 
\begin{eqnarray*}
\e \SBRA{\int_{0}^{t} \PAR{(\widehat{X_{s}}^{i, *})^{1 - \alpha_{i}} Z_{s}^{i} \pi_{s}^{i, *}}^{2}ds} < +\infty.
\end{eqnarray*}
For simplicity, we assume that $\pi_{t}^{i, *}$ is locally bounded, by requiring $\delta^{Z^{i}} \in L_{\loc}^{\infty}$. Let $p, \, q > 1$ such that $\frac{1}{p} + \frac{1}{q} = 1$. By Hölder's inequality, there exists a constant $C > 0$ such that
\begin{eqnarray} \label{eq:martingale_Holder}
\e \SBRA{\int_{0}^{t} \PAR{(\widehat{X_{s}}^{i, *})^{1 - \alpha_{i}} Z_{s}^{i} \pi_{s}^{i, *}}^{2}ds} \leq C \PAR{\e \SBRA{\int_{0}^{t}  (\widehat{X_{s}}^{i, *})^{2(1 - \alpha_{i})p}ds }}^{\frac{1}{p}} \PAR{\e \SBRA{\int_{0}^{t}  (Z_{s}^{i})^{2q}ds }}^{\frac{1}{q}}.
\end{eqnarray}
Note that, although the optimal relative wealth dynamics \eqref{xhatdyn} is linear, the integrability conditions ensuring that $\widehat{X_{t}}^{i, *}$ admits moments of order $2(1 - \alpha_{i})p$, for some $p > 1$, depend on the risk aversion regime. More precisely,
\begin{itemize}
    \item if $\alpha_{i} \in (0, 1)$, then $2(1-\alpha_i)>0$. In the drift of $(\widehat{X}_t^{i,*})^{2(1-\alpha_i)p}$, the aggregate competitors' consumption $\overline{c}_t^{(-i)}$ enters with a positive sign, whereas $c_t^{i,*}$ enters with a negative sign. Consequently, assuming that there exists $m_i > 2 (1 - \alpha_{i}) \theta_{i}$ such that
    \begin{eqnarray}
    \e \SBRA{\exp \PAR{\int_{0}^{t} m_{i} \bar{c_{s}}^{(-i)} ds}} < \infty, \quad \forall t \geq 0,
    \end{eqnarray} 
    is a convenient sufficient condition for the required positive-moment bounds. In particular, for any $r > 0$ such that $r \theta_{i} < m_i$, one obtains $(\widehat{X}^{i, *})^{r} \in \Hr_{\loc}^{1}.$ Choosing $p > 1$ close enough to $1$ so that $2(1 - \alpha_{i})p \theta_{i} < m_i$ yields the desired Hölder estimate.
    
    \item if $\alpha_{i} \in (1, \infty)$, then $2(1-\alpha_i)<0$ and the term $2(1-\alpha_i)\theta_i\,\overline{c}_t^{(-i)}\le 0$ does not impose additional constraints. However, $c_t^{i,*}$ now enters the drift of $(\widehat{X}_t^{i,*})^{2(1-\alpha_i)p}$ with a positive sign, and one is led to control, for some $p > 1$,
    \begin{eqnarray} \label{eq:martcondalpha>1}
        \e \SBRA{\exp \PAR{2(\alpha_{i} - 1)p \int_{0}^{t}  c_{s}^{i, *}ds}}, \quad \forall t \geq 0.
    \end{eqnarray}
    The solution \eqref{eq:explY} provides an explicit formula for the best-response consumption, namely 
    \begin{eqnarray}
    c_{t}^{i, *} = \frac{(\tilde{c_{t}}^{(-i)})^{\frac{\theta_{i}(\alpha_{i} - 1)}{\alpha_{i}}} (\E_{t}^{i})^{-1}}{\frac{1}{Y_{0}^{i}} + (\lambda - 1)\alpha_{i} \int_{0}^{t} (\E_{s}^{i})^{-1} (\tilde{c_{s}}^{(-i)})^{\frac{\theta_{i}(\alpha_{i} - 1)}{\alpha_{i}}}ds},
    \end{eqnarray}
    from which we obtain 
    \begin{eqnarray}
        \int_{0}^{t} c_{s}^{i, *}ds = \frac{1}{(\lambda - 1)\alpha_{i}} \log \PAR{1 + (\lambda - 1)\alpha_{i} Y_{0}^{i} \int_{0}^{t} (\E_{s}^{i})^{-1}(\tilde{c_{s}}^{(-i)})^{\frac{\theta_{i}(\alpha_{i} - 1)}{\alpha_{i}}}ds}.
    \end{eqnarray}
    Multiplying by $\rho$ and exponentiating leads 
    \begin{eqnarray*}
    \exp \PAR{\rho \int_{0}^{t} c_{s}^{i, *}ds} = \PAR{1 + (\lambda - 1)\alpha_{i} Y_{0}^{i} \int_{0}^{t} (\E_{s}^{i})^{-1}(\tilde{c_{s}}^{(-i)})^{\frac{\theta_{i}(\alpha_{i} - 1)}{\alpha_{i}}}ds}^{\frac{\rho}{(\lambda - 1)\alpha_{i}}}.
    \end{eqnarray*}
    Since we assume $\delta^{Z^{i}}$ and $\delta^{\phi^{i}}$ locally bounded, the process $(\E_{t}^{i})^{-1}$ has moments of all orders on every finite time horizon. In turn, assuming there exists $m_i> \frac{\theta_{i}(\alpha_{i} - 1)}{\alpha_{i}} \max \PAR{1, \frac{2(\alpha_{i} -1)}{(\lambda - 1)\alpha_{i}}}$ such that $(\tilde{c_{t}}^{(-i)})^{m_i} \in \Hr_{\loc}^{1}$ is sufficient to ensure that there exists $p$ close enough to $1$ such that \eqref{eq:martcondalpha>1} is finite.
\end{itemize}
Consequently, using asymmetric Hölder's weights $p$ and $q$, it remains to control moments of $Z_{t}^{i}$ of arbitrary order. The explicit dynamics \eqref{eq:dynY} can be expressed as a function of $c_{t}^{i, *}$ as 
\begin{eqnarray*}
 dL_{t}^{i} = L_{t}^{i} \PAR{K_{t}^{i} - c_{t}^{i, *}}dt + L_{t}^{i}\frac{1}{\alpha_{i}} \delta^{Z^{i}}(t).d\overline{W_{t}^{i}}.
\end{eqnarray*}
This linear SDE admits an explicit solution under boundedness of $\delta^{Z^{i}}, \delta^{\phi^{i}}$ and $c_{t}^{i, *} \in \Hr_{\loc}^{1}$, namely 
\begin{eqnarray*}
L_{t}^{i} = L_{0}^{i} \exp \PAR{\int_{0}^{t} (K_{s}^{i} - c_{s}^{i, *} - \frac{1}{2} \frac{\NRM{\delta^{Z^{i}}(s)}^{2}}{\alpha_{i}^{2}})ds + \int_{0}^{t} \frac{\delta^{Z^{i}}(s)}{\alpha_{i}}.d\overline{W_{s}^{i}}}.
\end{eqnarray*}
Since $c_{t}^{i, *}$ enters the exponential with a negative sign and the remaining coefficients are bounded, $L^{i}$ admits moments of all positive order. Consequently, $Z_{t}^{i} = (L_{t}^{i})^{\alpha_{i}} \in \Hr_{\loc}^{2q}$, for any $q \geq 1$. This yields finiteness of \eqref{eq:martingale_Holder}, which in turn ensures that, under the optimal control, the relative performance criterion $Q^{i}(t, \widehat{X_{t}}^{i, *})$ is a true martingale.

\end{preuve}

\begin{rmq}
\begin{itemize}
\item Boundedness of $\pi_{t}^{i, *}$ is not essential to establish existence and admissibility of the best-response allocation strategy. One can still obtain well-posedness and admissibility under the weaker assumption $\pi_{t}^{i, *} \in \Hr_{\loc}^{2}$. However, proving that the optimal criterion is a true martingale requires stronger moment estimates for the stochastic integral term, and local boundedness is a convenient sufficient assumption.

\item Existence of a solution to \eqref{eq:linearizedY} only requires pathwise integrability of the aggregated consumption factor $(\overset{\sim}{c_{t}}^{(-i)})^{\frac{\theta_{i} \PAR{1 - \alpha_{i}}}{\alpha_{i}}}$. By contrast, for admissibility of the resulting optimal consumption process $c_t^{i,*}$, one needs integrability assumptions on moments of the aggregated consumption factor with the corresponding exponent depending on whether $\alpha_i<1$ or $\alpha_i>1$.

\item The integrability conditions ensuring that the optimal criterion is a true martingale are driven by the consumption-dependent term in the optimal relative wealth dynamics \eqref{xhatdyn}. Similar regime-dependent requirements also appear in classical expected-utility maximization; see e.g. \cite{cheridito2011optimal}.

\end{itemize}
\end{rmq}

\section{The $n$-agent forward optimization problem} \label{sect:nagent}

Having characterized best-response forward relative performance criteria and the associated optimal strategies, we now investigate the existence of a Nash equilibrium in the $n$-player asset specialization problem with common noise. Informally, a Nash equilibrium is a set of $n$ strategies such that no agent can improve her performance being the only one deviating from the equilibrium. 

\subsection{Nash equilibrium for the $n$-player case under forward relative performance criteria}

In this section, we introduce Nash equilibrium under forward relative performance criteria. In the classical expected utility theory, a Nash equilibrium is characterized as $n$ strategies $(\pi_{t}^{i, *}, c_{t}^{i, *})_{i=1, ..., n}$ such that no manager can increase the expected utility of her performance metric by unilateral decision. In the forward setting, this notion is encoded through martingale optimality, as in \cite{dos2022forward}, \cite{zariphop24PUQR}.

\begin{defi}
Let $Q^{i}$ be a progressively measurable criterion generated by utility random fields $U^{i}$ and $V^{i}$, of the form
\begin{eqnarray*}
Q^{i}(t, \widehat{x}^{i}) = U^{i}(t, \widehat{x}^{i}) + \int_{0}^{t} V^{i}(s, \widehat{c_{s}}^{i} \widehat{x}^{i})ds.
\end{eqnarray*}
A forward Nash equilibrium consists in $n$-triples $(Q^{i}, \pi^{i, *}, c^{i, *})$ with $(\pi^{i, *}, c^{i, *}) \in \A$, for $i = 1, ..., n$, such that for each $i$:
\begin{itemize}
\item \textit{Martingale optimality -} Let managers $j \neq i$ act according $(\pi^{j, *}, c^{j, *})$. Then for any admissible $(\pi^{i}, c^{i})$, the process $Q^{i}(t, \widehat{X_{t}}^{i})$ is a local supermartingale. Let all managers $j=1, ..., n$ act along their optimal strategy $(\pi^{j, *}, c^{j, *})$. Then $Q^{i}(t, \widehat{X_{t}}^{i, *})$ is a local martingale.
\end{itemize}
\end{defi}

This definition is consistent with the classical notion of Nash equilibrium. Indeed, for any $i \in \BRA{1, ..., n}$, if all agents $j \neq i$ use their equilibrium strategies $\pi_{t}^{j, *}$, and if $Q^{i}(t, \widehat{X_{t}}^{i, *})$ is a true martingale while $Q^{i}(t, \widehat{X_{t}}^{i})$ is a true supermartingale for every admissible deviation of agent $i$, then
\begin{eqnarray*}
\e \SBRA{Q^{i}(t, \widehat{X_{t}}^{i, *})} = Q^{i}(0, \widehat{X_{0}^{i}}) \geq \e \SBRA{Q^{i}(t, \widehat{X_{t}}^{i})}, \quad \text{for all } t \geq 0.
\end{eqnarray*}
Hence no manager can improve her expected performance by unilateral deviation from the equilibrium. A corresponding statement holds in localized form when the optimal criterion is only a local martingale.

The characterization of Nash equilibria in the general case is challenging, especially for consumption. Consider a family of random utility field pairs $(U^{i}, V^{i})_{\BRA{i = 1, \dots, n}}$, each satisfying the consistency SPDE \eqref{utspde_first}. For given preferences, Nash equilibrium strategies are typically obtained through a fixed point argument based on the best-response formulas \eqref{opt_port} and \eqref{optcons}, in order to characterize mean equilibrium quantities , defined for $x = (x_{1}, \dots, x_{n}) \in \R^{n}$ by
\begin{eqnarray}
    \overline{\sigma \pi_{t}}(x) = \frac{1}{n} \sum_{k=1}^{n} \sigma_{k} \pi_{t}^{k}(x_{k}), \quad \tilde{c_{t}}(x) = \PAR{\prod_{k=1}^{n} c_{t}^{k}(x_{k})}^{\frac{1}{n}}
\end{eqnarray}
Define $D^{i}(t,x) = \sigma_{i} \theta_{i} \PAR{1 + \frac{U_{x}^{i}(t, x_{i})}{U_{xx}^{i}(t, x_{i})x_{i}}}$ for $i=1, \dots, n$. Then, the best-response portfolio \eqref{opt_port} takes the form 
\begin{eqnarray} \label{eq:portNashgen}
\pi_{t}^{i, *}(x_{i}) = \frac{1}{\nu_{i}^{2} + \sigma_{i}^{2} + D^{i}(t,x_{i})} \PAR{ D^{i}(t,x_{i}) \overline{\sigma \pi_{t}}(x) - \frac{1}{U_{xx}^{i}(t, x_{i})x_{i}} \PAR{\gamma_{x}^{i, W}(t, x_{i}) \nu_{i} + \gamma_{x}^{i, B}(t, x_{i}) \sigma_{i} + \mu_{i} U_{x}^{i}(t, x_{i})}}.
\end{eqnarray}
For $x = (x_{1}, \dots, x_{n}) \in \R^{n}$, let 
\begin{eqnarray*}
\psi_{n}^{U}(t,x) &=& \frac{1}{n} \sum_{k=1}^{n} \frac{\sigma_{k} D^{k}(t,x_{k})}{\nu_{k}^{2} + \sigma_{k}^{2} + D^{k}(t, x_{k})}, \\
\varphi_{n}^{U}(t,x) &=& -\frac{1}{n} \sum_{k=1}^{n}  \frac{\sigma_{k}}{\nu_{k}^{2} + \sigma_{k}^{2} + D^{k}(t,x_{k})} \frac{1}{U_{xx}^{k}(t, x_{k})x_{k}} \PAR{\gamma_{x}^{k, W}(t, x_{k}) \nu_{k} + \gamma_{x}^{k, B}(t, x_{k}) \sigma_{k} + \mu_{k} U_{x}^{k}(t, x_{k})}.
\end{eqnarray*}
Multiplying \eqref{eq:portNashgen} by $\sigma_{i}$ and averaging over $i$ yields 
\begin{eqnarray*}
    \overline{\sigma \pi_{t}}(x) = \psi_{n}^{U}(t, x) \overline{\sigma \pi_{t}}(x) + \varphi_{n}^{U}(t, x).
\end{eqnarray*}
Hence, if $\psi_{n}^{U}(t, x) \neq 1$ for all $x \in (\R^{+})^{n}$, $t \geq 0$, there exists a candidate Nash equilibrium portfolio, given by \eqref{eq:portNashgen}, where the aggregate portfolio rescaled by common noise weights is given by $\overline{\sigma \pi_{t}}(x) = \frac{\varphi_{n}^{U}(t, x)}{1 - \psi_{n}^{U}(t,x)}.$

Identifying equilibrium consumption is substantially more difficult. In \eqref{optcons} the competitors aggregate consumption $\overset{\sim}{c_{t}}^{(-i)}$ appears inside the inverse marginal utility from consumption, so the fixed point equation cannot be closed explicitly without further assumptions on $V^{i}$. We postpone this issue to Section \ref{sect:MF}, where analogous assumptions are introduced in a more tractable mean-field setting.

Finally, even after a candidate $(\overline{\sigma \pi_{t}}, \tilde{c_{t}})$ is obtained, proving admissibility of the resulting equilibrium remains delicate. The reason is the strong cross agent coupling generated by the consistency SPDE \eqref{utspde_first}, which propagates to both portfolio and consumption dynamics. In the following, we focus on separable CRRA utilities and show that admissibility of equilibrium consumption is already nontrivial in that setting.

\subsection{Nash equilibrium under forward relative performance criteria of CRRA type}

In this section, we study the $n$-player game for relative performance criteria generated by CRRA utility random fields with separable time and space dependence. Consider, for $i=1, ..., n$, utility random fields of the form 
\begin{eqnarray}
U^{(\alpha_{i})}(t, x) = Z_{t}^{i} u^{(\alpha_{i})}(x), \quad \text{and} \quad V^{(\alpha_{i})}(t, c) = \phi_{t}^{i} u^{(\alpha_{i})}(c), \quad \alpha_{i} \in (0,1) \cup (1, +\infty).
\end{eqnarray}

\vspace{0.3cm}
\noindent
To state the result, define
\begin{eqnarray} \label{eq:defpsi_nphi_n}
    \psi_{n} = \frac{1}{n} \sum_{k=1}^{n} \frac{\sigma_{k}^{2} \theta_{k}(1 - \frac{1}{\alpha_{k}})}{\nu_{k}^{2} + \sigma_{k}^{2}\PAR{1 + \theta_{k}(1 - \frac{1}{\alpha_{k}})}}, \quad \varphi_{n}(t) = \frac{1}{n} \sum_{k=1}^{n} \frac{\frac{\sigma_{k}}{\alpha_{k}} \PAR{\delta^{Z^{k}, W}(t) \nu_{k} + \delta^{Z^{k}, B}(t) \sigma_{k} + \mu_{k}}}{\nu_{k}^{2} + \sigma_{k}^{2} \PAR{1 + \theta_{k}(1 - \frac{1}{\alpha_{k}})}}.
\end{eqnarray}
\begin{eqnarray} \label{eq:defxi}
\xi_{n} = \sum_{k=1}^{n} \frac{\theta_{k}(1 - \alpha_{k}) }{(n-1) \alpha_{k} - (1 - \alpha_{k}) \theta_{k}}.
\end{eqnarray}
The next result establishes existence of a forward Nash equilibrium associated with the family $(U^{(\alpha_{i})}, V^{(\alpha_{i})})_{1 \leq i \leq n}$, and gives tractable admissibility regimes. We assume that each pair $(U^{(\alpha_{i})}, V^{(\alpha_{i})})$ solve the HJB-SPDE \eqref{utspde_first}, so that for all $i$, $Z^{i}$ satisfies \eqref{driftcond1}. Because consumption coupling enters both optimality and drift consistency, admissibility requires additional structure on $b^{\phi^{i}}$, and we identify two regimes depending on risk aversions.

\begin{thm} \label{thm:Nash}
Let $(U^{(\alpha_{i})}, V^{(\alpha_{i})})_{1 \leq i \leq n}$ be a family of pairs of utility random fields, and assume each pair $(U^{(\alpha_{i})}, V^{(\alpha_{i})})$ satisfies the HJB-SPDE \eqref{utspde_first}. If $\psi_{n} \neq 1$, then there exists a unique candidate Nash equilibrium strategy $(\pi_{t}^{i, *}, c_{t}^{i, *})$ given by
\begin{align} 
\pi_{t}^{i, *} &= \frac{\sigma_{i}}{\nu_{i}^{2} + \sigma_{i}^{2}\PAR{1 + \theta_{i}(1 - \frac{1}{\alpha_{i}})}} \SBRA{\sigma_{i} \theta_{i}(1 - \frac{1}{\alpha_{i}}) \frac{n}{n-1} \frac{\varphi_{n}(t)}{1 - \psi_{n}} + \frac{1}{\alpha_{i}} \PAR{\delta^{Z^{i}, W}(t) \nu_{i} + \delta^{Z^{i}, B}(t) \sigma_{i} + \mu_{i}}} \label{eqport} \\
c_{t}^{i, *} &= \PAR{\tilde{c_{t}}}^{- \frac{n \theta_{i}(1 - \alpha_{i})}{\alpha_{i}(n-1) - \theta_{i}(1 - \alpha_{i})}} \PAR{\frac{\phi_{t}^{i}}{Z_{t}^{i}}}^{\frac{n-1}{\alpha_{i}(n-1) - \theta_{i}(1 - \alpha_{i})}}, \label{eqcons}
\end{align}
where the geometric average of consumption rate processes
\begin{eqnarray} \label{eq:aggregatecompconsN}
\tilde{c_{t}} = \PAR{\prod_{k=1}^{n} \PAR{\frac{\phi_{t}^{k}}{Z_{t}^{k}}}^{\frac{1}{\alpha_{k}(n-1) - \theta_{k}(1 - \alpha_{k})}}}^{\frac{n-1}{n} \frac{1}{\xi_{n} + 1}}.
\end{eqnarray}
Assume in addition that for all agent $i \in \BRA{1, \dots, n}$, the drift of $\phi^{i}$ is of the form 
\begin{eqnarray} \label{eq:condphiNash}
    b^{\phi^{i}}(t, Z_{t}^{i}, \phi_{t}^{i}) = \overline{b^{\phi^{i}}}(t) - \lambda_{i} \frac{\alpha_{i}}{(\overset{\sim}{c_{t}}^{(-i)})^{\frac{\theta_{i}(1 - \alpha_{i})}{\alpha_{i}}}} \PAR{\frac{\phi_{t}^{i}}{Z_{t}^{i}}}^{\frac{1}{\alpha_{i}}},
\end{eqnarray}
where $\overline{b^{\phi^{i}}} \in L_{\loc}^{\infty}$, $\lambda_{i} \in \R$ and $\alpha_{i} > \frac{\theta_{i}}{n-1 + \theta_{i}}$. Then,
\begin{itemize}
    \item if $\alpha_{i} \in (0, 1)$ for all $i$, and $\lambda_{i} = 1$ for all $i$, then $(\pi_{t}^{i, *}, c_{t}^{i, *})_{1 \leq i \leq n}$ is an admissible forward Nash equilibrium.
    \item if $\alpha_{i} \in (1, \infty)$,  and $\lambda_{i} \geq 1$ for all $i$ (not necessarily equal across agents), then $(\pi_{t}^{i, *}, c_{t}^{i, *})_{1 \leq i \leq n}$ is an admissible forward Nash equilibrium.
\end{itemize}
\end{thm}

\begin{preuve}
The equilibrium portfolio is obtained via a fixed-point argument from \eqref{opt_port}. The equilibrium consumption component is derived similarly, but requires additional care because consumption dependent terms also enter the drift constraint \eqref{driftcond1}.

\vspace{0.3cm}
\noindent
\textbf{Equilibrium portfolio -} For CRRA separable wealth utilities $U^{(\alpha_{i})}(t, x) = Z_{t}^{i} u^{(\alpha_{i})}(x)$, the best-response portfolio \eqref{opt_port} rescaled by the common noise volatility weight $\sigma_{i}$ becomes
\begin{small}
\begin{eqnarray} \label{pi_gen_equil}
\sigma_{i} \pi_{t}^{i, *} = \frac{\sigma_{i}}{\nu_{i}^{2} + \sigma_{i}^{2}\PAR{1 + \theta_{i}(1 - \frac{1}{\alpha_{i}})}} \SBRA{\sigma_{i} \theta_{i}(1 - \frac{1}{\alpha_{i}}) \frac{n}{n-1} \overline{\sigma \pi_{t}} + \frac{1}{\alpha_{i}} \PAR{\delta^{Z^{i}, W}(t) \nu_{i} + \delta^{Z^{i}, B}(t) \sigma_{i} + \mu_{i}}}.
\end{eqnarray}
\end{small}
Averaging over $i \in \BRA{1, ..., n}$, we obtain
\begin{eqnarray*}
\overline{\sigma \pi_{t}} = \overline{\sigma \pi_{t}} \psi_{n} + \varphi_{n}(t),
\end{eqnarray*}
where $\psi_{n}$ and $\phi_{n}$ given by \eqref{eq:defpsi_nphi_n}. Hence, if $\psi_{n} \neq 1$
\begin{eqnarray*}
    \overline{\sigma \pi_{t}} = \frac{\varphi_{n}(t)}{1 - \psi_{n}}.
\end{eqnarray*}
Plugging this expression into \eqref{pi_gen_equil} yields the Nash equilibrium portfolio strategy.
If $\psi_{n} = 1$, the fixed-point equation is degenerate and there exists no Nash equilibrium.

\vspace{0.3cm}
\noindent
\textbf{Equilibrium consumption process -} The best-response consumption strategy \eqref{eq:optconsCRRAn} can be expressed as a function of $\tilde{c_{t}} = \PAR{\prod_{k=1}^{n} c_{t}^{k}}^{\frac{1}{n}}$ as
\begin{eqnarray} \label{optconsnagentp}
c_{t}^{i, *} = \PAR{\tilde{c_{t}}}^{- \frac{n \theta_{i}(1 - \alpha_{i})}{\alpha_{i}(n-1) - \theta_{i}(1 - \alpha_{i})}} \PAR{\frac{\phi_{t}^{i}}{Z_{t}^{i}}}^{\frac{n-1}{\alpha_{i}(n-1) - \theta_{i}(1 - \alpha_{i})}}.
\end{eqnarray}
Averaging over $i = 1, ..., n$, we obtain
\begin{eqnarray*}
\prod_{k=1}^{n} c_{t}^{k, *} = \prod_{k=1}^{n} \PAR{\frac{\phi_{t}^{k}}{Z_{t}^{k}}}^{\frac{n-1}{\alpha_{k}(n-1) - \theta_{k}(1 - \alpha_{k})}} \PAR{\tilde{c_{t}}}^{- n\sum_{k=1}^{n} \frac{\theta_{k}(1 - \alpha_{k})}{\alpha_{k}(n-1) - \theta_{k}(1 - \alpha_{k})}}.
\end{eqnarray*}
Using $\xi_{n}$ from \eqref{eq:defxi}, the geometric average of consumption processes $\tilde{c_{t}}$ then takes the form \eqref{eq:aggregatecompconsN}. Substituting in the best-response consumption \eqref{optconsnagentp} gives \eqref{eqcons}.

\noindent
\textbf{Admissibility -} Establishing well-posedness of the equilibrium consumption \eqref{eqcons} is significantly more delicate than in the best-response problem of Section \ref{sect:consutCRRA}. Indeed, under the consistency condition \eqref{driftcond1}, one can write the $n$-dimensional SDE for $Z$ and prove existence and uniqueness of a strong global solution. The main difficulty lies in proving that the resulting processes remain positive at all times, by computing directly the dynamics of the optimal consumption process. 
Let $R_{t}^{k} = \frac{\phi_{t}^{k}}{Z_{t}^{k}}$ for $k \in \BRA{1, \dots, n}$. Then the candidate equilibrium consumption \eqref{eqcons} can be written as
\begin{eqnarray*}
c_{t}^{i, *} = \prod_{k=1}^{n} (R_{t}^{k})^{\eta_{ik}},
\end{eqnarray*}
with exponents
\begin{align} \label{eq:eta}
\eta_{ik} = -\frac{(n-1)\,\theta_i(1-\alpha_i)}{(\xi_n+1)\,A_i\,A_k},
\quad \text{for } k\neq i, \quad \text{and } \eta_{ii} = \frac{n-1}{A_i}
\left(
1-\frac{\theta_i(1-\alpha_i)}{(\xi_n+1)A_i}
\right), 
\end{align}
where for $k=1,\dots,n$,
\[
D_k := \alpha_k (n-1) - \theta_k(1-\alpha_k) > 0, \quad \text{for } \alpha_{k} > \frac{\theta_{k}}{n-1+\theta_{k}}.
\]
Assuming that \eqref{eq:condphiNash} holds for all agents, and applying It\^o's formula, the equilibrium consumption vector $c_{t}^{*} = (c_t^{1,*},\dots,c_t^{n,*})$ is an $n$-dimensional Ito diffusion of the form
\begin{eqnarray*}
    dc_{t}^{*} = \beta^{c}(t, c_{t}^{*})dt + \gamma^{c}(t, c_{t}^{*}).d\overline{W_{t}},
\end{eqnarray*}
where $\overline{W_{t}}$ denotes the $n+1$ dimensional Brownian motion $(W_{t}^{1}, \dots, W_{t}^{n}, B_{t})$ and local characteristics $\beta^{c}(t, c) = \PAR{\beta^{c}_{1}(t, c), \dots, \beta^{c}_{n}(t, c)}^{\top}$ and $\gamma^{c}(t, c) \in \R^{n, n+1}$ are for $i=1, \dots, n$,
\begin{eqnarray*}
\beta_{i}^{c}(t, c) &=& c_{i}\PAR{\overline{\beta_{i}^{c}}(t) -  \sum_{k=1}^{n} (\lambda_{k} - 1) \eta_{ik} \alpha_{k} c_{k} + \frac{1}{2} \NRM{\delta^{c}_{i}(t)}^{2}} \\
\gamma^{c}_{i, .}(t, c) &=& c_{i} \delta^{c}_{i}(t)^{\top}.
\end{eqnarray*}
Here, $\overline{\beta_{i}^{c}}$ depends only on the exogenous coefficients $\overline{b^{Z^{k}}}, \overline{b^{\phi^{k}}}, \delta^{Z^{k}}$, $\delta^{\phi^{k}}$, while
\[
\delta^{c}_{i}(t)=\sum_{k=1}^n \eta_{ik}\big(\delta^{\phi^k}(t)-\delta^{Z^k}(t)\big).
\]
When attempting to apply Theorem 1.4 from \cite{lan2014new}, one obtains
\begin{eqnarray} \label{eq:dissipN}
\CRO{c, \beta^{c}(t, c)} = \sum_{i=1}^{n} c_{i}^{2} \PAR{\overline{\beta_{i}^{c}}(t) + \frac{1}{2} \NRM{\delta^{c}_{i}(t)}^{2}} - \sum_{i=1}^{n} \sum_{k=1}^{n} (\lambda_{k} - 1) \eta_{ik} \alpha_{k} c_{i}^{2} c_{k}.
\end{eqnarray}
The difficulty stems from the cubic interaction term in the second sum. Indeed,
\begin{itemize}
\item for $\alpha_{i} \in (0, 1)$, $\eta_{ik} < 0$ for all $k \neq i$, while $\eta_{ii} > 0$. The interaction has mixed signs and is not uniformly controllable in general. A tractable case is $\lambda_{k} = 1$ for all $k=1, \dots, n$, where interaction vanishes and the Nash equilibrium exists and is admissible.
\item for $\alpha_{i} \in (1, +\infty)$, and $1 + \xi_{n} > 0$, then $\eta_{ik} > 0$ for all $i, k$. Therefore, dissipativity of \eqref{eq:dissipN} holds for $\lambda_{k} \geq 1$, not necessarily equal across $k$, yielding existence of an admissible Nash equilibrium.
\end{itemize}
\end{preuve}

\begin{rmq}
Existence of an admissible Nash equilibrium is substantially more delicate when risk aversion regimes are heterogeneous across the population (e.g., some $\alpha_i<1$ and others $\alpha_j>1$). In that case, the coefficients in \eqref{eq:dissipN} generally have mixed signs, and dissipativity may fail. This motivates the mean-field formulation, in which population effects are captured through the conditional law of the state given the common noise.
\end{rmq}

\section{Mean-field forward optimization problem} \label{sect:MF}

Let $(B_{t})_{t \geq 0}$ and $(W_{t})_{t \geq 0}$ be two independent Brownian motions, living on a filtered probability space $\PAR{\Omega, \F, \f = \PAR{\F_{t}}_{t \geq 0}, \p}$. In order to model the continuum of agent of the mean-field optimization problem, we introduce a random type vector $\zeta$, independent of $B$ and $W$, whose distribution describes the proportion of the population following the corresponding preferences. We study a mean-field optimization problem under forward relative performance criterion in the It\^o setting. In particular, consider utility random fields of wealth and consumption $U$ and $V$ of the form
\begin{align*}
dU(t, x) &= \beta^{U}(t,x) dt + \gamma^{U, W}(t, x)dW_{t} + \gamma^{U, B}(t, x)dB_{t} \\
dV(t, c) &= \beta^{V}(t, c)dt + \gamma^{V, W}(t, c)dW_{t} + \gamma^{V, B}(t, c)dB_{t}.
\end{align*}
where $\beta^{U}, \beta^{V} \in \K^{m, \delta}_{\loc} $ and $ \gamma^{U, W}, \gamma^{U, B}, \gamma^{V, W}, \gamma^{V, B} \in \overline{\K}^{m, \delta}_{\loc}$. Compared with the random type vector specification in \cite{lacker2020many}, \cite{dos2022forward}, we distinguish static an dynamic preference heterogeneity in the forward setting with non-vanishing volatility. First, we define the typical random type vector collecting time independent preferences, including initial utility specifications, the competition parameter and market parameters, as 
\begin{eqnarray*}
    \zeta = \PAR{u_{0}, v_{0}, x_{0}, \theta, \mu, \nu, \sigma},
\end{eqnarray*}
with values on the type space
\begin{eqnarray*}
\Z = \PAR{\U_{std} \times \U_{std} \times (0, \infty) \times [0,1] \times \R \times \R^{+} \times \R^{+}}.
\end{eqnarray*}
Let $(\F_{t}^{MF})_{t \geq 0}$ be the smallest filtration such that $\zeta$ is $\F_{0}^{MF}$-measurable and $B$, $W$ are adapted, and let $(\F_{t}^{B})_{t \geq 0}$ denote the natural filtration generated by $B$. For general references on mean-field games with common noise, see \cite{carmona2016}, \cite{carmona2018probabilistic}.

To represent heterogeneity in time-varying preference inputs, we introduce a dynamic preference factor $\xi$, that includes the wealth utility volatility and the consumption utility local characteristics, namely
\begin{eqnarray*}
\xi_{t} = \PAR{\gamma_{t}^{U}, \beta_{t}^{V}, \gamma_{t}^{V}},
\end{eqnarray*}
where $\gamma_{t}^{U} = \PAR{\gamma^{U, W}(t, .), \gamma^{U, B}(t, .)}$ and $\gamma_{t}^{V} = \PAR{\gamma^{V, W}(t, .), \gamma^{V, B}(t, .)}$. Once $\gamma^{U}$ and the local characteristics of $V$ are fixed, the consistency SPDE determines the compatible drift $\beta^{U}$. At each time $t$, the law of $\xi_{t}$ describes the distribution of dynamic preference factors across the population. At this stage, $\xi$ is only an a priori formulation of dynamic preference heterogeneity: not every trajectory of $\xi_{t}$ yields a random utility pair $(U, V)$ defining a time-consistent forward criterion. The main challenge is to characterize the admissible subset of solution to the consistency SPDE, so that each admissible trajectory of $\xi_{t}$ induces an admissible forward criterion. We will do so in the general case under additional assumptions in Section \ref{sect:MF31}, and explicitly for separable CRRA preferences in Section \ref{sect:MF_CRRA}.

\subsection{Mean-field Nash equilibrium under forward relative performance criteria} \label{sect:MF31}

In this section, we first define the notion of mean-field Nash equilibrium under forward relative performance criteria, and then state our main equilibrium result. Compared with the best-response analysis of Section \ref{sect:bestresponsegen}, stronger structural and integrability assumptions are required. The reason is twofold: the optimal relative wealth dynamics are of McKean-Vlasov type with common noise, and the compatibility condition involves geometric population aggregates that need to be computed explicitly.

Under the self-financing condition, the generic agent's discounted wealth process solves
\begin{eqnarray} \label{wealthMF}
dX_{t}^{\pi, c} = \pi_{t} X_{t}^{\pi, c} \PAR{\mu dt + \nu dW_{t} + \sigma dB_{t}} - c_{t} X_{t}^{\pi, c} dt, \quad X_{0}^{\pi, c} = x_{0}, 
\end{eqnarray}
where $\pi_{t}$ stands for the fraction of wealth invested in the risky asset and $c_{t}$ is the rate of consumption per unit of wealth. The set of admissible strategies is defined as
\begin{multline}
    \A^{MF} = \left\{ (\pi_t, c_t) \, \F^{MF}-\text{progressively measurable process valued in $\R \times \R^{+, *}$,} \right. \\
    \left. \text{such that } \forall t \geq 0, \, \e \SBRA{\int_{0}^{t} (\ABS{\pi_{s}}^{2} + \ABS{c_{s}})ds} < \infty \right\}.
\end{multline}

In the following, we define mean-field forward relative performance criteria, that are evaluated along a relative performance metric. Let $(\overline{X_{t}})_{t \geq 0}$ and $(\overline{C_{t}})_{t \geq 0}$ be $\F_{t}^{B}$-adapted processes representing the geometric average wealth and the geometric average consumption of the continuum of agents. In the mean-field setting, the equilibrium strategy $(\pi^{*}, c^{*})$ should be typical of the population. In turn, population relative average wealth and consumption process at equilibrium are characterized by a fixed point relation, motivated by the conditional law of large numbers in the presence of common noise, see \cite{lacker2020many}, \cite{dos2022forward}. Heuristically, consider a finite-population approximation of the mean-field model in which each agent $k$ follows the same strategy $(\pi^{*}, c^{*})$, so that $X^{k}$ solves \eqref{wealthMF}. The empirical geometric mean, analogous to \eqref{relativewealth}, is
\begin{eqnarray*}
\overline{X_{t}}^{(n)} = \exp \PAR{\frac{1}{n} \sum_{k=1}^{n} \log X_{t}^{k}}.
\end{eqnarray*}
Hence, for each fixed $t$, conditionally on the common noise filtration $\F_{t}^{B}$, all agents face i.i.d copies of the same control problem. Therefore the random variables $(\log X_{t}^{k})_{k \geq 1}$ are $i.i.d$, and denoting $X_{t}^{*}$ the optimal wealth process resulting from $(\pi^{*}, c^{*})$, the conditional law of large numbers yields
\begin{eqnarray*}
\frac{1}{n} \sum_{k=1}^{n} \log X_{t}^{k} \underset{n \to \infty}{\to} \e \SBRA{\log X_{t}^{*} | \F_{t}^{B}}, \quad \text{a.s}.
\end{eqnarray*}
Exponentiating gives $\overline{X_{t}} = \exp \PAR{\e \SBRA{\log X_{t}^{*} | \F_{t}^{B}}}$, which characterizes relative geometric average wealth. The geometric average consumption term is characterized in the same way. Time consistency of the preference criterion is then enforced by the martingale optimality principle.

\begin{defi} \label{MFequi}
Let $(\overline{X_{t}})_{t \geq 0}$ and $(\overline{C_{t}})_{t \geq 0}$ be $\F^{B}$-adapted processes representing the geometric average wealth and the geometric average consumption of the continuum of agents. 
For an admissible strategy $(\pi_t, c_t)$, let $X^{\pi, c}$ solve \eqref{wealthMF}. Consider an $\F^{MF}$-progressively measurable random field $Q$ of the form
\begin{eqnarray}\label{eq:MF_frpp}
Q(t, \widehat{x_{t}}) = U(t, \widehat{x_{t}}) + \int_{0}^{t} V(s, \widehat{c_{s}} \widehat{x_{t}})ds, 
\end{eqnarray}
with relative metric given by $\widehat{c_{t}} = \displaystyle\frac{c_{t}}{(\overline{C_{t}})^{\theta}},$  $\widehat{x_{t}} = \displaystyle\frac{x_{t}}{(\overline{X_{t}})^{\theta}}$, where $U$, $V$ are utility random fields from wealth and consumption respectively. $Q$ is a mean-field forward relative performance criterion if
\begin{enumerate}
    \item \textit{Martingale optimality -} For any admissible strategy $(\pi_t, c_t)$, $Q\PAR{t, \displaystyle\frac{X_{t}^{\pi, c}}{(\overline{X_{t}})^{\theta}}}$ is a local supermartingale and there exists an admissible strategy $(\pi_{t}^{*}, c_{t}^{*})$, such that $Q \PAR{t, \displaystyle\frac{X_{t}^{\pi^{*}, c^{*}}}{(\overline{X_{t}})^{\theta}}}$ is a local martingale.
    \item \textit{Compatibility -} The strategy $(\pi_{t}^{*}, c_{t}^{*}) \in \A^{MF}$ is a mean-field Nash equilibrium if
    \begin{eqnarray} \label{compatcond}
        \ACC{\overline{X_{t}} &= \exp \PAR{\e \SBRA{\log X_{t}^{*} | \F_{t}^{B}}} \\
        \overline{C_{t}} &= \exp \PAR{\e \SBRA{\log c_{t}^{*} | \F_{t}^{B}}}}.
    \end{eqnarray}
    We call the equilibrium \textit{strong} if $(\pi_{t}^{*}, c_{t}^{*})$ is $\F_{0}^{MF}$-measurable.
    \end{enumerate}
\end{defi}

Additional regularity assumptions, analogous to Assumptions \ref{ass:admport} and \ref{ass:admcons}, are needed to ensure admissibility of the resulting mean-field equilibrium strategy. These conditions are stronger than in the best-response setting, because the optimal relative wealth process $\widehat{X_{t}^{*}}$ solves a McKean-Vlasov SDE with common noise. The first assumption requires a linear form for the wealth utility derivative ratios arising in the equilibrium portfolio.

\begin{ass} \label{ass:MFadm}
There exists $\F^{MF}$-adapted processes $A_{t}^{1}, \, A_{t}^{2}, \, A_{t}^{3} \in \Hr_{\loc}^{2}$ such that for all $t \geq 0$, $x > 0$:
\begin{eqnarray}
\frac{\gamma_{x}^{U, W}(t, x)}{U_{xx}(t, x)} = A_{t}^{1}x, \quad \frac{\gamma_{x}^{U, B}(t, x)}{U_{xx}(t, x)} = A_{t}^{2}x, \quad \frac{U_{x}(t, x)}{U_{xx}(t, x)} = A_{t}^{3}x.
\end{eqnarray}
We also require that a.s for all $t \geq 0$, $\e \SBRA{\frac{\theta \sigma^{2}}{\nu^{2} + \sigma^{2}}} + \e \SBRA{\theta \sigma^{2} A_{t}^{3} | \F_{t}^{B}} - 1 \neq 0.$
\end{ass}

Secondly, one structural assumption linking the compound of marginal utilities from consumption and wealth is required for the compatibility condition \eqref{compatcond} on consumption to be explicitly solvable. We formulate it as a state independence property of the mean-field equilibrium consumption, which holds for separable CRRA utilities studied in the next section.

\begin{ass} \label{ass:structconsMF}
There exist an $\F^{MF}$-adapted process $K_{t} \in \Hr_{\loc}^{1}$ and an $\F_{0}^{MF}$-measurable real valued random variable $p$ such that a.s for all $t \geq 0$, $x > 0$:
\begin{eqnarray*}
V_{x}^{-1}(t, U_{x}(t, x)(\overline{C_{t}})^{\theta}) = K_{t} x (\overline{C_{t}})^{p}.
\end{eqnarray*}
Additionally, the random bound $K_{t}$ is required to satisfy for all $t \geq 0$,
\begin{eqnarray} \label{eq:admconsMF}
    \e \SBRA{\int_{0}^{t} K_{t} \exp \PAR{\e \SBRA{\log(K_{t}) | \F_{t}^{B}}}^{{\frac{p+ \theta}{1 - \e \SBRA{p+\theta}}}}} < \infty.
\end{eqnarray}
\end{ass}
Although somewhat technical, condition \eqref{eq:admconsMF} is precisely the integrability required to guarantee admissibility of the resulting mean-field equilibrium consumption. For later use, define 

\begin{multline} \label{eq:defpsiMF}
\psi(t, x) = U_{x}(t, x) x \PAR{- \theta \overline{\mu \pi_{t}} + \frac{\theta}{2} \overline{\Sigma \pi_{t}^{2}} + \frac{\theta^{2}}{2} \overline{\sigma \pi_{t}}^{2} + \theta \overline{c_{t}}} + \frac{1}{2}U_{xx}(t, x) x^{2} (\theta \overline{\sigma \pi_{t}})^{2} - \gamma_{x}^{U, B}(t, x) \theta \overline{\sigma \pi_{t}} x.
\end{multline}
where bar quantities denote conditional expectations with respect to the common noise filtration, namely 
\begin{eqnarray}
\overline{\mu \pi_{t}} = \e \SBRA{\mu \pi_{t} | \F_{t}^{B}}, \quad \overline{\sigma \pi_{t}} = \e \SBRA{\sigma \pi_{t} | \F_{t}^{B}}, \quad \overline{\Sigma \pi_{t}^{2}} = \e \SBRA{(\nu^{2} + \sigma^{2})\pi_{t}^{2} | \F_{t}^{B}}, \quad \overline{c_{t}} = \e \SBRA{c_{t} | \F_{t}^{B}}.
\end{eqnarray}
Given the dynamics of the wealth process \eqref{wealthMF}, the compatibility condition \eqref{compatcond} yields the dynamics of the average process $\overline{X_{t}}$. Applying Itô-Ventzel's formula then leads to a mean-field consistency SPDE, analogous to \eqref{utspde_first}. We also derive the candidate optimal strategy in a general setting. 

\begin{thm} \label{thm:MFSPDE}
Let $(U,V)$ be a candidate random utility system satisfying Assumptions \ref{ass:MFadm} and \ref{ass:structconsMF}. Assume that $(U, V)$ solves
\begin{multline} \label{eq:MFSPDE}
dU(t,x) = \PAR{- \psi(t, x) + \frac{1}{2} U_{xx}(t, x) x^{2} (\nu^{2} + \sigma^{2}) (\pi_{t}^{*}(x))^{2} - \overset{\sim}{V}(t, (\overline{C_{t}})^{\theta} U_{x}(t, x))}dt \\
+ \gamma^{U, W}(t, x)dW_{t} + \gamma^{U, B}(t, x)dB_{t},
\end{multline}
where
\begin{multline} \label{portMF}
\pi_{t}^{*}(x) = \frac{1}{\nu^{2} + \sigma^{2}} \left( \theta \sigma \overline{\sigma \pi_{t}} - \frac{1}{U_{xx}(t, x) x} \PAR{\gamma_{x}^{U, W}(t,x) \nu + \gamma_{x}^{U, B}(t, x) \sigma  + (\mu - \theta \sigma \overline{\sigma \pi_{t}})U_{x}(t,x)} \right).
 \end{multline}
 \begin{flalign}
c_{t}^{*}(x) = \frac{V_{x}^{-1} \PAR{t, U_{x}(t, x) (\overline{C_{t}})^{\theta}} (\overline{C_{t}})^{\theta}}{x}. \label{FOCconsMF}
 \end{flalign}
Then $Q$ defined by \eqref{eq:MF_frpp} is a mean-field forward relative performance criterion, and the policy $(\pi_{t}^{*}, c_{t}^{*})$ is a mean-field Nash equilibrium in the sense of Definition \ref{MFequi}. In particular, the aggregate quantities $\overline{\sigma \pi_{t}}$ and $\overline{C_{t}}$ are uniquely determined and compatible with the optimal strategy $(\pi_{t}^{*}, c_{t}^{*}) \in \A^{MF}$.
\end{thm}

\begin{preuve}
The proof is in three steps. First, using the compatibility condition \eqref{compatcond}, we derive the It\^o decomposition of $\overline{X_{t}}.$ Second, we apply It\^o-Ventzel's formula to $Q(t, \frac{X^{\pi, c}}{(\overline{X^{\pi, c}})^{\theta}})$, and analyze the martingale optimality condition of Definition \ref{MFequi}. Third, as in Section \ref{sect:FRPP}, the first-order conditions lead to an explicit form of the mean-field Nash equilibrium.

\vspace{0.3cm}
\noindent
\textbf{Average wealth process -} Let $(\pi_{t}, c_{t}) \in \F_{t}^{MF}$. By \eqref{compatcond}, we restrict to processes $(\overline{X_{t}}^{\pi, c})_{t \geq 0}$ satisfying $\overline{X_{t}}^{\pi, c} = \exp \e \SBRA{\log X_{t}^{\pi, c} | \F_{t}^{B}}$. Applying It\^o's formula
\begin{align}
\overline{X_{t}}^{\pi, c} &= \exp \e \SBRA{\log X_{t}^{\pi, c} | \F_{t}^{B}} \nonumber \\
    &= \exp \e \left[ \log x_{0} + \int_{0}^{t}(\mu \pi_{s} - \frac{1}{2} \pi_{s}^{2}(\nu^{2} + \sigma^{2}))ds + \int_{0}^{t} \nu \pi_{s}dW_{s} \right. \nonumber \\
    & \pushright{ \left. \hfill + \int_{0}^{t} \sigma \pi_{s}dB_{s} - \int_{0}^{t} c_{s}ds |\F_{t}^{B} \right] }\label{espbruitcomm} \\
    &= \exp \PAR{\overline{\log x_{0}} + \int_{0}^{t} (\overline{\mu \pi_{s}} - \frac{1}{2} \overline{\Sigma \pi_{s}^{2}})ds + \int_{0}^{t} \overline{\sigma \pi_{s}}dB_{s} - \int_{0}^{t} \overline{c_{s}}ds} \label{espint} \\
    &= \overline{x_{0}} + \int_{0}^{t} \eta \overline{X_{s}}ds + \int_{0}^{t} \overline{\sigma \pi_{s}} \overline{X_{s}}^{\pi, c}dB_{s} - \int_{0}^{t} \overline{c_{s}}\overline{X_{s}}^{\pi, c}ds, \label{Xbardyn}
\end{align}
where
\begin{eqnarray*}
\overline{\eta_{t}} = \overline{\mu \pi_{t}} - \frac{1}{2}(\overline{\Sigma \pi_{t}^{2}} - \overline{\sigma \pi_{t}}^{2}), \quad \overline{x_{0}} = \exp \e \SBRA{\log x_{0}}, \quad \overline{\mu \pi_{t}} = \e \SBRA{\mu \pi_{t} | \F_{t}^{B}}, \quad \overline{c_{t}} = \e \SBRA{c_{t} | \F_{t}^{B}}.
\end{eqnarray*}
Since admissible controls $(\pi_{t}, c_{t})$ are $\F_{t}^{MF}$- progressively measurable, each $(\pi_{s}, c_{s})$ for $s \leq t$ is measurable with respect to $\sigma(\zeta, (B_{u})_{0 \leq u \leq s}, (W_{u})_{0 \leq u \leq s})$, hence independent of common noise increments $(B_{u} - B_{s})_{s \leq u \leq t}$, which justifies the equality between \eqref{espbruitcomm} and \eqref{espint}. For any admissible strategy $(\pi, c)$, set $\widehat{X}^{\pi, c} = \displaystyle\frac{X^{\pi, c}}{(\overline{X}^{\pi, c})^{\theta}}$. From \eqref{Xbardyn}, the dynamics of $\overline{X}^{\pi, c}$ are given by
\begin{eqnarray*}
\frac{d\overline{X_{t}}^{\pi, c}}{\overline{X_{t}}^{\pi, c}} = \overline{x_{0}} + (\eta -\overline{c_{s}})dt + \overline{\sigma \pi_{s}}dB_{t}.
\end{eqnarray*}
Applying It\^o's formula
\begin{multline} \label{eq:relwealthMF}
    \frac{d \widehat{X}^{\pi, c}}{\widehat{X_{t}}^{\pi, c}} = \PAR{\mu \pi_{t} - \theta \overline{\mu \pi_{t}} + \frac{\theta}{2} \overline{\Sigma \pi_{t}^{2}} + \frac{\theta^{2}}{2} \overline{\sigma \pi_{t}}^{2} - \theta \sigma \pi_{t} \overline{\sigma \pi_{t}}}dt \\
    + \nu \pi_{t} dW_{t} + (\sigma \pi_{t} - \theta \overline{\sigma \pi_{t}})dB_{t} - (c_{t} - \theta \overline{c_{t}})dt,
\end{multline}
with initial condition $\widehat{X_{0}}^{\pi, c} = \displaystyle\frac{x_{0}}{(\overline{x_{0}})^{\theta}}$. 

\vspace{0.3cm}
\noindent
\textbf{Forward relative performance criterion -} Applying It\^o-Ventzel's formula, we obtain the dynamics of the relative performance criterion along the relative wealth process $Q(t, \widehat{X_{t}}^{\pi, c})$:
\begin{align*}
dQ(t, \widehat{X_{t}}^{\pi, c}) & = (\beta^{U}(t, \widehat{X_{t}}^{\pi, c}) + V(\widehat{c_{t}} \widehat{X_{t}}^{\pi, c}))dt + \gamma^{U, W}(t, \widehat{X_{t}}^{\pi, c})dW_{t}  + \gamma^{U, B}(t, \widehat{X_{t}}^{\pi, c})dB_{t} + U_{x}(\widehat{X_{t}}^{\pi, c})d\widehat{X_{t}}^{\pi, c} \\ 
&\quad + \frac{1}{2} U_{xx}(t, \widehat{X_{t}}^{\pi, c}) d\CRO{\widehat{X_{t}}^{\pi, c}} + \CRO{\gamma_{x}^{U, W}(t, \widehat{X_{t}}^{\pi, c})dW_{t} + \gamma_{x}^{U, B}(t, \widehat{X_{t}}^{\pi, c})dB_{t}}, \widehat{X_{t}}^{\pi, c} \\
&= \PAR{\beta^{U}(t, \widehat{X_{t}}^{\pi, c}) + V(t, \widehat{c_{t}} \widehat{X_{t}}^{\pi, c})}dt + \gamma^{U, W}(t, \widehat{X_{t}}^{\pi, c})dW_{t} + \gamma^{U, B}(t, \widehat{X_{t}}^{\pi, c})dB_{t}   \\
&\quad \left. + U_{x}(t, \widehat{X_{t}}^{\pi, c}) \widehat{X_{t}}^{\pi, c} \PAR{\mu \pi_{t} - \theta \overline{\mu \pi_{t}} + \frac{\theta}{2} \overline{\Sigma \pi_{t}^{2}} + \frac{\theta^{2}}{2} \overline{\sigma \pi_{t}}^{2} - \theta \sigma \pi_{t} \overline{\sigma \pi_{t}} -(c_{t} - \theta \overline{c_{t}})} \right)dt \\
&\quad  + \frac{1}{2}U_{xx}(t, \widehat{X_{t}}^{\pi, c})(\widehat{X_{t}}^{\pi, c})^{2} \PAR{(\nu \pi_{t})^{2} + (\sigma \pi_{t} - \theta \overline{\sigma \pi_{t}})^{2}}dt  \\
&\quad + U_{x}(t, \widehat{X_{t}}^{\pi, c}) \widehat{X_{t}}^{\pi, c} \PAR{\nu \pi_{t}dW_{t} + (\sigma \pi_{t} - \theta \overline{\sigma \pi_{t}})dB_{t}} \\
&\quad + \PAR{\gamma_{x}^{U, W}(t, \widehat{X_{t}}^{\pi, c}) \nu \pi_{t} \widehat{X_{t}}^{\pi, c} + \gamma_{x}^{U, B}(t, \widehat{X_{t}}^{\pi, c}) \PAR{\sigma \pi_{t} - \theta \overline{\sigma \pi_{t}}}}dt.
\end{align*}

\vspace{0.3cm}
\noindent
\textbf{Best response strategy -} The first order conditions for martingale optimality give
\begin{align*}
\pi_{t}^{*}(x) &= \frac{1}{\nu^{2} + \sigma^{2}} \left( \theta \sigma \overline{\sigma \pi_{t}} - \frac{1}{U_{xx}(t, x) x} \PAR{\gamma_{x}^{U, W}(t, x) \nu + \gamma_{x}^{U, B}(t, x) \sigma + (\mu - \theta \sigma \overline{\sigma \pi_{t}})U_{x}(t, x)} \right) \\
c_{t}^{*}(x) &= \frac{V_{x}^{-1} \PAR{t, U_{x}(t, x) (\overline{C_{t}})^{\theta}} (\overline{C_{t}})^{\theta}}{x}. 
\end{align*}
With $\psi$ defined in \eqref{eq:defpsiMF}, the drift of $Q(t, \widehat{X_{t}}^{\pi, c})$ takes the form 
\begin{align*}
\drift Q(t, \widehat{X_{t}}^{\pi, c}) & = \beta^{U}(t, \widehat{X_{t}}^{\pi, c}) + \psi(t, \widehat{X_{t}}^{\pi, c}) + \frac{1}{2}U_{xx}(t, \widehat{X_{t}}^{\pi, c}) (\widehat{X_{t}}^{\pi, c})^{2}(\nu^{2} + \sigma^{2})(\pi_{t}^{2} - 2 \pi_{t} \pi_{t}^{*}(\widehat{X_{t}}^{\pi, c}))  \\
& \quad + \overset{\sim}{V}(t, (\overline{C_{t}})^{\theta} U_{x}(t, \widehat{X_{t}}^{\pi, c})) + V(t, \widehat{c_{t}} \widehat{X_{t}}^{\pi, c}) - c_{t} \widehat{X_{t}}^{\pi, c} U_{x}(t, \widehat{X_{t}}^{\pi, c}) \\
& \quad - \overset{\sim}{V}(t, (\overline{C_{t}})^{\theta} U_{x}(t, \widehat{X_{t}}^{\pi, c})).
\end{align*}
The drift condition ensuring that $Q$ is a mean-field forward relative performance criterion becomes
\begin{eqnarray} \label{eq:MFdriftcondgen}
    \beta(t, x) = - \psi(t, x) + \frac{1}{2} U_{xx}(t, x) x^{2} (\nu^{2} + \sigma^{2}) (\pi_{t}^{*}(x))^{2} - \overset{\sim}{V}(t, (\overline{C_{t}})^{\theta} U_{x}(t, x)),
\end{eqnarray}
which is equivalent to the mean-field HJB-SPDE \eqref{eq:MFSPDE}. Define $m_{t} = \Lr \PAR{\widehat{X_{t}^{*}} | \F_{t}^{B}}$ and the conditional mean $\overline{m_{t}} := \int y m_{t}(dy)$. Note that evaluating the relative wealth process \eqref{eq:relwealthMF} along the optimal strategy $(\pi_{t}^{*}, c_{t}^{*})$ makes the bar quantities $\overline{\eta_{t}}, \, \overline{\mu \pi_{t}}, \, \overline{c_{t}}, \, \overline{C_{t}}$ all depend on the conditional law $\overline{m_{t}}$. We make this dependence explicit in the following.

\noindent
\textbf{Mean equilibrium characterization -} Multiplying \eqref{portMF} by $\sigma$ and taking conditional expectation with respect to the common noise filtration, a fixed-point argument yields the conditional mean portfolio equilibrium
\begin{eqnarray}
\overline{\sigma \pi_{t}}(\overline{m_{t}}) = \frac{1}{\e \SBRA{\frac{\theta \sigma^{2}}{\nu^{2} + \sigma^{2}} + \frac{\theta \sigma^{2} U_{x}(t, \widehat{X_{t}^{*}})}{U_{xx}(t, \widehat{X_{t}^{*}}) \widehat{X_{t}^{*}}} | \F_{t}^{B}} - 1} \PAR{\e \SBRA{\frac{\sigma \nu \gamma_{x}^{U, W}(t, \widehat{X_{t}^{*}}) + \sigma^2 \gamma_{x}^{U, B}(t, \widehat{X_{t}^{*}})+ \sigma \mu U_{x}(t, \widehat{X_{t}^{*}})}{U_{xx}(t, \widehat{X_{t}^{*}}) \widehat{X_{t}^{*}}} | \F_{t}^{B}}}.
\end{eqnarray}
Hence, equilibrium exists provided that $\e \SBRA{\frac{\theta \sigma^{2}}{\nu^{2} + \sigma^{2}}} + \e \SBRA{\frac{\theta \sigma^{2} U_{x}(t, x)}{U_{xx}(t,x) x} | \F_{t}^{B}} - 1 \neq 0$, which is ensured by Assumption \ref{ass:MFadm}. We next make explicit the compatible aggregate consumption process $\overline{C_{t}}$ from \eqref{compatcond}. Taking the logarithm of the mean-field optimal consumption candidate \eqref{FOCconsMF} and conditioning on $\F_{t}^{B}$, we obtain
\begin{eqnarray} \label{eq:genimpfixpointbarC}
\log\PAR{\overline{C_{t}}(\overline{m_{t}})(1 - \theta)} = \e \SBRA{\log \PAR{V_{x}^{-1}(t, U_{x}(t, \widehat{X_{t}^{*}})(\overline{C_{t}}(\overline{m_{t}}))^{\theta})} - \log(\widehat{X_{t}^{*}}) | \F_{t}^{B}}.
\end{eqnarray}
Under Assumption \ref{ass:structconsMF}, there exists an $\F^{MF}$-adapted process and an $\F_{0}^{MF}$-measurable random variable $p$ such that
\begin{eqnarray*}
\e \SBRA{\log(c_{t}^{*}) | \F_{t}^{B}} = \e \SBRA{\log(K_{t} (\overline{C_{t}}(\overline{m_{t}}))^{p + \theta} | \F_{t}^{B}}.
\end{eqnarray*}
Together with \eqref{compatcond}, the above equation admits a unique solution $\overline{C_{t}}$,  independent of the conditional law, given by 
\begin{eqnarray} \label{eq:explbarCMF}
    \overline{C_{t}} = e^{\frac{1}{1 - \e \SBRA{p + \theta}}} \e \SBRA{\log(K_{t}) | \F_{t}^{B}}.
\end{eqnarray}

\noindent
\textbf{Admissibility -} The mean-field candidate equilibrium strategy \eqref{portMF} and \eqref{FOCconsMF} depend both on the state variable $x$ through utility related terms, and to its conditional law with respect to the common noise through $\overline{\sigma \pi_{t}}(\overline{m_{t}})$. We make the dependence of the mean-field equilibrium strategy in $\overline{m_{t}}$ explicit in the remaining of this proof. Admissibility of the mean-field equilibrium consumption reduces to proving that $\pi_{t}^{*}(\widehat{X_{t}^{*}},\overline{m_{t}})$ and $c_{t}^{*}(\widehat{X_{t}^{*}}, \overline{m_{t}})$ are such that the McKean-Vlasov dynamics with common noise of the relative optimal wealth process $\widehat{X_{t}^{*}}$ admits a unique global strong solution. From \eqref{eq:relwealthMF}, $\widehat{X_{t}^{*}}$ satisfies the McKean-Vlasov SDE 
\begin{eqnarray} \label{eq:MKVoptrelwealth}
d \widehat{X_{t}^{*}} = \beta^{X}(t, \widehat{X_{t}^{*}}, m_{t})dt + \gamma^{X, W}(t, \widehat{X_{t}^{*}}, m_{t})dW_{t} + \gamma^{X, B}(t, \widehat{X_{t}^{*}}, m_{t})dB_{t},
\end{eqnarray}
with coefficients written in terms of the conditional mean, 
\begin{multline}
\beta^{X}(t, x, m) = x\left( (\mu - \theta \sigma \overline{\sigma \pi_{t}}(\overline{m_{t}}) \pi_{t}^{*}(x, \overline{m_{t}}) - \theta \overline{\mu \pi_{t}}(\overline{m_{t}})  \right. \\
     \left. + \frac{\theta}{2} \overline{\Sigma \pi_{t}^{2}}(\overline{m_{t}}) + \frac{\theta^{2}}{2} (\overline{\sigma \pi_{t}}(\overline{m_{t}}))^{2} - \PAR{c_{t}^{*}(x, \overline{m_{t}}) - \theta \overline{c_{t}}(\overline{m_{t}})} \right) 
\end{multline}
\vspace{-0.8cm}
\begin{align}
\gamma^{X, W}(t, x, m) = x \nu \pi_{t}(x, \overline{m_{t}}) , \quad 
\gamma^{X, B}(t, x, m) = x \PAR{ \sigma \pi_{t}(x, \overline{m_{t}}) - \theta \overline{\sigma \pi_{t}}(\overline{m_{t}})}.
\end{align}
Under assumption \ref{ass:structconsMF}, the mean equilibrium consumption process $\overline{C_{t}}$ is uniquely determined by \eqref{eq:explbarCMF}, and is thus independent of the law $m$. Consequently, under Assumption \ref{ass:MFadm} and \ref{ass:structconsMF}, $\beta^{X}$ and $\gamma^{X}$ are independent from the conditional law $m$ and satisfy the coercivity and sublinear Lipschitz growth properties in the state variable. In turn, \eqref{eq:MKVoptrelwealth} reduces to a standard SDE, and Theorem 1.4 in \cite{lan2014new} yields existence of a strong global solution. This proves admissibility of the equilibrium strategy.

\begin{rmq}
\begin{itemize}
    \item Controlling coercivity and Lipschitz property of $\gamma^{X, B}$ requires control of the cross term in $x \pi_{t}^{*}(x, \overline{m}) \overline{\sigma \pi_{t}}(\overline{m})$. This motivates Assumption \ref{ass:MFadm}, which ensures that $\overline{\sigma \pi_{t}}$ is independent of the conditional law $\overline{m}$. In particular this is satisfied in the separable CRRA case considered in Section \ref{sect:MF_CRRA}.
    \item Existence and uniqueness results for McKean-Vlasov SDE with common noise, e.g \cite{kumar2022well}, typically require joint Lipschitz and linear growth bounds in $(x, m)$. Assumption \ref{ass:structconsMF} is essential to ensure that $\overline{C_{t}}$ can be solved explicitly, while Assumptions \ref{ass:MFadm} and \ref{ass:structconsMF} reduce the problem to a standard SDE with random coefficients. This allows us to lighten the strict monotonicity in space requirements from \cite{kumar2022well}, to stochastic Lipschitz properties.
    \item Stating this verification result under weaker hypothesis seems difficult because of the non linearity of the mean-field equilibrium consumption \eqref{FOCconsMF}. One must first establish that equation \eqref{eq:genimpfixpointbarC} admits a unique fixed-point solution $\overline{C_{t}}$, despite the heavy coupling induced by $\widehat{X_{t}^{*}}$ which itself depends on $\overline{C_{t}}$. Additionally, applying Theorem 2.1 from \cite{kumar2022well} would also require joint Lipschitz property with respect to space and law of the map $(x, m) \mapsto (V_{x})^{-1} \PAR{t, U_{x}(t, x)(\overline{C_{t}}(\overline{m})^{\theta}} (\overline{C_{t}}(\overline{m})^{\theta}$, which is unclear with an implicit definition of the mean equilibrium consumption.    
\end{itemize}
    
\end{rmq}
\end{preuve}

\subsection{Mean-field Nash equilibrium for separable utilities of CRRA type} \label{sect:MF_CRRA}

In this section, we study the existence of a mean-field Nash equilibrium for the optimization problem associated with a forward relative performance criterion $Q$ of the form \eqref{eq:MF_frpp}, generated by utility random fields of CRRA type. As in the $n$-player case considered in Section \ref{sect:consutCRRA}, analogous calculations show that if a utility pair $(U^{\alpha}, V)$ satisfies the mean-field SPDE \eqref{eq:MFSPDE}, and $U^{\alpha}$ has separable power form, then the consumption utility $V$ must also be of separable CRRA type, with the same risk aversion parameter $\alpha$. Consequently, \eqref{eq:MFSPDE} reduces to a drift condition on the time component of $U^{\alpha}$, expressed in terms of the time component of $V$.

\begin{ass}[Mean-field CRRA] \label{ass_MF_CRRA}
Consider a random field $Q$ of the form
\begin{eqnarray*}
Q(t, \widehat{x}) = U^{(\alpha)}(t, \widehat{x}) + \int_{0}^{t} V^{(\alpha)}(s, \widehat{c_{s}} \widehat{x})ds, 
\end{eqnarray*}
where $U^{(\alpha)}$ and $V^{(\alpha)}$ are CRRA utility random fields
\begin{eqnarray} \label{eq:sep_utMF}
U^{(\alpha)}(t, x) = Z_{t} u^{(\alpha)}(x), \quad \text{and} \quad V^{(\alpha)}(t, c) = \phi_{t} u^{(\alpha)}(c),
\end{eqnarray}
with $u^{(\alpha)}(x) = \frac{x^{1 - \alpha}}{1 - \alpha}$, $\alpha \in \A^{\alpha} := (0, 1) \cup (1, +\infty)$. Assume processes $(Z_{t})_{t \geq 0}$, $(\phi_{t})_{t \geq 0}$ follow log-normal dynamics
\begin{eqnarray} \label{dynZMF}
    dZ_{t} = Z_{t} \PAR{b^{Z}(t, Z_{t}, \phi_{t}) dt + \delta^{Z, W}(t)dW_{t} + \delta^{Z, B}(t) dB_{t}}, \quad Z_{0} = z_{0}. \\
    d\phi_{t} = \phi_{t} \PAR{b^{\phi}(t, Z_{t}, \phi_{t}) dt + \delta^{\phi, W}(t)dW_{t} + \delta^{\phi, B}(t) dB_{t}}, \quad \phi_{0} = \phi_{0}. \nonumber
\end{eqnarray}
\end{ass}
\noindent
Under Assumption \ref{ass_MF_CRRA}, the mean-field HJB-SPDE \eqref{eq:MFSPDE} is equivalent to 
\begin{eqnarray} \label{driftcondMF_prop}
b^{Z}(t, Z_{t}, \phi_{t}) = \overline{b^{Z}}(t) - (1 - \alpha) \theta \overline{c_{t}} - \frac{\alpha}{\PAR{\overline{C_{t}}}^{\frac{\theta (1- \alpha)}{\alpha}}} \PAR{\frac{\phi_{t}}{Z_{t}}}^{\frac{1}{\alpha}},
\end{eqnarray}
where the $Z$ independent component is
\begin{multline}
\overline{b^{Z}}(t) = (1 - \alpha) \delta^{Z, B}(t) \theta \overline{\sigma \pi_{t}} - (1 - \alpha) \PAR{- \theta \overline{\mu \pi_{t}} + \frac{\theta}{2} \overline{\Sigma \pi_{t}^{2}} + \frac{\theta^{2}}{2} \overline{\sigma \pi_{t}}^{2}} \\
+ \frac{1}{2} \alpha (1 - \alpha) \PAR{(\theta \overline{\sigma \pi_{t}})^{2} - (\nu^{2} + \sigma^{2})(\pi_{t}^{*})^{2}}.
\end{multline}
For convenience in the calculation, we include the dynamic preference factor, composed of $Z$ and $\phi$ local characteristics (except $\overline{b^{Z}}$), in the random type vector, which in turn takes the form
\begin{eqnarray*}
    \zeta = \PAR{\delta^{Z, W}, \delta^{Z, B}, b^{\phi}, \delta^{\phi, W}, \delta^{\phi, B}, z_{0}, \phi_{0}, x_{0}, \alpha, \theta, \mu, \nu, \sigma},
\end{eqnarray*}
with values on the type space
\begin{eqnarray*}
\Z = \PAR{L^{\infty}_{\loc}(\R^{+, *})}^{2} \times \Hr^{1}_{\loc}(\R) \times \PAR{L^{\infty}_{\loc}(\R^{+, *})}^{2} \times (0, \infty) \times (0, \infty) \times (0, \infty) \times \A^{\alpha} \times [0,1] \times \R \times \R^{+} \times \R^{+}.
\end{eqnarray*}
\noindent
Next, we derive a candidate mean-field Nash equilibrium for the class of separable CRRA utilities satisfying Assumption \ref{ass_MF_CRRA}. For the following, define
\begin{eqnarray} \label{defquant_MF}
\psi^{\sigma} = \e \SBRA{\theta (1 - \frac{1}{\alpha}) \frac{\sigma^{2}}{\nu^{2} + \sigma^{2}}}, \quad \varphi^{\sigma}(t) = \e \SBRA{ \frac{\sigma}{\alpha} \frac{ (\delta^{Z, W}(t) \nu + \delta^{Z, B}(t) \sigma + \mu)}{\nu^{2} + \sigma^{2}}}, \quad
K_{\alpha, \theta} = - \frac{\frac{\theta(1 - \alpha)}{\alpha}}{1 + \e \SBRA{\frac{\theta (1 - \alpha)}{\alpha}}}.
\end{eqnarray}
 
\begin{thm} \label{MFthm}
Let $(U^{(\alpha)}, V^{(\alpha)})$ be a pair of utility random fields satisfying Assumptions \ref{ass_MF_CRRA} and the mean-field HJB-SPDE \eqref{eq:MFSPDE}. If $\psi^{\sigma} \neq 1$, then the candidate mean-field Nash equilibrium strategy $(\pi_{t}^{*}, c_{t}^{*})$ takes the form
\begin{align} 
 \pi_{t}^{*} &= \frac{1}{\nu^{2} + \sigma^{2}} \PAR{\theta \sigma (1 - \frac{1}{\alpha}) \frac{\varphi^{\sigma}(t)}{1 - \psi^{\sigma}} + \frac{1}{\alpha}(\delta^{Z, W}(t) \nu + \delta^{Z, B}(t) \sigma + \mu)}, \label{optportMF} \\
c_{t}^{*} &=  \exp \PAR{\e \SBRA{ \frac{1}{\alpha} \log \PAR{ \frac{\phi_{t}}{Z_{t}}} | \F_{t}^{B}}}^{K_{\alpha, \theta}} \PAR{\frac{\phi_{t}}{Z_{t}}}^{\frac{1}{\alpha}}. \label{MFconsdyn}
\end{align}
Assume in addition that the drift of $\phi$ has a mean-reverting component, with the same dependence on mean consumption as the drift of $Z$, namely
\begin{eqnarray} \label{eq:condphiMF}
    b^{\phi}(t, Z_{t}, \phi_{t}) = \overline{b^{\phi}}(t) - \lambda \SBRA{(1 - \alpha)\theta \overline{c_{t}} + \alpha c_{t}^{*}}
\end{eqnarray}
where $\lambda \geq 1$ is a deterministic preference parameter, common among the continuum of agents. Then $(\pi_{t}^{*}, c_{t}^{*}) \in \A^{MF}$ is a well defined mean-field Nash equilibrium.
\end{thm}

\noindent
Note that the mean-field equilibrium $(\pi_{t}^{*}, c^{*}_{t})$ is stochastic and $\F_{t}^{MF}$-progressively measurable. In particular, Theorem \ref{MFthm} provides a sufficient condition \eqref{eq:condphiMF} ensuring that the preference criterion is dynamically consistent and the optimal strategies admissible. Therefore, the admissible random type vector in the separable CRRA framework is obtained by replacing $b^{\phi}$ with its unrestricted part $\overline{b^{\phi}}$, that is 
\begin{eqnarray*}
    \zeta = \PAR{\delta^{Z, W}, \delta^{Z, B}, \overline{b^{\phi}}, \delta^{\phi, W}, \delta^{\phi, B}, z_{0}, \phi_{0}, x_{0}, \alpha, \theta, \mu, \nu, \sigma}.
\end{eqnarray*}
According to Theorem \ref{MFthm}, each realization $\zeta \in \Z$ defines an admissible pair of CRRA random utilities $(U^{(\alpha)}, V^{(\alpha)})$ which generates a mean-field forward relative performance criterion.

\begin{preuve}
\vspace{0.3cm}
\textbf{Candidate mean-field Nash equilibrium strategies -} For separable CRRA forward criteria, the optimal portfolio \eqref{portMF} is
\begin{eqnarray*}
    \pi_{t}^{*} = \frac{1}{\nu^{2} + \sigma^{2}} \PAR{\theta \sigma \overline{\sigma \pi_{t}}(1 - \frac{1}{\alpha}) + \frac{1}{\alpha}(\delta^{Z, W}(t) \nu + \delta^{Z, B}(t) \sigma + \mu)}.
\end{eqnarray*}
Multiplying by $\sigma$ and taking the conditional expectation with respect to $\F_{t}^{B}$ gives
\begin{eqnarray*}
    \overline{\sigma \pi_{t}} = \psi^{\theta} \overline{\sigma \pi_{t}} + \varphi^{\sigma}(t),
\end{eqnarray*}
Since the coefficients above are $\F_{0}^{MF}$-measurable, they are independent of $\F_{t}^{B}.$ Hence, if $\psi^{\sigma} \neq 1$, there exists a mean-field portfolio equilibrium $\pi_{t}^{*}$ given by \eqref{portMF} with $\overline{\sigma \pi_{t}^{*}} = \frac{\varphi^{\sigma}(t)}{1 - \psi^{\sigma}}$. For consumption, the candidate mean-field Nash equilibrium strategy in the separable CRRA case is 
\begin{eqnarray*}
c_{t}^{*} = \frac{1}{\PAR{\overline{C_{t}}}^{\frac{\theta(1 - \alpha)}{\alpha}}} \PAR{\frac{\phi_{t}}{Z_{t}}}^{\frac{1}{\alpha}}.
\end{eqnarray*}
Substituting the explicit expression of $\overline{C_{t}}$ from \eqref{eq:explbarCMF} yields \eqref{MFconsdyn}.

\noindent
\textbf{Admissibility -} For greater generality, assume that the drift of the time-dependent part of the consumption utility is of the form
\begin{eqnarray*}
        b^{\phi}(t, Z_{t}, \phi_{t}) = \overline{b^{\phi}}(t) - (1 - \alpha)\theta \overline{\lambda} \overline{c_{t}} - \lambda \alpha c_{t}^{*},
\end{eqnarray*}
where $\overline{\lambda}$ and $\lambda$ are real-valued random variables which may be added to the random type vector. After It\^o calculations similar to those in Section \ref{sect:consutCRRA}, setting $Y_{t} = \PAR{\frac{\phi_{t}}{Z_{t}}}^{\frac{1}{\alpha}}$, we obtain
\begin{eqnarray*}
d \log (Y_{t}) = \PAR{a_{t} + \frac{(1 - \alpha)\theta}{\alpha} (1 - \overline{\lambda}) \overline{c_{t}} - (\lambda - 1) c_{t}^{*}}dt + \delta^{Y, W}(t)dW_{t} + \delta^{Y, B}(t) dB_{t},
\end{eqnarray*}
with 
\begin{align*}
    a_{t} &= \frac{1}{\alpha} \PAR{\overline{b^{\phi}}(t) - \overline{b^{Z}}(t)}(t) - \frac{1}{2}\PAR{(\delta^{\phi, W}(t))^{2} + (\delta^{\phi, B}(t))^{2} - (\delta^{Z, W}(t))^{2} - (\delta^{Z, B}(t))^{2}} \\
    \delta^{Y, W}(t) &= \frac{1}{\alpha} \PAR{\delta^{\phi, W}(t) - \delta^{Z, W}(t)}, \quad \delta^{Y, B}(t) = \frac{1}{\alpha} \PAR{\delta^{\phi, B}(t) - \delta^{Z, B}(t)}.
\end{align*}
Taking conditional expectation with respect to $\F_{t}^{B}$, and denoting $\overline{a_{t}} = \e \SBRA{a_{t} | \F_{t}^{B}}$ and $\overline{\delta^{Y, B}}(t) = \e \SBRA{\delta^{Y, B}(t) | \F_{t}^{B}}$ gives 
\begin{eqnarray*}
d \e \SBRA{\log (Y_{t}) | \F_{t}^{B}} = \PAR{\overline{a_{t}} + \e \SBRA{\frac{(1 - \alpha)\theta}{\alpha} (1 - \overline{\lambda}) \overline{c_{t}} \F_{t}^{B}} - \e \SBRA{(\lambda - 1) c_{t}^{*} | \F_{t}^{B}}}dt + \overline{\delta^{Y, B}}(t) dB_{t}.
\end{eqnarray*}
If $\overline{\lambda}$ and $\lambda$ are common across the population, thus deterministic, they can be taken outside expectations so 
\begin{eqnarray*}
d \e \SBRA{ \log (Y_{t}) | \F_{t}^{B}} = \PAR{\overline{a_{t}} + \overline{c_{t}} \PAR{\e \SBRA{\frac{(1 - \alpha)\theta}{\alpha}} (1 - \overline{\lambda}) - (\lambda - 1)}}dt + \overline{\delta^{Y, B}}(t) dB_{t}.
\end{eqnarray*}
Since $c_{t}^{*} = \exp \PAR{ \log(Y_{t}) + K_{\alpha, \theta} \e \SBRA{ \log (Y_{t}) | \F_{t}^{B}}} $, the candidate equilibrium consumption satisfies
\begin{align} \label{eq:MFconsdynpreuve}
    d \log (c_{t}^{*}) = \PAR{ a_{t} + K_{\alpha, \theta} \overline{a_{t}} + \eta \overline{c_{t}} - (\lambda - 1) c_{t}^{*}}dt + \delta^{Y, W}(t)dW_{t} + \PAR{\delta^{Y, B}(t) + K_{\alpha, \theta} \overline{\delta^{Y, B}}(t)} dB_{t},
\end{align}
where the coefficient in front of conditional mean consumption simplifies as
\begin{align}
    \eta = K_{\alpha, \theta} (\overline{\lambda} - \lambda).
\end{align}
Hence, $c_{t}^{*}$ follows a common noise McKean-Vlasov SDE, with non positive quadratic drift when $\lambda \geq 1$. Additionally, setting $\lambda = \overline{\lambda}$ kills the dependence of these dynamics on the conditional law, and the equation reduces to a standard SDE. Existence and uniqueness of a strong global solution $c_{t}^{*}$ follows from the same inverse-transform argument as in the proof of Theorem \ref{thm:optstratpowerN}, which proves admissibility of the mean-field Nash consumption strategy. In this case, $c_{t}^{*}$ is stochastic and solves 
\begin{multline}
    d c_{t}^{*} = c_{t}^{*}\PAR{ a_{t} + K_{\alpha, \theta} \overline{a_{t}} - (\lambda - 1) c_{t}^{*} + \frac{1}{2}(\delta^{Y, W})^{2}(t) + \frac{1}{2}\PAR{\delta^{Y, B}(t) + K_{\alpha, \theta} \overline{\delta^{Y, B}}(t)}^{2}}dt \\
    + c_{t}^{*}\delta^{Y, W}(t)dW_{t} + c_{t}^{*}\PAR{\delta^{Y, B}(t) + K_{\alpha, \theta} \overline{\delta^{Y, B}}(t)} dB_{t},
\end{multline}
Define $\Gamma_{t} = \exp \PAR{- \displaystyle\int_{0}^{t} (a_{s} + K_{\alpha, \theta} \overline{a_{s}})ds - \int_{0}^{t} \delta^{Y, W}(s)dW_{s} - \int_{0}^{t} \PAR{\delta^{Y, B}(s) + K_{\alpha, \theta} \overline{\delta^{Y, B}}(s) } dB_{s}}$. Then, the mean-field Nash equilibrium consumption admits the explicit representation 
\begin{eqnarray}
    c_{t}^{*} = \frac{c_{0}^{*} \Gamma_{t}^{-1}}{ 1 + (\lambda - 1)c_{0}^{*} \int_{0}^{t} \Gamma_{s}^{-1}ds}, \quad \text{with} \quad c_{0}^{*} = \exp \PAR{\e \SBRA{\frac{1}{\alpha} \log \PAR{\frac{\phi_{0}}{z_{0}}}}}^{K_{\alpha, \theta}} \PAR{\frac{\phi_{0}}{z_{0}}}^{\frac{1}{\alpha}}.
\end{eqnarray}
\end{preuve}

\begin{rmq}
\begin{itemize}
\item Allowing $\lambda$ to vary among the population introduces heterogeneity in the law-dependent term of the $c_{t}^{*}$ dynamics. Existence of a well defined mean-field Nash equilibrium in this case remains open. 
\item When $\eta \neq 0$, existence of a strong global solution to \eqref{eq:MFconsdynpreuve} is unclear. Indeed, setting $I_{t}^{*} = \frac{1}{c_{t}^{*}}$ linearizes the equation, but the law dependence becomes $\e \SBRA{ \frac{1}{I_{t}^{*}} | \F_{t}^{B}}$, for which the global Lipschitz conditions used in \cite{kumar2022well} cannot be ensured. 
\end{itemize}
\end{rmq}

\subsection{Examples of mean-field Nash equilibrium} \label{sect:MFexamples}

Next, we analyze the equilibrium consumption strategy \eqref{MFconsdyn} under two simplifying assumptions. First, imposing a proportional relationship between wealth and consumption utility yields a tractable setting in which the optimal strategy is $\F_{0}^{MF}-$ measurable, providing an example of a strong mean-field Nash equilibrium. Second, motivated by \cite{dos2022forward}, we consider a non-linear relation between wealth and consumption utilities, by linking their time-varying components through a power-type relation. In that case, Theorem \ref{MFthm} ensures well-posedness, and the equilibrium consumption admits a closed-form logistic representation.
 
\paragraph{Proportional consumption-wealth preference -} Since the ratio $\frac{\phi_{t}}{Z_{t}}$ appears in the mean-field equilibrium consumption, and therefore also in the mean-field HJB-SPDE \eqref{eq:MFSPDE}, a natural simplification is to assume that this relative preference for consumption versus wealth is specified ex ante.

\begin{ass}\label{ass_strong_MF}
$\phi_{t} = k(t) Z_{t}$, for some continuous function $k : \R^{+} \to \R^{+, *}$.
\end{ass}
In the mean-field setting, $k$ must be considered as a random variable, valued in the space of continuous strictly positive functions, representing the distribution of this proportional factor over the population. We add this parameter to the type vector $\zeta,$ requiring that $\e \SBRA{\log(k(t))} < + \infty$ for all $t \geq 0$.

\begin{prop}
Under Assumptions \ref{ass_MF_CRRA} and \ref{ass_strong_MF}, if $\psi^{\sigma} \neq 1$, there exists a strong mean-field Nash equilibrium strategy $(\pi_{t}^{*}, c_{t}^{*})$, $\F_{0}^{MF}$-measurable, given by
\begin{align} 
 \pi_{t}^{*} &= \frac{1}{\nu^{2} + \sigma^{2}} \PAR{\theta \sigma (1 - \frac{1}{\alpha}) \frac{\varphi^{\sigma}(t)}{1 - \psi^{\sigma}} + \frac{1}{\alpha}(\delta_{W}^{Z}(t) \nu + \delta_{B}^{Z}(t) \sigma + \mu)}, \label{optportMFstrong} \\
c_{t}^{*} &= \exp \PAR{ K_{\alpha, \theta} \e \SBRA{\frac{\log(k(t))}{\alpha}}}k(t)^{\frac{1}{\alpha}}.  \label{optconsMFstrong}
\end{align}
\end{prop}
Under Assumption \ref{ass_strong_MF}, forward relative performance criteria with non zero volatility produce a strong Nash equilibrium with explicit time dependence: $\pi^{*}$ varies over time through the volatility coefficients $\delta^{Z}$ and its mean via $\varphi{\sigma}$, while $c^*$ varies over time through the proportional consumption-wealth preference parameter $k$. By contrast, in the zero-volatility forward setting, \cite{dos2022forward} obtain a strong Nash equilibrium with constant portfolio, while consumption remains time-dependent, also driven by the ratio $\frac{\phi_{t}}{Z_{t}}$, in a less transparent way.

\paragraph{Relative market-consumption preference -} Second, we consider a non-linear power-type coupling between the time-varying components of the consumption and wealth utilities.

\begin{ass} \label{ass:kappa}
$\phi_{t} = Z_{t}^{1 - \kappa}$, where $\kappa \in \R$ 
 is the \textit{risk relative consumption preference parameter}.
\end{ass}

\noindent
Under Assumption \ref{ass:kappa}, Theorem \ref{MFthm} applies and yields a well-defined mean-field Nash equilibrium with explicit consumption, provided that $\kappa \leq 0$ and is common across the population. This restriction also appears in \cite{dos2022forward}, as it enables decoupling of the optimal consumption. Our calculations indicate that allowing $\kappa$ to be random introduces heterogeneity in the conditional law dependence of the equilibrium consumption, and admissibility becomes unclear. To state the equilibrium in closed form, define the stochastic exponential
\begin{eqnarray}
\Gamma^{A}_{t} = \exp \PAR{-\int_{0}^{t} b^{A}(s)ds - \int_{0}^{t} \delta^{A, W}(s) dW_{s} - \int_{0}^{t} \delta^{A, B}(s) dB_{s}},
\end{eqnarray}
with
\begin{align} \label{driftAkappa}
b^{A}(t) &= \frac{\kappa}{2 \alpha} \PAR{2 \alpha - 1 - \kappa(\alpha - 1)} \NRM{\delta^{Z}(t)}^{2} - \frac{\kappa}{\alpha} \overline{b^{Z}}(t) + K_{\alpha, \theta} \PAR{\frac{\kappa^{2}}{2} \e \SBRA{\frac{1}{\alpha}} \e \SBRA{\NRM{\delta^{Z}(t)}^{2}} - \kappa \e \SBRA{\frac{\overline{b^{Z}}(t)}{\alpha}}} \\
\delta^{A, B}(t) &= - \frac{\kappa}{\alpha} \delta^{Z, B}(t) - \kappa K_{\alpha, \theta} \e \SBRA{\frac{1}{\alpha}} \e \SBRA{\delta^{Z, B}(t)}, \quad  \delta^{A, W}(t) = - \frac{\kappa}{\alpha} \delta^{Z, W}(t). \label{volAWkappa}
\end{align}

\begin{prop} \label{prop:CRRA rel}
Suppose Assumptions \ref{ass_MF_CRRA} and \ref{ass:kappa} hold with $\kappa \leq 0$ common to all agents. If $\psi^{\sigma} \neq 1$, then there exists a unique mean-field Nash equilibrium strategy $(\pi_{t}^{*}, c_{t}^{*}) \in \A^{MF}$. The optimal portfolio is
\begin{eqnarray} \label{optportMFkappa}
     \pi_{t}^{*} = \frac{1}{\nu^{2} + \sigma^{2}} \PAR{\theta \sigma (1 - \frac{1}{\alpha}) \frac{\varphi^{\sigma}(t)}{1 - \psi^{\sigma}} + \frac{1}{\alpha}(\delta^{Z, W}(t) \nu + \delta^{Z, B}(t) \sigma + \mu)},
\end{eqnarray}
\vspace{0.3cm}
and the equilibrium consumption process $(c_{t}^{*})_{t \geq 0}$ satisfies
\begin{eqnarray}
c_{t}^{*} = \frac{c_{0}^{*} (\Gamma_{t}^{A})^{-1}}{1 + \ABS{\kappa} c_{0}^{*} \int_{0}^{t} (\Gamma^{A}_{s})^{-1} ds}, \quad \text{with} \quad c_{0}^{*} = \exp \PAR{\e \SBRA{\frac{1}{\alpha} \log \PAR{\frac{\phi_{0}}{z_{0}}}}}^{K_{\alpha, \theta}} \PAR{\frac{\phi_{0}}{z_{0}}}^{\frac{1}{\alpha}}.
\end{eqnarray}
\end{prop}

Observe that the resulting Nash equilibrium remains time-dependent, and that the optimal consumption is stochastic. Utility volatility appears as an additional preference input from the agent, scaling both the optimal portfolio and consumption. This is the main novelty of the non-zero volatility setting for forward relative performance, compared with the deterministic equilibrium consumption obtained in \cite{dos2022forward}, which is characterized by an ordinary differential equation.

\section{Discussion of the equilibrium} \label{sect:discMF}

In this section, we analyze the equilibrium arising from the mean-field optimization problem under forward relative performance criteria. In particular, we show that our framework naturally recovers the single-agent optimization problem of \cite{el2018consistent} when $\theta = 0$, and the locally riskless setting of \cite{dos2022forward} when utility volatility vanishes. We analyze optimal consumption in detail under Assumptions \ref{ass_strong_MF} and \ref{ass:kappa}, and highlight the implication of non-zero volatility.

\subsection{Investment strategy}

\begin{figure}[ht]
\centering
\begin{minipage}[t]{0.48\textwidth}
        \includegraphics[scale=0.55]{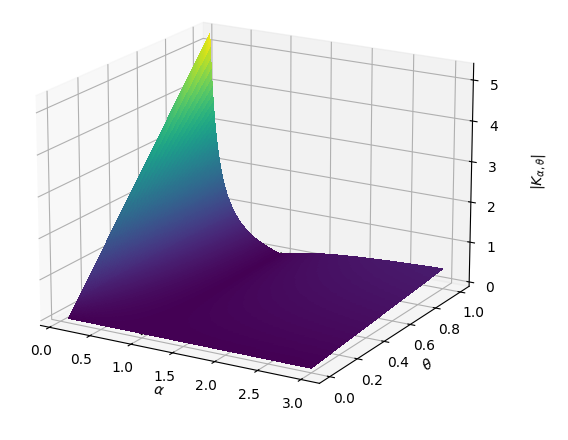}
        \caption{$|K_{\alpha, \theta}|$ function of risk aversion $\alpha$ and competition parameter $\theta$}
    \label{fig:surf_K}
\end{minipage} \hfill
\begin{minipage}[t]{0.48\textwidth}
        \includegraphics[scale=0.55]{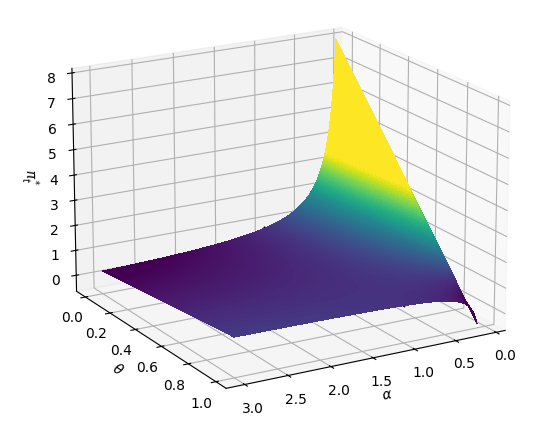}
        \caption{Optimal investment for $\mu = 0.3$, $\sigma = 1$, $\delta^{Z, B} = \e \SBRA{\delta^{Z, B}} = 0.5$, $\e \SBRA{\frac{1}{\alpha}} = 2$ and $\e \SBRA{\theta} = 0.7.$}
    \label{surf_opt_port}
\end{minipage}
\end{figure}

In the CRRA relative performance setting, the mean-field equilibrium strategy depends on time only through the volatility parameters of the wealth utility, namely $\delta_{W}^{Z}, \delta_{B}^{Z}$ and their expectations. More precisely, the optimal investment $\pi_{t}^{*}$ can be written as the sum of $\pi_{t}^{1, *}$ and $\pi_{t}^{2, *}$ where: 
\begin{eqnarray} \label{eq:opt_port_decomp}
\pi_{t}^{1, *} = \frac{1}{\nu^{2} + \sigma^{2}} \theta \sigma \PAR{1 - \frac{1}{\alpha}} \frac{\varphi^{\sigma}(t)}{1 - \psi^{\sigma}} \quad \text{and} \quad \pi_{t}^{2, *} = \frac{\delta^{Z, W}(t) \nu + \delta^{Z, B}(t) \sigma + \mu}{\alpha (\nu^{2} + \sigma^{2})}.
\end{eqnarray}
In contrast, the deterministic CRRA relative utility framework of \cite{lacker2019mean}, \cite{lacker2020many}, and the forward CRRA relative performance model with zero volatility in \cite{dos2022forward} yield time-independent optimal investment strategies. In our setting, the effect of competition is captured by $\pi_{t}^{1, *}$. As in \cite{lacker2019mean}, the optimal portfolio given by \eqref{optportMF} is decreasing in the competition parameter $\theta$: the more an agent values relative rather than absolute performance, the smaller her optimal portfolio allocation. As $\theta$ goes to zero, the optimal portfolio converges to the limit process $\pi_{t}^{2, *}$ which corresponds to the optimal strategy in the non-competitive framework, see \cite{el2018consistent}. In other words, in presence of competition, the optimal strategy is corrected by $\pi_{t}^{1, *}$, whose magnitude is given by the coefficient $K_{\alpha, \theta}$ defined in \eqref{defquant_MF}, which we represent in Figure \ref{fig:surf_K}. Observe that, for a fixed risk aversion, the magnitude of the correction $K_{\alpha, \theta}$ grows linearly in $\theta$. Additionally, the sign of $K_{\alpha, \theta}$ is the same as that of $\alpha - 1$. Indeed, the mean population competition and risk aversion parameters act only as scaling factors, whereas individual preference parameters determine whether the correction is positive or negative. Figure \ref{fig:surf_K} shows that extremely risk-averse behavior leads to a slight increase of the optimal portfolio allocation for competitive agents. By contrast, under risk-seeking behavior, the effect is more sensitive to competition preferences: for less competitive agents, the optimal allocation significantly increases, whereas for more competitive agents the optimal allocation decreases. In particular, for $\theta$ close to $1$ and $\alpha$ close to $0$, the allocation may become negative, meaning that competitive, risk-seeking agents tend to go short, in line with observations in \cite{lacker2019mean}.

\vspace{0.3cm}
As $\alpha$ decreases to $0$, the time-dependent quantity $\pi_{t}^{2, *}$ increases. Thus, the portfolio allocation tends to be more sensitive to individual utility volatility preferences for risk-seeking agent. Moreover, $\pi_{t}^{1, *}$ also decreases as $\alpha$ goes to $0$, so a risk-seeking agent will see her optimal investment strategy more strongly influenced by competition. In Figure \ref{surf_opt_port}, we plot the equilibrium portfolio strategy as a function of $\alpha$ and $\theta$ in the single stock model, that is assuming $\nu, \sigma, \mu$ deterministic with $\nu=0$. We refer to \cite{lacker2019mean} for extended discussions of this optimal portfolio equilibrium.

In the mean-field portfolio equilibrium, non-zero utility volatility acts as a time-varying scaling for both the benchmark equilibrium portfolio $\pi_{t}^{2, *}$, and the competitive correction term $\pi_{t}^{1, *}$ through expectations in $\varphi^{\sigma}$. Hence, compared with the zero volatility benchmark, non-zero volatility does not change the economic mechanism of competition, rather, it rescales the portfolio level according to the representative agent's utility sensitivity to environmental noise. Note also that the mean-field Nash equilibrium portfolio is independent of the relative importance the agents give to wealth utility compared to consumption utility. By contrast, we will see in Section \ref{sect:4cons} that the optimal consumption rate depends heavily on this preference, characterized by the ratio $\frac{\phi_{t}}{Z_{t}}$.

\subsection{Consumption strategy} \label{sect:4cons}

Our attention now shifts to the optimal consumption process $c_{t}^{*}$ given by \eqref{MFconsdyn}. This equilibrium strategy is random, driven by the process $\frac{\phi_{t}}{Z_{t}}$, which represents the weight the representative agent gives to her utility from consumption, relatively to her wealth utility. 

\paragraph{Proportional consumption-wealth preference equilibrium -} First, assume that every agent has a utility from consumption proportional to their utility from wealth, as in Assumption \ref{ass_strong_MF}. Let $\phi_{t} = KZ_{t},$ for some positive real random variable $K$ such that $\e \SBRA{\log(K)} < +\infty$, that we assume time-independent for simplicity, and add in the type vector $\zeta$. The optimal consumption process $c_{t}^{*}$ is then a function of $K$, the competition parameter $\theta$ and the risk aversion $\alpha$, given by
\begin{eqnarray} \label{optconsMFproportionnel}
    c_{t}^{*} = \exp \PAR{ \frac{\log(K)}{\alpha} + K_{\alpha, \theta} \e \SBRA{\frac{\log(K)}{\alpha}}}.
\end{eqnarray}

This equilibrium consumption $c_{t}^{*}$ is independent of utility volatility. Nevertheless, this case is worth highlighting because it makes the role of the ratio $\frac{\phi_{t}}{Z_{t}}$ explicit in the consumption equilibrium. In other words, the relative weight of wealth utility over consumption utility drives the equilibrium consumption rate over time. In this specification, $K$ is $\F_{0}^{MF}$-measurable, and therefore constitutes a primitive preference input for both the representative agent and the population. By contrast, the equilibria in \cite{lacker2020many} and \cite{dos2022forward} take a weighted exponential form, in which the effect of this relative preference is less transparent. This gives our approach a clear economic interpretation of equilibrium outcomes.

\begin{figure}[ht]
    \centering
    \makebox[\textwidth][c]{\includegraphics[scale=0.55]{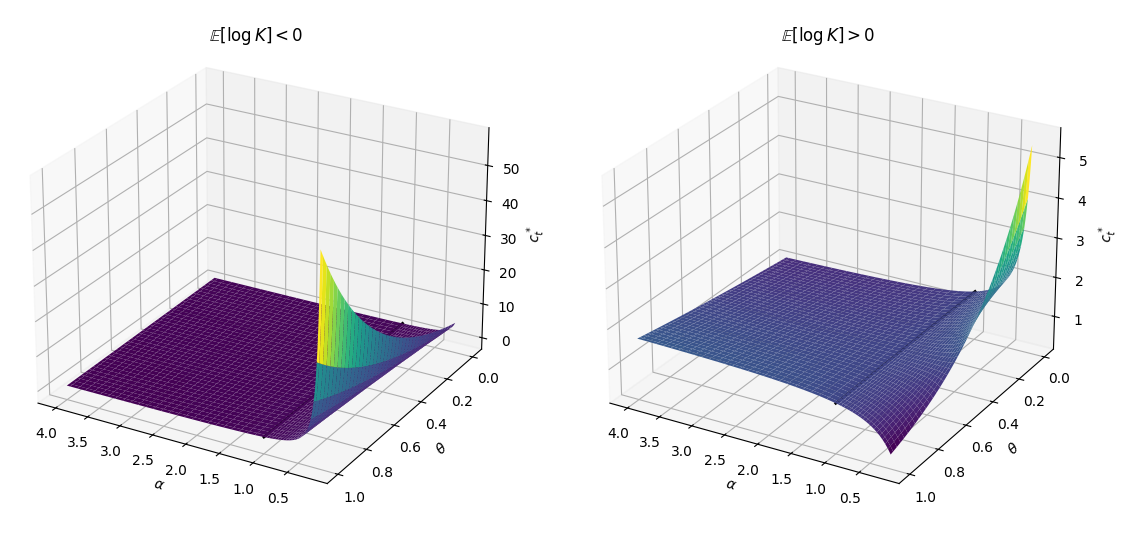}}
    \caption{Optimal consumption for $K=1.4$ and $\ABS{ \e \log(K)} = 0.5$}
    \label{fig:surfKpos}
\end{figure}

\begin{figure}[ht]
    \centering
    \makebox[\textwidth][c]{\includegraphics[scale=0.55]{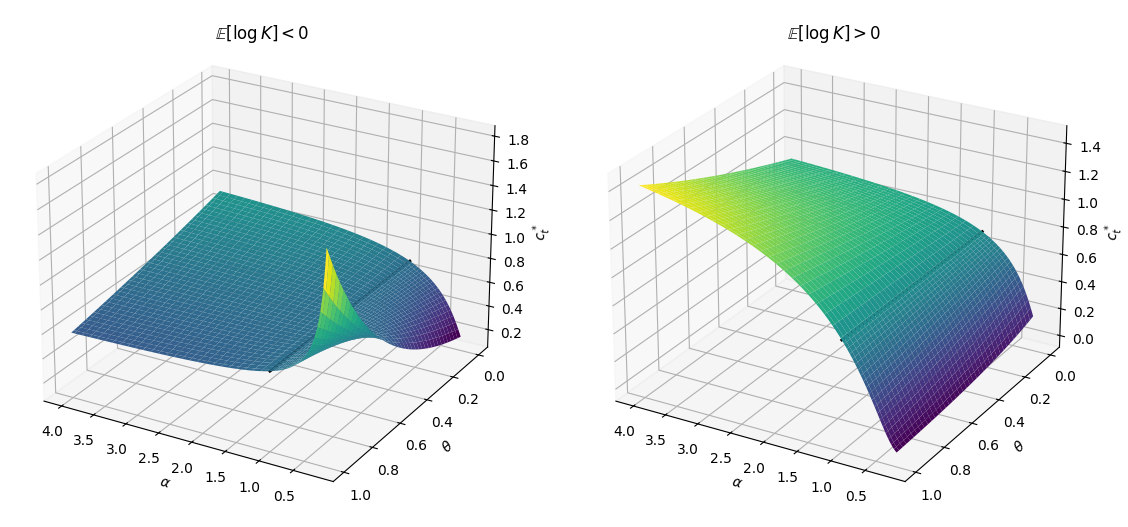}}
    \caption{Optimal consumption for $K=0.7$ and $\ABS{ \e \log(K)} = 0.5$}
    \label{fig:surkKneg}
\end{figure}

The benchmark case $\e \SBRA{\log(K)} = 0$ corresponds to a population that is, on average, neutral between consumption and wealth utility. It coincides with the situation with no competition, i.e $\theta = 0$, where the competitive correction vanishes so that the mean-field Nash equilibrium becomes $c_{t}^{*, NC} = \exp \PAR{ \frac{\log(K)}{\alpha}}.$ In other words, the equilibrium consumption \eqref{optconsMFproportionnel} is the optimal rate of consumption $c_{t}^{*, NC}$, corrected in the presence of competition by a term reflecting the population's aggregate preference. Since $c_{t}^{*, NC}$ depends on the sign of $\log(K)$, while the competitive correction depends on its average $\e \SBRA{\log(K)}$, we represent $c_{t}^{*}$ as a function of $(\alpha, \theta)$ in each sign configuration in Figure \ref{fig:surfKpos} and Figure \ref{fig:surkKneg}.

If the representative agent prefers consumption to wealth utility that is $\log(K) > 0$, then $c_{t}^{*, NC} \underset{\alpha \to 0}{\to} + \infty$, reflecting that risk seeking agents tends to consume more. By contrast, optimal consumption of risk averse agents converges to the proportional market consumption preference parameter $K$ as the risk aversion parameter goes to $1$. Conversely, if the representative agent places a higher weight on wealth than on consumption, i.e $\log(K) < 0$, then $c_{t}^{*, NC} \underset{\alpha \to 0}{\to} 0$, reflecting that the optimal consumption of risk seeking agents goes to $0$, while the limit as $\alpha$ goes to $1$ remains unchanged.

The magnitude of the correction is again determined by the quantity $K_{\alpha, \theta}$, while its direction is determined by the sign of the product $K_{\alpha, \theta} \e \log(K)$. In turn, the effect of competition on optimal consumption switches with the sign of both $\e \log(K)$ and $K_{\alpha, \theta}$. If $\e \log(K) < 0$ and $\alpha < 1$ (so that $K_{\alpha, \theta} > 0$), then $c_{t}^{*}$ increases with $\theta$. Conversely, if $\e \log(K) > 0$ and $\alpha > 1$ (so that $K_{\alpha, \theta} < 0$), then $c_{t}^{*}$ decreases with $\theta$. 
This emphasizes that the consumption equilibrium depends explicitly on the relative weight assigned to wealth versus consumption utilities, both at the individual level and at the population level.

\paragraph{Relative market-consumption preference parameter -} We now turn to the power type specification of the relative market-risk consumption preference, namely $\phi_{t} = Z_{t}^{1 - \kappa}$ as in Assumption \ref{ass:kappa}. When $\kappa \leq 0$ is common to all agents, the equilibrium consumption dynamics reduces to a logistic-type diffusion with explicit strong global solution, as detailed in Proposition \ref{prop:CRRA rel}. Long term behavior results are available when local characteristics $b^{A}$, $\delta^{A, B}$ and $\delta^{A, W}$ are time independent, see \cite{jiang2005note}, \cite{giet2015logistic}. We then introduce the parameter $q = - \frac{b^{A}}{\NRM{\delta^{A}}^{2}}$ which governs the asymptotic result stated below.

\begin{prop} \label{prop:cons_asympt}
Under Assumptions \ref{ass_MF_CRRA} and \ref{ass:kappa}, with $\kappa \leq 0$ common to the whole population, the mean-field Nash equilibrium is well defined, and the optimal consumption rate is explicitly given by
\begin{eqnarray} \label{eq:eq_cons_prop41}
c_{t}^{*} = \frac{c_{0}^{*}\exp \PAR{b^{A}t + \delta^{A}.d\overline{W_{t}}}}{1 - \kappa c_{0}^{*} \displaystyle \int_{0}^{t} \exp \PAR{b^{A}s + \delta^{A}.d\overline{W_{s}}}ds}, \quad \text{with} \quad c_{0}^{*} = \exp \PAR{\e \SBRA{\frac{1}{\alpha} \log \PAR{\frac{\phi_{0}}{z_{0}}}}}^{K_{\alpha, \theta}} \PAR{\frac{\phi_{0}}{z_{0}}}^{\frac{1}{\alpha}}
\end{eqnarray}
\begin{itemize}
\item If $q < 0$, $c_{t}^*$ converges in law towards a Gamma distribution with shape $-2q$ and scale $\frac{\NRM{\delta^{A}}^{2}}{2\ABS{\kappa}}$.
\item If $q > 0$, the diffusion goes to $0$ almost surely.
\end{itemize}
\end{prop}
Setting $\delta^{A} = 0$ recovers the deterministic equilibrium structure in Corollary 4.5 of \cite{dos2022forward}. Once the wealth utility volatility is non zero, the equilibrium consumption becomes stochastic, and its long run behavior no longer coincides with the deterministic benchmark.

\begin{figure}[ht]
\centering
\begin{minipage}[t]{0.48\textwidth}
        \includegraphics[scale=0.45]{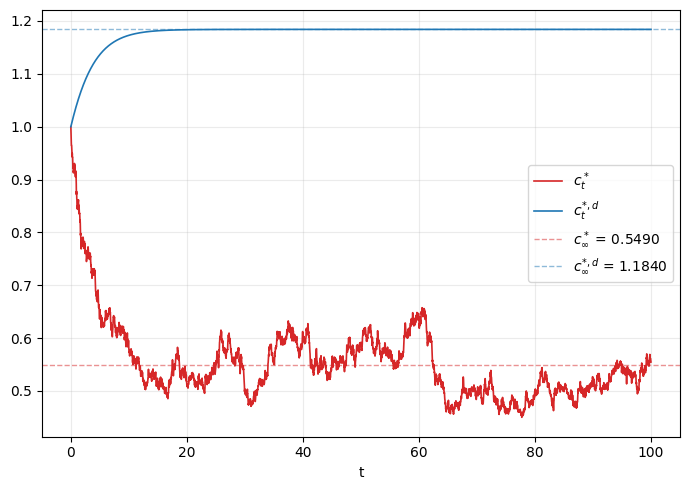}
        \caption{Example trajectory of the equilibrium consumption, $\alpha = 1.65$, $\theta = 0.65$, $\delta^{Z, B} = 0.25$}
    \label{fig:traj_cons1}
\end{minipage} \hfill
\begin{minipage}[t]{0.48\textwidth}
        \includegraphics[scale=0.45]{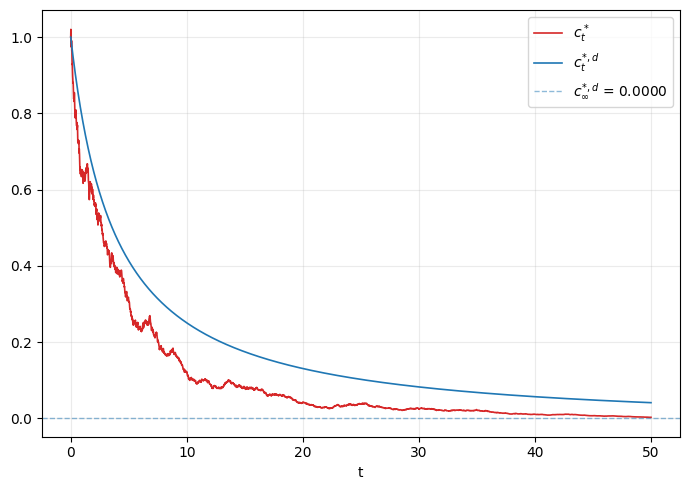}
        \caption{Example trajectory of the equilibrium consumption, $\alpha = 0.55$, $\theta = 0.35$, $\delta^{Z, B} = 0.25$}
    \label{fig:traj_cons2}
\end{minipage}
\end{figure}

\paragraph{Numerical setup -} To compare with \cite{dos2022forward}, we work in the single stock model with deterministic $\mu, \, \sigma$, $\nu=0$ and $\delta^{Z, W} = 0$. For numerical tests, we allow $\alpha > 1$ (extreme risk aversion), and sample $\alpha \sim \U(0, 3)$, $\alpha \neq 1$ and $\theta \sim \U(0, 1)$. For clarity, $\delta^{Z, B}$ is taken common across agents. Figures \ref{fig:traj_cons1} and \ref{fig:traj_cons2} show representative trajectories of $c_{t}^{*}$, together with the deterministic benchmark $c_{t}^{*, d}$ from \cite{dos2022forward}. We observe that these equilibria differ as soon as $\delta^{Z, B} \neq 0$.

\begin{figure}[ht]
    \centering
    \makebox[\textwidth][c]{\includegraphics[scale=0.4]{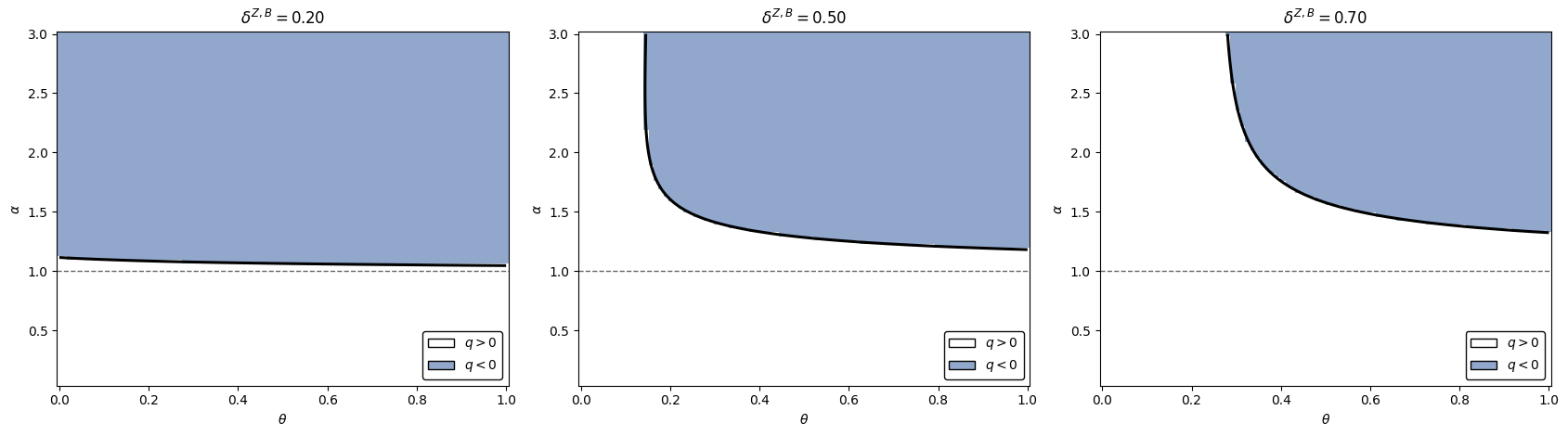}}
    \caption{Sign region for $q$}
    \label{fig:sign_region_q}
\end{figure}

\noindent
\textbf{Regime map -} From Proposition \ref{prop:cons_asympt}, the regime is governed by the sign of $q$, so equivalently by that of $b^{A}$. This coefficient compares $\overline{b^{Z}}/\alpha$ which represent the drift of wealth utility when there is no consumption, to its average among the population, with specific coefficients capturing competition and risk aversion preferences. When $\delta^{Z}$ is non zero, additional terms involving $\delta^{Z}$ and $\e \SBRA{\delta^{Z}}$ enter this comparison. Figure \ref{fig:sign_region_q} displays the corresponding sign regions for some volatility values.

\noindent
For small $\delta^{Z, B}$, the switching boundary is close to \(\alpha\approx 1\), with mild sensitivity to \(\theta\). As $\delta^{Z}$ increases, the $q=0$ boundary is reshaped and for a fixed $\alpha$, keeping $q > 0$ requires a higher competition preference parameter $\theta$.

\begin{figure}[ht]
    \centering
    \makebox[\textwidth][c]{\includegraphics[scale=0.6]{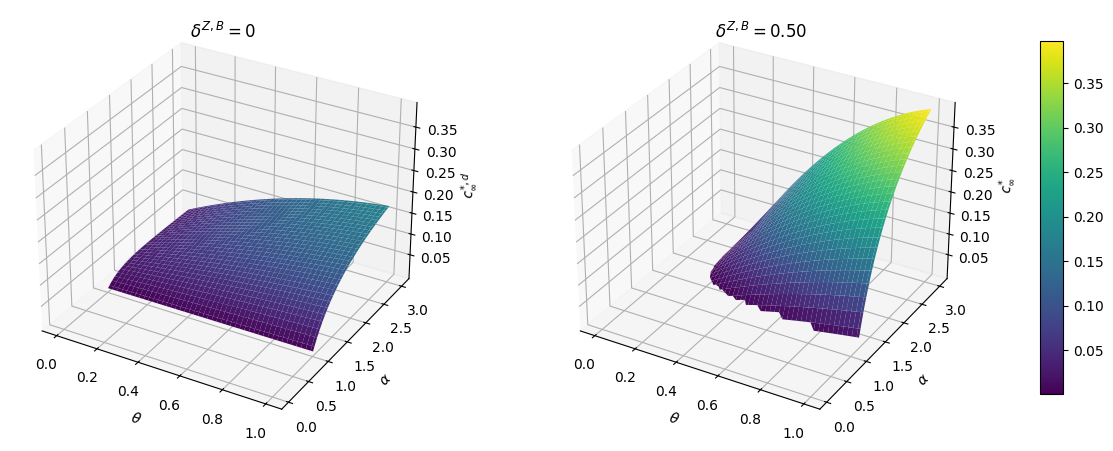}}
    \caption{Asymptotic consumption level function of $(\alpha, \theta)$}
    \label{fig:asymptcons}
\end{figure}

\paragraph{Asymptotic level -} When $q < 0$, Proposition \ref{prop:cons_asympt} yields that the equilibrium consumption converges in distribution to a Gamma distribution, whose expectation is given by 
\begin{eqnarray}
\e \SBRA{c_{\infty}^{*}} = - \frac{b^{A}}{\kappa}.
\end{eqnarray}
The asymptotic level of the deterministic benchmark $c_{t}^{*, d}$ in \cite{dos2022forward} has the same form, where $b^{A}$ is evaluated at $\delta^{Z} = 0$, denoted $b^{A}|_{\delta^{Z} = 0}$. Figure \ref{fig:asymptcons} compares asymptotic consumption levels as a functions of $\alpha, \, \theta$. Both $c_{\infty}^{*}$ and $c_{\infty}^{*, d}$ increase with $\alpha$ and $\theta$, but they differ substantially in magnitude and support.

\noindent
The difference is structural: the quantity $\overline{b^{Z}}$ depends on $\delta^{Z, B}$ both directly and through the equilibrium portfolio $\pi_{t}^{*}$. Therefore, introducing non-zero volatility modifies both the portfolio strategy level and the equilibrium consumption. In addition, opposite regimes can occur: there are parameter regions where $c_{t}^{*} \to 0$ while $c_{t}^{*, d}$ converges to a strictly positive level. This regime disagreement corresponds to the area below $q = 0$, and above $\alpha = 1$ in Figure \ref{fig:sign_region_q}.

\section*{Conclusion} 

Time consistency of forward relative performance criteria is characterized by an HJB-SPDE, with additional competition-related terms compared with the single agent forward framework of \cite{el2008backward}. In particular, this SPDE reveals a link between wealth and consumption utility, as discussed in Remark \ref{rmq:linkVU}. In the separable CRRA framework, and under an additional structural assumption on the drift of the consumption utility, we derive the best-response strategy in Theorem \ref{SPDEprop}, as well as the mean-field equilibrium in Theorem \ref{thm:MFSPDE}. Allowing for non-zero preference volatility provides an additional economically meaningful degree of freedom in equilibrium strategies. It rescales both the optimal portfolio, through the wealth utility volatility, and the optimal consumption, through the joint dynamics of wealth and consumption utilities. From an inverse preference viewpoint, echoing Remark \ref{rmq:prefdesign}, the wealth utility volatility can be used as a calibration parameter to match a prescribed target portfolio strategy. Theorem \ref{MFthm} then yields the associated equilibrium consumption, while ensuring that the resulting preference criterion remains time consistent. Our results contribute to a better understanding of the impact of non-vanishing utility volatility on the structure of forward preferences and the associated optimal policies. Extending our verification results to an exhaustive characterization of admissible forward preference criteria remains a challenging problem, as it requires a deeper understanding of the admissible class of solutions to the HJB-SPDE.

\paragraph*{Acknowledgements} The authors are grateful the Department of Mathematics at the National University of Singapore for their hospitality, as this work emerged from our visit to Professor Chao Zhou in May and June 2023. The authors also thank Professor Thaleia Zariphopoulou, as well as the anonymous referees, for their careful reading and valuable comments, which substantially improved the manuscript.

\newpage
\bibliographystyle{plain}
\bibliography{biblio}

\end{document}